\documentclass[12pt]{article}
\usepackage[utf8]{inputenc}
\usepackage[english]{babel}
\usepackage{amsmath,bm}
\usepackage{amsfonts}
\usepackage{amssymb}
\usepackage{graphicx}
\usepackage[left=2cm,right=2cm,top=2cm,bottom=2cm]{geometry}
\usepackage{setspace}
\usepackage{verbatim}
\usepackage{subfig}
\usepackage{xcolor}

\usepackage[ruled,vlined]{algorithm2e}

\usepackage{mathtools}	

\usepackage{authblk}

\usepackage{graphicx}
\usepackage{float}
\usepackage{subfig}

\makeatletter
\newbox\abstract@box
\renewenvironment{abstract}
{\global\setbox\abstract@box=\vbox\bgroup
	\hsize=\textwidth\linewidth=\textwidth
	\small
	\begin{center}%
		{\bfseries \abstractname\vspace{-.5em}\vspace{\z@}}%
	\end{center}%
	\quotation}
{\endquotation\egroup}
\expandafter\def\expandafter\@maketitle\expandafter{\@maketitle
	\ifvoid\abstract@box\else\unvbox\abstract@box\if@twocolumn\vskip1.5em\fi\fi}
\makeatother

\newcommand{\bluecom}[1]{{}}

\begin{document}
	
\title{An implicit FFT-based method for wave propagation in elastic heterogeneous media}

\author[1]{R. Sancho}
\author[1]{V. Rey de Pedraza}
\author[2,3]{P. Lafourcade}
\author[4]{R.A. Lebensohn \thanks{lebenso@lanl.gov}}
\author[1,5]{J. Segurado\thanks{javier.segurado@upm.es}}
\affil[1]{\small{Departamento de Ciencia de Materiales, ETSI Caminos, Canales y Puertos, Universidad Politécnica de Madrid, C/ Profesor Aranguren s/n, 28040 Madrid, España}}
\affil[2]{\small{CEA, DAM, DIF, F-91297 Arpajon, France}}
\affil[3]{\small{Universit\'e Paris-Saclay, CEA, LMCE, 91680 Bruy\`eres-le-Ch\^atel, France}}
\affil[4]{\small{Los Alamos National Laboratory, MS B216, Los Alamos, NM, 87544, USA}}
\affil[5]{\small{ Fundaci\'on IMDEA Materiales, C/ Eric Kandel 2, 28906, Getafe, Madrid, España}}
		
\renewcommand\Authands{ and }
	
\date{\small{\today}}
	
\providecommand{\keywords}[1]
{
\small	
\textbf{\textit{Keywords---}} #1
}
	
\begin{abstract}
An FFT-based algorithm is developed to simulate the propagation of elastic waves in heterogeneous $d$-dimensional rectangular shape domains. The method allows one to prescribe the displacement as a function of time in a subregion of the domain, emulating the application of Dirichlet boundary conditions on an outer face.
Time discretization is performed using an unconditionally stable beta-Newmark approach. The implicit problem for obtaining the displacement at each time step is solved by transforming the equilibrium equations into Fourier space and solving the corresponding linear system with a preconditioned Krylov solver. The resulting method is validated against analytical solutions and compared with implicit and explicit finite element simulations and with an explicit FFT approach. The accuracy of the method is similar to or better than that of finite elements, and the numerical performance is clearly superior, allowing the use of much larger models. To illustrate the capabilities of the method, some numerical examples are presented, including the propagation of planar, circular, and spherical waves and the simulation of the propagation of a pulse in a polycrystalline medium.
\vspace{1cm}
\end{abstract}

\maketitle
			
\keywords{FFT-based homogenization, micromechanics, wave propagation, polycrystals, acoustics, elastodynamics, scattering}
\footnote{Article accepted in Computer Methods in Applied Mechanics and Engineering}
\clearpage
\section*{Notation}

\begin{tabular}{ll}
Vector and tensor notation\\
$\mathbf{x},\boldsymbol{\xi},\mathbf{u}$  & Vectors $x_{i},\xi_{i},q_{i}$\\
$\boldsymbol{\sigma,\epsilon},\mathbf{g}$ & Second-order tensors $\sigma_{ij},\epsilon_{ij},g_{ij}$\\
$\mathbb{C},\mathbb{G},\mathbb{K}$ & Fourth-order tensors $C_{ijkl},G_{ijkl},K_{ijkl}$\\
$\mathbf{A}=\boldsymbol{\tau}\cdot \mathbf{F}$ & Dot product $A_{ij}=\tau_{ip}F_{pj}$\\
$e=\boldsymbol{\sigma}:\boldsymbol{\varepsilon}$ & Double dot product $e=\sigma_{ij}\varepsilon_{ij}$ \\
$\mathbf{I}$ & Second-order identity tensor $I_{ij}=\delta_{ij}$ \\ \\
Differential operators\\
$v(x)=\mathcal{A}(u(x))$ & Linear differential operator \\
$\mathbf{U}=\nabla \mathbf{u}$ & Gradient of a vector field $U_{ij}=\frac{\partial u_i}{\partial X_j}$ \\
$\boldsymbol{\varepsilon}=\nabla^s \mathbf{u}$ & Symmetric gradient of a vector field $\varepsilon_{ij}=\frac{1}{2} \left( \frac{\partial u_i}{\partial X_j}+\frac{\partial u_j}{\partial X_i}\right)$ \\
$\mathbf{F}=\nabla\cdot \boldsymbol{\sigma}$ & Divergence of tensor field $F_i=\frac{\partial \sigma_{ij}}{\partial X_j}$ \\\\
Fourier Transforms and convolutions\\
$\hat{f} = \mathcal{F}(f)$ & Fourier transform of $f$\\
$ f = \mathcal{F}^{-1} (\hat{f})$ & Inverse Fourier transform of $\hat{f}$\\
$G \ast P$ & Convolution operation \\
\end{tabular}
\clearpage

\section{Introduction}

The effect of the microstructure on the mechanical response of heterogeneous materials has been thoroughly studied from both the experimental and modeling point of view. In the case of quasistatic behavior, it is well-known that the microstructure affects the stiffness, yield, and fracture. In the case of dynamic excitation, the microstructure also plays a fundamental role. For low energies, the microstructure has a strong effect on the acoustic response of the material, as has been reported in many studies in polycrystals \cite{Bathia1959},  composites \cite{tanaka2000band,Mariatti2016}, porous \cite{attenborough1982acoustical}, or architectured materials \cite{fleck2010}. For higher energies, the microstructure also controls the nonlinear response and failure of these materials \cite{CHEN2012218,lieberman2016microstructural,CLAYTON20054613}.

Modeling the dependency of the mechanical response with the microstructure at the mesoscopic scale relies on continuum micromechanics, in which the domain under study incorporates the microstructure but the continuum hypothesis still applies. In the quasistatic regime, computational homogenization \cite{Milton2003,Segurado2018} and many multiscale models \cite{matouvs2017review} are based on numerical simulation of the mechanical response of representative volume elements (RVE) of the microstructure, usually under periodic boundary conditions. For dynamical problems, computational approaches based on representative unit cells can still be used in simulations performed in the frequency domain, using the Bloch wave formalism \cite{kushwaha1993,Hussein2009,Segurado2021}. However, in many other cases, such as the study of the propagation of short waves or impulses, the simulations require a dynamic time-domain solution. Moreover, in these cases the study of periodic cells is no longer valid, and full domains that include the microstructure are used \cite{vanpamel2017,Huang2020} or multiscale simulations are considered \cite{SRIDHAR2018414}.

The resolution of dynamical problems in domains which explicitly represnt the microstructure is computationally very demanding and requires the use of efficient and stable numerical approaches. Moreover, wave propagation can lead to spurious oscillations in the vicinity of wave fronts, and various studies have been dedicated to the development of efficient time-integration schemes that help to  reduce these oscillations \cite{Bathe2007, Noh2013, Park2012}. Implicit time integration schemes are potentially the most appropriate choice because they can be unconditionally stable allowing large time increments. Due to their interest, implicit schemes have been recently introduced to study wave propagation with numerical schemes alternative to standard FE such as isogeometrical methods \cite{ZAMPIERI2021114047} or discontinuous Galerkin \cite{Kronbichler2016}. However, implicit solvers require to solve a full mechanical problem at each time step and their computational cost in general very expensive. Therefore, numerical methods such as finite elements (FE) or finite differences tend to overcome this issue using an explicit integration, and almost every attempt to study wave propagation and other types of dynamical problems explicitly considering the microstructure is based on explicit integration schemes. In fact, many researchers have made an important effort to develop new and more efficient explicit-based approaches using standard FE \cite{ZHANG2021113811, AURICCHIO20122}, discontinuous Galerkin \cite{Kronbichler2016,STANGLMEIER2016748}, or isogeometric methods \cite{ZAMPIERI2019776}. However, explicit FE methods are conditionally stable, so the time step is still limited in terms of the highest frequency of the system. Therefore, the maximum time step is inversely proportional to the discretization level. In the particular case of full-field micromechanical models, the requirement of fine meshes to represent the complex microstructures imposes extremely small time steps, so the resulting models are still very expensive. An example of this extreme computational demand is the study of wave propagation and attenuation in polycrystals, massively parallelized in large computer clusters \cite{vanpamel2017,Huang2020}.  

In computational homogenization for the quasistatic regime, FFT-based homogenization approaches are nowadays a mature and very extended alternative to FE models, see the articles \cite{Lucarini2021,schneider2021} for a comprehensive review. The main benefit of this technique, introduced for micromechanics in the seminal work of Moulinec and Suquet \cite{MS94}, is its excellent numerical performance, which allows one to simulate very complex microstructures with fine discretizations in a fraction of the time needed for FE models. Furthermore, the FFT methods are based on voxelized grids, which allows the construction of RVEs in a very simple manner synthetically or using tomographic data as input \cite{lucarini2022adaptation}. The method has been successfully applied to a wide range of micromechanical quasistatic problems including non-linear homogenization \cite{DEGEUS2017412}, phase field fracture \cite{CHEN2019167,ERNESTI2020112793}, multiscale modeling \cite{NKOUMBOUKAPTCHOUANG2022114921,GIERDEN2021113566}, etc.

However, the use of FFT based methods for dynamical problems is very scarce. One of the first and few approaches that use spectral solvers for elastodynamic problems and that analyze the advantage of implicit over explicit time-integration schemes is the work presented by Zampieri and Tagliani \cite{zampieri1997numerical}.  Further studies on wave propagation methods relying on spectral approaches have been performed by the same group, including the use of spectral finite elements with explicit  \cite{ZAMPIERI2006308,ZAMPIERI2006b} or implicit  \cite{ZAMPIERI20062649} integration methods. All these approaches have been developed for the wave equation with constant coefficients to represent a homogeneous medium. An alternative direction within spectral approaches was proposed by Amlani and Bruno \cite{AMLANI2016333}, who developed an elastodynamic solver based on the Fourier Continuation method (\cite{bruno2010high}). The proposed method is capable of considering nonperiodic domains and boundary conditions and, due to its spectral nature, is free from dispersion errors. However, the method proposed is explicit (is implemented using a fourth-order Adams-Bashforth method), and therefore the time step for stable solutions decreases linearly with the number of discretization points, still being very expensive for coarsediscretizations.  Regarding the extension of FFT-based homogenization to elastodynamic problems in heterogeneous microstructures, the only work available in the literature (to the authors knowledge) was a very interesting work proposed recently by Morin et. al. \cite{Morin2021}. The focus of this work was the study of fracture propagation in heterogeneous media including dynamical effects. In this work the load was prescribed by macroscopic homogeneous strain/stress history. However, the approach proposed relied on explicit dynamics, with the computational cost associated for fine discretizations. Moreover, in \cite{Morin2021} as well as in almost every computational homogenization study, standard periodic boundary conditions are imposed, discarding the introduction of a local perturbation to study wave propagation. The only exception to periodicity, to the authors knowledge, is the introduction of Dirichlet boundary conditions in  elastic problems proposed in \cite{Gelebart2020}.

In summary, an FFT-based method for elastodynamic problems in heterogeneous materials which allows to prescribe a local excitation and which relies on implicit integration to eliminate the inverse dependency of the time increment size with discretization is still missing.

The main objective of this work is to develop a novel and efficient approach for wave propagation in heterogeneous media based on implicit time integration and FFT solvers. The method will allow us to prescribe the displacement as a function of time on a subregion, emulating time-dependent Dirichlet boundary conditions, and allowing us to solve propagation of pulses or impacts in a heterogeneous solid. The method will be applied for studying the propagation of waves in polycrystals to illustrate the potential of the technique proposed for these studies.

The article is organized as follows. Section 2 presents the method and the associated algorithms. In Section 3, the accuracy and numerical efficiency of the method will be evaluated against analytical solutions and finite element simulations. Section 4 presents some numerical examples in heterogeneous media, and finally, conclusions and open issues will be given in Section 5.

\section{Theory and numerical approach}

\subsection{One dimensional case}\label{subsec:one_dim}

\begin{figure}[h]
    \centering
    \includegraphics[width = 0.7\textwidth]{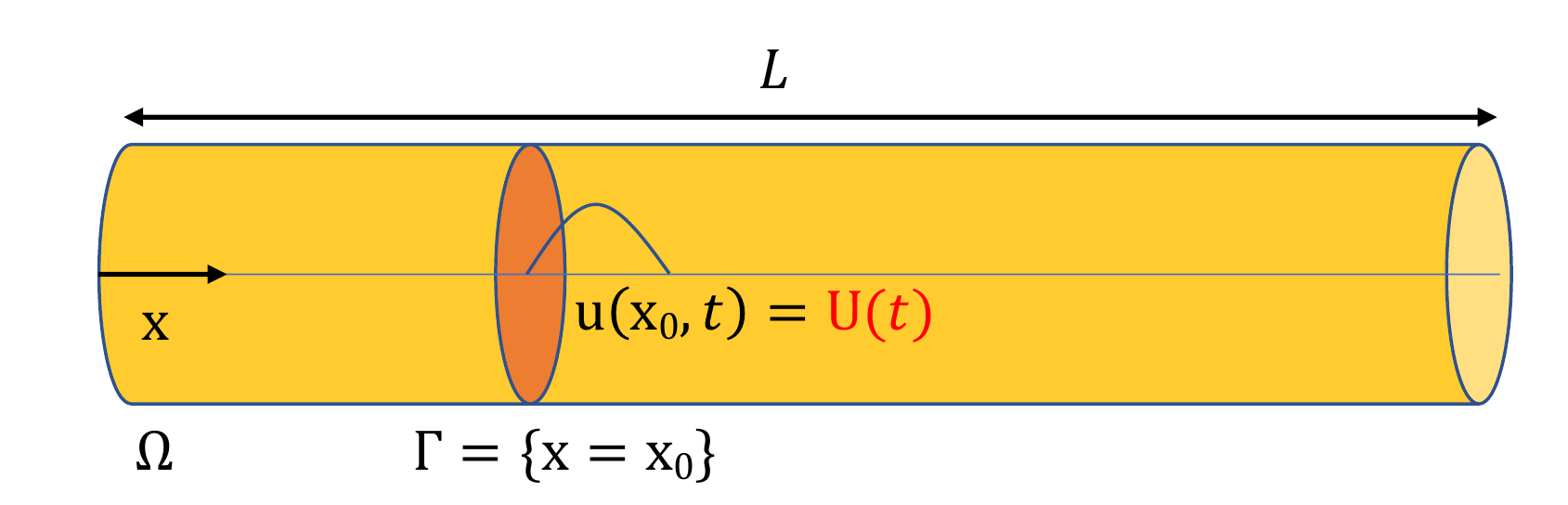}
    \caption{1D-problem definition.}
    \label{fig:1D-problem}
\end{figure}

The objective is to solve the propagation of elastic waves in a one-dimensional domain $\Omega=\{0\leq x \leq L\}$ occupied by a heterogeneous linear elastic medium. The wave is introduced in the domain by perturbing a planar region $\Gamma$, perpendicular to the propagation direction, and which in the one dimensional representation corresponds to a single point, here named $x_0$ (Figure \ref{fig:1D-problem}). The perturbation consists in a prescribed displacement $U(t)$ , so that $u(x\in \Gamma, t) =u(x_0, t)= U(t)$. Longitudinal waves are considered to particularize the problem, but the equations can be directly adapted for shear waves.
The heterogeneous elastic domain is defined by a spatial distribution of Young's modulus $E(x)$ and density $\rho(x)$. Periodic boundary conditions apply for all fields involved, as a requirement of the Fourier-based approach. 

The prescription of a displacement at the region $\Gamma$ is made by introducing an artificial body force field whose value is set to meet the displacement condition. Let $f(x,t)$ be a singular body force density field concentrated around $\Gamma$ (which corresponds to the point $x_0$). The field $f(x,t)$ can then be represented using a Dirac delta function
\begin{equation}\label{eq:f_in1D}
f(x,t)=F(t)\delta_\Gamma(x) = F(t) \delta(x-x_0)
\end{equation}
such that 
$$\int_{0}^L f(x,t)  \mathrm{d}x = \int_{0}^L F(t)\delta(x-x_0)\mathrm{d}x=F(t).$$ 

The value of the total force prescribed at every time, $F(t)$, is obtained to fulfill the displacement condition at that point. Within this framework, the linear momentum balance for the problem can be written as
\begin{equation}
\label{eq:1D_consv_lm}
\left\{ 
\begin{array}{c}
\frac{\partial}{\partial x} \left[\sigma(x,t) \right] + f(x,t)  =  \rho(x) \ddot{u}(x,t) \\
u(x_0,t)=U(t) \quad;\quad u(x,t),\sigma(x,t) \text{periodic} \\
u(x,0) = u_0 (x) \quad;\quad \dot{u}(x,0) = \dot{u}_0 (x)
\end{array} \right.
\end{equation}
where $u_0 (x)$ and $\dot{u}_0 (x)$ correspond to the initial value of the displacement and velocity fields respectively and the body force field $f(x,t)$ is given by Eq. \eqref{eq:f_in1D}.

The constitutive equation is the one-dimensional Hooke's law, which defines the relation between the stress $\sigma$ and the strain $\varepsilon$ at every point of the domain,
\begin{equation}
\label{eq:1D_const-eq}
\sigma(x,t)=E(x)\varepsilon(x,t)= E(x)\frac{\partial}{\partial x} u(x,t).
\end{equation}

The unknowns of the problem for a time interval $[0,T]$ are the field $u(x,t) :\Omega\times[0,T]\rightarrow \mathbb{R}$ and the function $F(t):  [0,T] \rightarrow \mathbb{R}$.

\subsubsection{Time discretization and integration scheme}
The time interval $[0,T]$ is discretized in $N$  time increments of size $\Delta t=T/N$, such that

\begin{equation}
\label{eq:time_discrt}
\left\{ 
\begin{array}{c}
t_n = n \cdot \Delta t \\
u(x,t_n)=u_n(x) \quad;  \dot{u}(x,t_n)=\dot{u}_n (x) \quad; \ddot{u}(x,t_n)=\ddot{u}_n (x)
\end{array} \right.
\end{equation}
Eq. \eqref{eq:1D_consv_lm} has to be integrated in time, and
the implicit Newmark integration scheme is used for this purpose. The velocity $\dot{u}$ and the acceleration $\ddot{u}$ of the time step $n$ are expressed as
\begin{equation}
\label{eq:Newmark}
\begin{split}
	&\dot{u}_{n}=\dot{u}_{n-1} + \Delta t (1-\gamma)\ddot{u}_{n-1} + \Delta t \gamma \ddot{u}_{n} \\
	&\ddot{u}_{n}=\frac{1}{\beta \Delta t^2} \left( u_{n}-u_{n-1}- \Delta t \dot{u}_{n-1} - \Delta t ^2 (0.5-\beta)\ddot{u}_{n-1} \right)
\end{split}
\end{equation}
where $\beta$ and $\gamma$ are two non-negative real numbers that define the version of the Newmark-$\beta$ method. The values $\beta = 0.25$ and $\gamma = 0.5$ are used to define the implicit  unconditionally stable integration used.

\subsubsection{Continuum solution}

Introducing the acceleration $\ddot{u}_n$ for a time $t_n$, from Eq. \eqref{eq:Newmark} in the conservation of linear momentum (Eq. \eqref{eq:1D_consv_lm}) leads to

\begin{equation}
\label{eq:1D_consv_lm_Newmark}
\frac{\mathrm{d}}{\mathrm{d}x}\left[E(x)\frac{\mathrm{d}}{\mathrm{d}x}u_{n}\right]+f_n(x) = 
\frac{\rho(x)}{\beta \Delta t^2} \left[u_{n}-u_{n-1}- \Delta t\dot{u}_{n-1} - \Delta t^2 (0.5-\beta)\ddot{u}_{n-1} \right],
\end{equation}
where $f_n$ stands for the force density field prescribed at time $t_n$ and, for brevity, $u_n (x)$ has been written as $u_n$. Grouping the terms (unknown left, known right) and multiplying by the factor $- \beta \Delta t^2$ allows us to write the equation as 
\begin{equation}
\label{eq:1D_consv_lm_Newmark_ord}
-\beta \Delta t^2 \frac{\mathrm{d}}{\mathrm{d}x}[E(x)\frac{\mathrm{d}}{\mathrm{d}x}u_{n}] + \rho(x) u_{n} = \beta \Delta t^2 f_n(x) + \rho(x) \left[ u_{n-1} + \Delta t \dot{u}_{n-1} + \Delta t^2 (0.5-\beta) \ddot{u}_{n-1} \right].
\end{equation}

Equation \eqref{eq:1D_consv_lm_Newmark_ord} is a differential equation of the displacement at time $t_n$, and will be solved using Fourier transforms. To this aim, the derivatives therein have been computed using the definition of the derivatives in Fourier space of a scalar field $V(x)$:
\begin{equation}
    \label{eq:Fourier_derivatives}
    \mathcal{F}\lbrace V'(x)\rbrace = i \xi \mathcal{F}\lbrace V(x)\rbrace = i \xi \hat{V}(\xi)
\end{equation}

\noindent where $\mathcal{F}$ corresponds to the Fourier transform, $\mathcal{F}\lbrace V(x) \rbrace=\hat{V}$, $i$ is the imaginary unit and $\xi$ is the spatial frequency. Equation \eqref{eq:1D_consv_lm_Newmark_ord} in Fourier space reads

\begin{equation}
\label{eq:Fourier_1D_lm}
\begin{split}
	&-\beta \Delta t^2 i\xi \mathcal{F} \left\lbrace E(x) \mathcal{F}^{-1} \left\lbrace i \xi \hat{u}_n \right\rbrace \right\rbrace + \mathcal{F} \left\lbrace \rho(x) \mathcal{F}^{-1} \{\hat{u}_n\} \right\rbrace  = \\
&	\mathcal{F} \left\lbrace \rho(x) \left(u_{n-1}+\Delta t \dot{u}_{n-1} + \Delta t ^2 (0.5-\beta)\ddot{u}_{n-1} \right)\right\rbrace + \beta \Delta t^2 \mathcal{F}(f_n) ,
\end{split}
\end{equation}

\noindent where the Fourier transform of the force density field has an explicit form as
$$\mathcal{F}(f_n)=\hat{f}_n=F_n e^{-i\xi x_0}.$$
The left-hand side of Eq. \eqref{eq:Fourier_1D_lm} corresponds to a linear operator $\hat{\mathcal{A}}(\cdot)$ in Fourier space that acts on the displacement field defined in Fourier space.

\begin{equation}
\label{eq:A-operator}
\hat{\mathcal{A}}(\hat{u}_n) = - \beta \Delta t^2 i\xi \mathcal{F} \left\lbrace E(x) \mathcal{F}^{-1} \left\lbrace i \xi \hat{u}_n \right\rbrace \right\rbrace+ \mathcal{F} \left\lbrace \rho(x) \mathcal{F}^{-1} \hat{u}_n \right\rbrace
\end{equation}

On the right-hand side, let $\hat{b}_{n-1}$ be a term that groups the Fourier transform of all fields from the previous step,
\begin{equation}
\label{eq:b_Fourier-term}
\hat{b}_{n-1} = \mathcal{F} \left\lbrace \rho(x) \left(u_{n-1}+\Delta t \dot{u}_{n-1} + \Delta t ^2 (0.5-\beta)\ddot{u}_{n-1} \right)\right\rbrace.
\end{equation}
Therefore, Eq. \eqref{eq:Fourier_1D_lm} can be condensed as;
\begin{equation}
\label{eq:Fourier_1D_lm_compact}
\hat{\mathcal{A}}(\hat{u}_{n}) = \beta \Delta t^2 \hat{f}_n +\hat{b}_{n-1} 
\end{equation}

The displacement field $\hat{u}_n$ in Fourier space is computed by solving two linear algebraic problems. The first provides $\hat{u}_{b}$,
\begin{equation}\label{eq:lin1D_1}
\hat{\mathcal{A}}(\hat{u}_{b})=\hat{b}_{n-1},
\end{equation}
and the second one provides the effect of a force $F$ on $x_0$, $\hat{u}_f$, 
\begin{equation}\label{eq:lin1D_2}
\hat{\mathcal{A}}(\hat{u}_f)= \hat{f}_n.
\end{equation}
Note that the solution of this problem provides the Green's function of the operator $\mathcal{A}(\cdot)$ for a force $F_n$ on the point $x_0$, $g(x,x_0)$ which is defined as 
\begin{equation}\label{eq:Greens-function}
u_f(x)= g(x,x_0) F_n
\end{equation}
In the case of a homogeneous medium, the linear operator $\mathcal{A}$ defines a Helmholtz equation with constant coefficients, and its Green's function for an infinite medium has an analytical expression in real space, which is shown in appendix \ref{anex:green}. For a heterogeneous medium, the Green's function $g(x,x_0)$ can be obtained numerically and stored once at the beginning of the simulation by solving this linear problem.

The solution of the displacement can then be expressed as
\begin{equation}
\label{eq:Sol_Fourier}
\hat{u}_n = \beta \Delta t^2 F_n \hat{u}_f +\hat{u}_{b} 
\end{equation}
which in real space corresponds to
\begin{equation}
\label{eq:Sol_Real}
u_n =  \beta \Delta t^2 F_n \mathcal{F}^{-1} \lbrace \hat{u}_f \rbrace + u_{b} =  \beta \Delta t^2 F_n g(x,x_0) + \mathcal{F}^{-1} \lbrace \hat{u}_{b} \rbrace 
\end{equation}

Introducing in Eq. \eqref{eq:Sol_Real} the prescribed displacement at $x_0$, $u_n(x=x_0) = U(t_n)$ provides an equation to obtain the force $F_n$,
\begin{equation}
\label{eq:L-determination}
F_n = \dfrac{1}{\beta \Delta t^2} g^{-1}(x_0,x_0) \left(U(t_n) - u_{b}(x_0)\right)
\end{equation}
where $g^{-1}(x_0,x_0)$ is the inverse of the Green's function on $x_0$, which represents the unit force to be applied at a point to produce a unit displacement. It is important to note that $g(x,x)$ is not singular at $x_0$, as it is shown in the appendix \ref{anex:green}.

Finally, once $F_n$ is determined by the previous equation 
, the displacement field $u_n(x)$ is obtained from Eq. \eqref{eq:Sol_Real}.

\subsubsection{Spatial discretization and solution algorithm}
In order to numerically solve the problem, the domain $\Omega$ will be divided into $N$ equal segments, and the values of the functions will be represented by their value in the center of each segment,
\begin{eqnarray}
E(x)  \rightarrow E(x_m)=E^m,\nonumber \\
\rho(x)\rightarrow \rho(x_m)= \rho^m \nonumber \\
u(t=t_n,x=x_m) \rightarrow u_n^m,
\end{eqnarray}
with $x_m=\frac{L}{N}(\frac{1}{2} + m)$ and $m\in[0,N-1]$. 
The point force defined using a delta function, $F_n=\int_\Omega f_n \delta(x-x_0) \mathrm{d}\Omega$, is introduced in the discrete version as
$$f_n(x^m) = \left\{ \begin{array}{c}
F_n \ \text{if} \ x^m=x_0 \\ 
0 \ \text{if} \ x^m\neq x_0 
\end{array}
\right.
$$

The Fourier transform and its inverse will be approximated with the discrete Fourier transform that can be computed with the efficient FFT algorithm. The corresponding $N$ discrete frequencies in Fourier space are 
\begin{equation}\label{frequencies}
\xi_k=\left\{ \begin{array}{c} 
\frac{2\pi}{L} (k-\frac{(N-1)}{2}) \quad \text{ for $N$ odd}\\
\frac{2\pi}{L} (k-\frac{N}{2}) \quad \text{ for $N$ even}
\end{array}
\right.  \text{ for } \quad k =0, \dots, N-1 
\end{equation}

The discrete version of the problem defined in Eqs. 10, 12 and 14 will be solved numerically using Krylov linear solvers. In particular, the linear equations will be solved in real space using the conjugate gradient. The resulting algorithm is given in Algorithm \ref{alg1D}.

\begin{algorithm}[H]\label{alg1D}
\SetAlgoLined
\caption{FFT-based elastodynamic algorithm --1D--}\label{alg:1D_basic_implementation}
\KwData{Microstructure: $E(x), \rho(x)$,\\
Prescribed fields on $\Gamma$: $U(t, x_0), \dot{U}(t, x_0), \ddot{U}(t, x_0)$, \\
Initial conditions: $u_0(x), \dot{u}_0(x)$;
Time discretization [$t_0, n \Delta t, ..., t_f $] and space discretization $x = [0, x_m, ..., L]$, for $n$ and $m$ $\in \mathbb{Z}$; \\
Newmark constants, $\beta=0.25$, $\gamma=0.5$}
\KwResult{$u(t, x), \dot{u}(t, x), \ddot{u}(t, x)$, $F(t)$}	
Form function $\hat{\mathcal{A}}()$ eq. \eqref{eq:A-operator} and compute $\ddot{u}_0(x)=\mathcal{F}^{-1}\left( i\xi \mathcal{F} \left\lbrace E(x) \mathcal{F}^{-1} \left\lbrace i \xi \hat{u}_0 \right\rbrace \right\rbrace \right)$;\\
$g(x, x_0) = \mathcal{F}^{-1}\left\lbrace \hat{\mathcal{A}}^{-1}(e^{-i\xi x_0}) \right\rbrace$ ; $\hat{\mathcal{A}}^{-1}(\cdot)$ using CG ($\hat{\mathcal{A}}, \hat{\delta}(x=x_0)$) \;
\While{$t_n \leq t_f$}{
$t_n = t_{n-1} + \Delta t $\;
$\hat{b}_{n} = \mathcal{F}\left\lbrace\rho(x) \left(u_{n-1} (x)+ \Delta t \dot{u}_{n-1}(x) + \Delta t^2 (0.5-\beta) \ddot{u}_{n-1}(x) \right)\right\rbrace$\;
$\hat{u}_{b_n} = \mathcal{F}^{-1}\left\lbrace \hat{\mathcal{A}}^{-1}(\hat{b}_n) \right\rbrace$; $\hat{\mathcal{A}}^{-1}(\cdot)$ using CG ($\hat{\mathcal{A}}, \hat{b}_n$)\;
$F_n = \dfrac{U(t_n)-u_{b_n}(x_0)}{\beta \Delta t^2 g (x_0, x_0)}$\;
$u_n(x) = \beta \Delta t^2 F_n g (x, x_0) + u_{b_n}(x)$\;
$\ddot{u}_n(x) = \frac{1}{\beta \Delta t^2} \left( u_{n}(x)-u_{n-1}(x)- \Delta t \ \dot{u}_{n-1} (x)- \Delta t ^2 (0.5-\beta)\ddot{u}_{n-1}(x) \right)$\;
$\dot{u}_n(x) = \dot{u}_{n-1}(x) + \Delta t (1-\gamma)\ddot{u}_{n-1} (x)+ \Delta t \ \gamma \ddot{u}_{n}(x)$\;
}
\end{algorithm}

\vspace{3mm}

{\bf Solution for a homogeneous material:} In the case of a homogeneous material, Young's modulus and density are constant along the space material $E(x) = E$ and $\rho(x)=\rho$. In this case, $\hat{\mathcal{A}}^{-1}$ has a closed-form expression, 
\begin{equation}
\label{eq:A_oper_inv_homog}
\hat{\mathcal{A}}^{-1}(\cdot) = \dfrac{1}{(\beta \Delta t^2 E \xi^2 + \rho)} (\cdot),
\end{equation}
and there is no need to use CG to solve Eqs. \eqref{eq:lin1D_1} and \eqref{eq:lin1D_2}.
\clearpage

\subsection{Multidimensional case}
\label{nD}
\subsubsection{Problem statement}

\begin{figure}[H]
    \centering
    \includegraphics[width=0.6\textwidth]{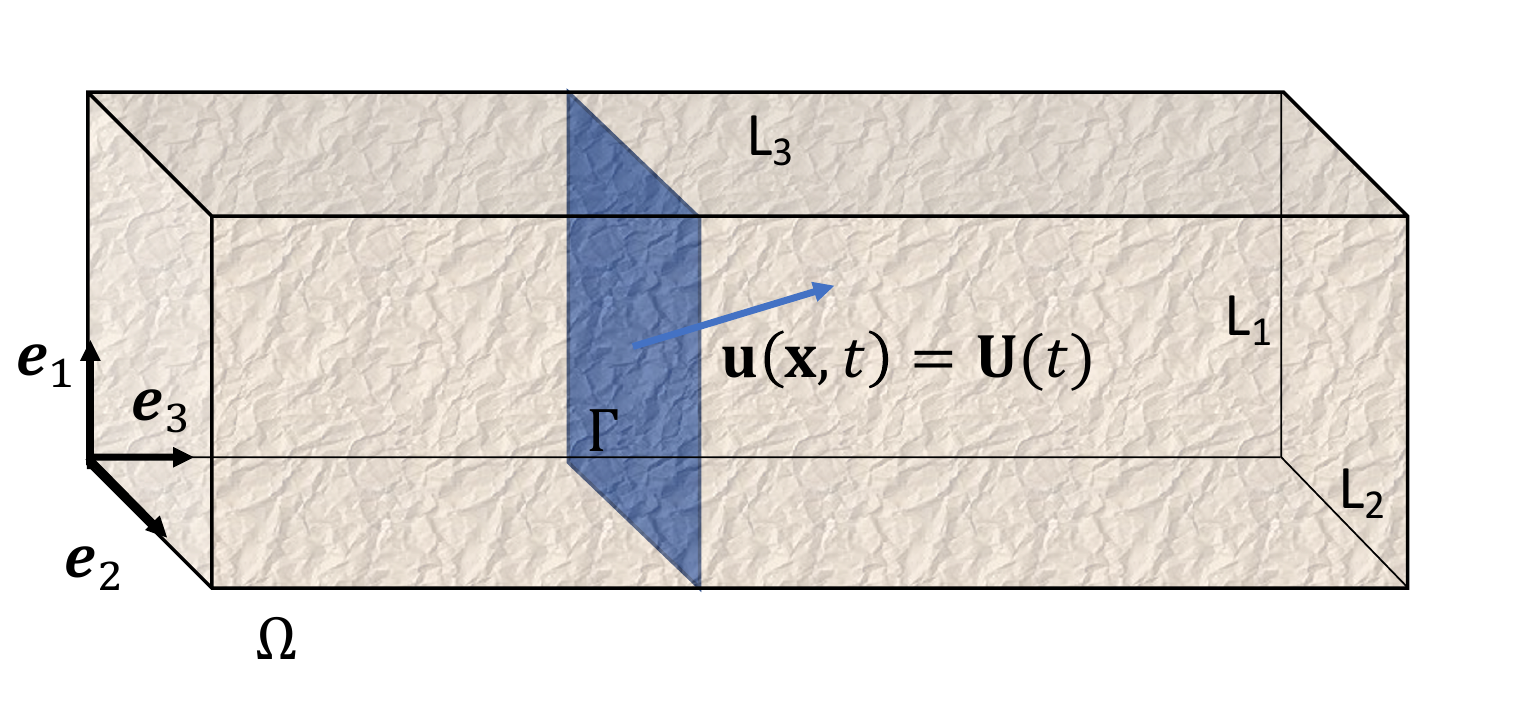}
    \caption{3D periodic domain $\Omega$ and the plane $\Gamma$ in which displacement is prescribed as $\mathbf{U}(t)$}
    \label{Fig:3Ddomain}
\end{figure}

Let $\Omega$ be a periodic prismatic domain in $d=2$ or $d=3$ dimensions, $\Omega= \left\{\mathbf{x} \ |\  0\leq x_i\leq L_i ; 1\leq i \leq d \right\}$, Fig. \ref{Fig:3Ddomain}. The domain is occupied by a heterogeneous linear elastic medium, whose microstructure is characterized by the spatial distribution of the stiffness tensor and density $\mathbb{C}(\mathbf{x})$, $\rho(\mathbf{x})$, respectively.  
The problem consists in finding the periodic displacement field for a time interval $[0, T]$, $\mathbf{u}(\mathbf{x},t):\Omega \ \mathrm{x} \ [0, T]:\rightarrow \mathbb{R}^d$ such that the problem is in equilibrium, and the value of $\mathbf{u}(\mathbf{x},t)=\mathbf{U}(t)$ is prescribed in a lower-dimensional embedded manifold $\Gamma$ in $\mathbf{R}^d$, a curve if $d=2$ and a surface if $d=3$. To impose the displacement value at $\Gamma$, a force field $\mathbf{F}(\mathbf{x},t)$ that acts only on $\mathbf{x}\in \Gamma$ is introduced and its value is calculated at each time to satisfy the displacement condition on $\Gamma$. Note that since $\mathbf{u}$ is set as periodic, the macroscopic strain $\mathbf{E}=<\nabla^s \mathbf{u}>$ becomes equal to zero.

The starting point for solving the problem is the conservation of the linear momentum.
\begin{equation}
\label{eq:ND_consv_lm}
\left\{ 
\begin{array}{c}
\nabla \cdot \left[\boldsymbol{\sigma} (\mathbf{(x}, t)\right] + \mathbf{f}(\mathbf{x},t)  = \rho(\mathbf{x}) \ddot{\mathbf{u}}(\mathbf{x},t) \\
\mathbf{u}(\mathbf{x}\in \Gamma,t)=\mathbf{U}(t) \quad; \quad \mathbf{u}(\mathbf{x},t), \boldsymbol{\sigma} (\mathbf{x}, t) \ \mathrm{periodic} \\
\mathbf{u}(\mathbf{x}, 0) = \mathbf{u}_0(\mathbf{x}) \quad ; \quad \dot{\mathbf{u}}(\mathbf{x}, 0) = \dot{\mathbf{u}}_0(\mathbf{x}) , \quad \mathbf{x} \in \Omega
\end{array} \right.
\end{equation}
where $\mathbf{f}(\mathbf{x},t)$ is a force density field defined in $\Omega$, with dimensions of force per unit area or force per unit volume for $d=2$ and $d=3$, respectively. This force density field is zero everywhere but on the points lying on $\Gamma$ and can be expressed as
\begin{equation}\label{eq:f}
 \mathbf{f}(\mathbf{x},t)  = \mathbf{F}(\mathbf{x},t) \delta_\Gamma(\mathbf{x}),
\end{equation}
where $\delta_\Gamma(\mathbf{x})$ denotes the delta function extended to $\Gamma$ \cite{ONURAL200618} and $\mathbf{F}(\mathbf{x},t)$ denotes a force density applied in $\Gamma$, with  dimensions of force per unit length or force per unit volume, for $d$ = 2 and $d=3$, respectively. The properties of $\delta_\Gamma(\mathbf{x})$ are equivalent to the original Dirac delta,
 $$\int_\Omega \mathbf{f}(\mathbf{x},t)\delta_\Gamma(\mathbf{x})  \mathrm{d}\Omega=\int_\Gamma \mathbf{F}(\mathbf{x},t) \mathrm{d}\Gamma.$$
The value of $\mathbf{F(x},t)$ along $\Gamma$ is obtained to satisfy the prescribed displacement on $\Gamma$
\begin{equation}
\mathbf{u}(\mathbf{x},t) = \mathbf{U}(t) \quad \text{for} \quad \mathbf{x}\in \Gamma,\end{equation}
where $\mathbf{U}(t)$ is the time function defining the value of the displacement vector of all the points in $\Gamma$.
The constitutive equation is the Hooke law for a general heterogeneous medium, characterized by the spatial distribution of its fourth-order stiffness tensor $\mathbb{C}(\mathbf{x})$,
\begin{equation}
\boldsymbol{\sigma}(\mathbf{x})=\mathbb{C}(\mathbf{x}): \boldsymbol{\varepsilon}(\mathbf{x}, t) = \mathbb{C}(\mathbf{x}): \nabla^s \left[\mathbf{u}(\mathbf{x},t) \right]
\end{equation}

\subsubsection{Continuum solution}

As in the one-dimensional case, the Newmark integration scheme is applied to the linear momentum conservation (Eq. (\ref{eq:ND_consv_lm})) by substituting the acceleration definition given by Eq. \eqref{eq:Newmark}. The linear momentum for time $t_n$ can then be written as a function of the displacement field $\mathbf{u}_{n}$ and the force density $\mathbf{f}_n$ in step $n$ (current step) and displacement $\mathbf{u}_{n-1}$, velocity $\dot{\mathbf{u}}_{n-1}$ and acceleration $\ddot{\mathbf{u}}_{n-1}$ in the previous step $n-1$.

\begin{equation}
\label{eq:ND_consv_lm_Newmark}
\nabla \cdot \left[\mathbb{C}(\mathbf{x}): \nabla ^s \mathbf{u}_n\right] + \mathbf{f}_n (\mathbf{x})  = \frac{\rho(\mathbf{x})}{\beta \Delta t^2} \left[\mathbf{u}_{n}-\mathbf{u}_{n-1}- \Delta t\dot{\mathbf{u}}_{n-1} - \Delta t^2 (0.5-\beta)\ddot{\mathbf{u}}_{n-1} \right].
\end{equation}
\noindent
Grouping the terms and multiplying by $-\beta \Delta t^2$ results in
\begin{equation}
\label{eq:ND_consv_lm_Newmark_ord}
-\beta \Delta t^2 \nabla \cdot \left[\mathbb{C}(\mathbf{x}): \nabla ^s \mathbf{u}_n\right] + \rho(\mathbf{x})\mathbf{u}_{n} =  \beta \Delta t^2 \mathbf{f}_n (\mathbf{x})  + \rho(\mathbf{x}) \left[\mathbf{u}_{n-1}+ \Delta t\dot{\mathbf{u}}_{n-1} + \Delta t^2 (0.5-\beta)\ddot{\mathbf{u}}_{n-1} \right].
\end{equation}
Naming the unknown fields $\mathbf{f}_n,\mathbf{u}_n$ as $\mathbf{f},\mathbf{u}$, respectively, Eq. \eqref{eq:ND_consv_lm_Newmark_ord} is transformed into Fourier space, leading to
\begin{equation}
\label{eq:Fourier_ND_lm}
\begin{split}
&	-\beta \Delta t^2 i\xi_k \mathcal{F} \left\lbrace C_{jklm}(\mathbf{x}) \mathcal{F}^{-1} \left\lbrace \dfrac{i}{2} \left( \hat{u}_{l} \xi_m + \hat{u}_{m} \xi_l  \right) \right\rbrace \right\rbrace + \mathcal{F} \left\lbrace \rho(\mathbf{x}) \mathcal{F}^{-1} \{\hat{u}_{j}\} \right\rbrace  = \\
&	\mathcal{F} \left\lbrace \rho(\mathbf{x}) \left(u^{n-1}_{j}+\Delta t \dot{u}^{n-1}_{j} + \Delta t ^2 (0.5-\beta)\ddot{u}^{n-1}_{j} \right)\right\rbrace + \beta \Delta t^2 \mathcal{F} \left\lbrace f_j  \right\rbrace,
\end{split}
\end{equation}
where $\boldsymbol{\xi}$ is the frequency vector. 

As in the one-dimensional case, the left-hand side of Eq. \eqref{eq:Fourier_ND_lm} corresponds to a linear operator $\hat{\mathcal{A}}(\cdot)$ in Fourier space. Therefore, Eq. \eqref{eq:Fourier_ND_lm} can be written as

\begin{equation}
\label{eq:Fourier_ND_lm_compact}
\hat{\mathcal{A}}\left(\hat{\mathbf{u}}\right) = \beta \Delta t^2 \widehat{\mathbf{f}} + \widehat{\mathbf{b}},
\end{equation}
where $\widehat{\mathbf{f}}$ is the Fourier transform of the force density field at time $t=t_n$ (Eq. (\ref{eq:f})) and $\hat{\mathcal{A}}\left(\hat{\mathbf{u}}\right)$ and $\hat{\mathbf{b}}$ are defined as (in components):
\begin{equation}
\label{eq:Fourier_ND_A}
\hat{\mathcal{A}}\left(\hat{u}\right)_j = -\beta \Delta t^2 i\xi_k \mathcal{F} \left\lbrace C_{jklm}(\mathbf{x}) \mathcal{F}^{-1} \left\lbrace \dfrac{i}{2} \left( \hat{u}_{l} \xi_m + \hat{u}_{m} \xi_l  \right) \right\rbrace \right\rbrace + \mathcal{F} \left\lbrace \rho(\mathbf{x}) \mathcal{F}^{-1} \{\hat{u}_{j}\} \right\rbrace
\end{equation}
\begin{equation}
\label{eq:Fourier_ND_b}
\hat{b}_j = \mathcal{F} \left\lbrace \rho(\mathbf{x}) \left(u^{n-1}_{j}+\Delta t \dot{u}^{n-1}_{j} + \Delta t ^2 (0.5-\beta)\ddot{u}^{n-1}_{j} \right)\right\rbrace.
\end{equation}

Due to the linearity of the operator $\hat{\mathcal{A}}\left(\hat{\mathbf{u}} \right)$, the problem can be split into two linear systems whose solutions will be added later. The first linear problem is
\begin{equation}\label{eq:lin_prob1}
    \hat{\mathcal{A}}\left(\hat{\mathbf{u}}_b \right)=\hat{\mathbf{b}}
\end{equation}
whose solution provides the field $\mathbf{u}_b=\mathcal{F}^{-1} \left\lbrace \hat{\mathbf{u}}_b\right\rbrace$. The second problem consists in finding the effect of the force density field $\mathbf{f}$, defined in Eq. \eqref{eq:f}, for a given value of the force $\mathbf{F}$ acting on $\Gamma$ at time $t=t_n$,
\begin{equation}\label{u_p}
\hat{\mathcal{A}}\left(\hat{\mathbf{u}_f} \right)=\hat{\mathbf{f}}.
\end{equation}
The solution $\mathbf{u}_n$ in real space corresponds to the sum of the solution of the two problems (Eqs. \eqref{eq:lin_prob1} and \eqref{u_p}),
\begin{equation}
\label{eq:Sol_Real_ND}
\mathbf{u}_n(\mathbf{x}) = \beta \Delta t^2 \mathbf{u}_f(\mathbf{x})+\mathbf{u}_b(\mathbf{x}).
\end{equation}

At this stage, the value of the force density $\mathbf{F}$ acting at each point of $\Gamma$ is not yet known and should be obtained to meet the prescribed displacement. To obtain its value, Green's functions will be used. Green's function of the operator $\mathcal{A}$ provides the effect of a point force $\mathbf{P}$ concentrated on $\mathbf{x}'$ on the displacement at a point $\mathbf{x}$,
$$
\mathbf{u}(\mathbf{x}) = \mathbf{g} (\mathbf{x},\mathbf{x}')\cdot \mathbf{P}(\mathbf{x}').
$$
The effect of a field $\mathbf{f}$ is then
\begin{equation}\label{eq:conv3D}
\mathbf{u}(\mathbf{x}) = \int_\Omega \mathbf{g}(\mathbf{x},\mathbf{x}')\cdot \mathbf{f}(\mathbf{x}') \mathrm{d}\Omega = \int_\Gamma \mathbf{g}(\mathbf{x},\mathbf{x}')\cdot \mathbf{F}(\mathbf{x}') \mathrm{d}\Gamma = (\mathbf{g}\ast \mathbf{F})_\Gamma,
\end{equation}
where the convolution is performed at the points on $\Gamma$.
The value of the force distributed on $\Gamma$ must be determined at each time step to satisfy the value of the prescribed displacement for the points in $\mathbf{x}\in\Gamma$,
\begin{equation}
\label{eq:new}
\mathbf{u}_n(\mathbf{x}\in \Gamma) =\mathbf{U}(t)=\beta \Delta t^2 (\mathbf{g}\ast \mathbf{F}(\mathbf{x},t))_\Gamma
+ \mathbf{u}_b (\mathbf{x}\in \Gamma) 
\end{equation}
and can be obtained inverting the previous equation,
\begin{equation} \label{eq:deconvolution}
\text{find} \ \mathbf{F}(\mathbf{x},t) \ | \ (\mathbf{g}\ast \mathbf{F}(\mathbf{x},t))_\Gamma  =\frac{1}{\beta \Delta t^2}\left(\mathbf{U}(t)-\mathbf{u}_b(\mathbf{x}) \right)\ \text{for} \ \mathbf{x}\in \Gamma \rightarrow \mathbf{F}(\mathbf{x},t)=\frac{1}{\beta \Delta t^2}\mathbf{g}^{-1}\ast  \left(\mathbf{U}(t)-\mathbf{u}_b(\mathbf{x})\right)
\end{equation}
where $\mathbf{g}^{-1}$ stands for the inverse of Green's function.
Solving Eq. \eqref{eq:deconvolution} can be easily handled when the manifold $\Gamma$ is discretized. In this case,  $\mathbf{g}(\mathbf{x,x}')$ can be replaced by a matrix that links the forces and displacement of every point in $\Gamma$.
This matrix is formed by solving Eq. \eqref{u_p} for each point in $\Gamma$ and its inverse, $\mathbf{g}^{-1}$, is calculated and only once and then is stored and used during the simulation. At each time step, $\mathbf{F}(\mathbf{x},t)$ is obtained using that inverse matrix. Finally, the force field $\mathbf{f}(\mathbf{x},t)=\delta_\Gamma\mathbf{F}(\mathbf{x},t)$ is used as the RHS of Eq. \eqref{u_p} to obtain $u_f$, and the displacement field is obtained using Eq. \eqref{eq:Sol_Real_ND}.

\subsubsection{Spatial discretization and algorithm}
Focusing on the numerical resolution of the three-dimensional case, the periodic domain is discretized in $N_1\cdot N_2\cdot N_3$ equispaced voxels, whose centers are given by
\begin{equation}
\mathbf{x}_{i,j,k}=\left(\frac{L_1}{N_1}(\frac{1}{2}+i),\frac{L_2}{N_2}(\frac{1}{2}+j),\frac{L_3}{N_3} (\frac{1}{2}+k)   \right)
\end{equation}
with $(i,j,k) \in [0,N_1-1]\times[0,N_2-1]\times[0,N_3-1]$. The values of the properties of the material in each voxel $\mathbb{C}(\mathbf{x}_{i,j,k})$ and $\rho(\mathbf{x}_{i,j,k})$ correspond to the property of the material that occupies that position. The displacement vector is discretized using its value in the center of each voxel, $\mathbf{u}(\mathbf{x}_{i,j,k})$.

The discrete fields in Fourier space have the same dimension as their counterparts in real space. The corresponding $N_1\cdot N_2 \cdot N_3$ discrete frequencies in Fourier space are  
\begin{equation}\label{eq:freq3D}
\boldsymbol{\xi}_{k_1,k_2,k_3}=\left( {\xi}_{k_1}, {\xi}_{k_2}, {\xi}_{k_3}\right)
\end{equation}
 being the frequencies in each direction $\xi_{k_1},\xi_{k_2},\xi_{k_3}$ defined as in the 1D case, (Eq. \ref{frequencies}). 
 
 After discretization, the linear operator $\hat{\mathcal{A}}$ defined in equation \eqref{eq:Fourier_ND_A} is replaced by its discrete counterpart in which the Fourier and inverse Fourier transforms are substituted by discrete Fourier transforms, and the frequencies used are the discrete set defined in equation \eqref{eq:freq3D}. The equations of the type
$ \hat{\mathcal{A}}(\mathbf{u})=\mathbf{b}$, as Eqs. (\eqref{eq:lin_prob1},\eqref{u_p}),
 correspond to a linear system of equations, which can be solved efficiently using Krylov iterative solvers. In the present case, the conjugate gradient can be used because of the Hermitian nature of the operator. Moreover, as proposed in \cite{Lucarini2019c} and \cite{Segurado2021} for similar problems, preconditioners can be built that strongly improve numerical performance. The proposed preconditioner is the exact solution of the linear equation for a homogeneous medium with average volume stiffness and density $\mathbb{C}=\overline{\mathbb{C}(\mathbf{x})}$ and $\rho=\overline{\rho(\mathbf{x})}$, and corresponds to 
 \begin{equation}\label{eq:precond}
     \hat{\mathcal{M}}(\hat{\mathbf{u}})= \left(\beta \Delta t^2 \mathbf{K(\boldsymbol{\xi})}+\overline{\rho}\mathbf{I}\right)^{-1}\cdot \hat{\mathbf{u}}(\boldsymbol{\xi}),
 \end{equation}
where $\mathbf{I}$ is the identity tensor and $\mathbf{K}(\boldsymbol{\xi})$ is the acoustic tensor, defined for each frequency as
 \begin{equation}\label{eq:acous}
 K_{ij}(\boldsymbol{\xi})=\overline{C}_{ikjm}\xi_k\xi_m .
  \end{equation}
 The preconditioner, Eq.\eqref{eq:precond}, is formed and stored once at the beginning of the simulation. Its computation is not expensive computationally since it only requires the inversion of a 3$\times$3 matrix for each frequency. This step is indeed equivalent to form the Gamma operator in a standard homogenization FFT simulation. The use of the preconditioner within the conjugate gradient implies just adding a multiplication of the current value of $\hat{\mathbf{u}}$ by $\hat{\mathcal{M}}$ at each iteration.

For obtaining the force distribution on $\Gamma$ at every time step, the Eq. \eqref{eq:deconvolution} have to be solved. To perform this efficiently, the discrete inverse Green operator will be constructed once at the beginning of the simulation and used in subsequent time steps. Let $\mathbf{x}_p$ be the position of a point lying in $\Gamma$, $1\leq p \leq N_p$ and let $\mathbf{f}_p^j(\mathbf{x})$, be a force density field defined as
\begin{equation}\label{eq:f_discrete}
 \mathbf{f}_p^j(\mathbf{x} )= \left\{ 
\begin{array}{cc}
 \mathbf{e}_j & \text{if} \ \mathbf{x}=\mathbf{x}_p \\
\mathbf{0} & \text{elsewhere}
\end{array}
\right. 
\end{equation}
for every point $p$ and direction $j=1,2,3$. The displacement result of applying this field is obtained by solving the linear system,
\begin{equation}\label{eq:solve_green}
\hat{\mathcal{A}}(\hat{\mathbf{u}}_{(p,j)})=\hat{\mathbf{f}}^j_p\end{equation}
where $\mathbf{u}_{(p,j)}(\mathbf{x})=\mathcal{F}^{-1}(\hat{\mathbf{u}}_{(p,j)})$. 
The result of solving the problem for each applied field around $p$ allows to form the Green's tensor for each other point $q$ lying on $\Gamma$, $\mathbf{g}(q,p), p,q=1,N_p$, with dimensions $3\times3$ which is defined as
\begin{equation}\label{eq:g}
g^{ij}(q,p) =u^i_{(p,j)}(\mathbf{x}_q)
\end{equation}
The assembly of all $\mathbf{g}(p,q)$ matrices defines a new matrix, $\mathbf{G}$, with dimension $3N_p\times 3N_p$ as
\begin{equation}\label{eq:G}
G_{IJ} = g^{ij}(q,p) \ \mathrm{with} \ I=3\cdot(p-1)\cdot(i-1)\ ,\ J=3\cdot(q-1)\cdot(j-1) \quad p,q=1,...,N_p;\ i,j=1,2,3 
\end{equation}
which express the relation between displacement and the forces of the $N_p$ nodes in $\Gamma$. The convolution restricted to the manifold $\Gamma$ defined in Eq. \eqref{eq:conv3D} is therefore written, for the discrete case, as the matrix vector product in Eq. \eqref{eq:solve_green2}.
\begin{equation}\label{eq:solve_green2}
\begin{bmatrix}
u^1_1 \\
u^2_1\\
u^3_1\\
\vdots \\
u^1_{Np} \\
u^2_{Np} \\
u^3_{Np} \\
\end{bmatrix}=\begin{bmatrix}
\mathbf{g}(1,1) & \cdots & \mathbf{g} (1,N_p) \\
\vdots & \ddots & \vdots \\
\mathbf{g} (N_p,1) & \cdots & \mathbf{g} (N_p,N_p)
\end{bmatrix}
\begin{bmatrix}
F^1_1 \\
F^2_1\\
F^3_1\\
\vdots \\
F^1_{Np} \\
F^2_{Np} \\
F^3_{Np} \\
\end{bmatrix}
\end{equation}
The matrix $\mathbf{G}$ is non-singular but is in general non-symmetric for a heterogeneous medium, since in this case the Green's operator is not translation invariant. $\mathbf{G}$ is computed once at the beginning of the simulation and then it is inverted using a direct method to obtain $\mathbf{G}^{-1}$, which will be stored. Finally, from the prescribed displacement for all nodes in $\Gamma, \mathbf{u}(\mathbf{x}, t) = \mathbf{U}(t)$, the vector $\mathbf{V}(t)$ can be formed;
\begin{equation}\label{eq:Vt}
\mathbf{V}=\begin{bmatrix}
U^1(t)-{u_b^1}_1 \\
U^2(t)-{u_b^2}_1 \\
U^3(t)-{u_b^3}_1 \\
\vdots \\
U^1(t)-{u_b^1}_{N_p} \\
U^2(t)-{u_b^2}_{N_p} \\
U^3(t)-{u_b^3}_{N_p} \\
\end{bmatrix}
\end{equation}
and the vector $\mathbf{F}$ for the points in $\Gamma$ can be calculated as

\begin{equation}
\label{eq:L-determination_ND}
\mathbf{F} (t) = \frac{1}{\beta \Delta t^2} \mathbf{G}^{-1}\cdot \mathbf{V}(t)
\end{equation}
and the solution $\mathbf{u}_n(\mathbf{x})$ is obtained for every point of equation \eqref{eq:Sol_Real_ND} once $\mathbf{F}$ is obtained.

The algorithm described in this section is presented in the Algorithm box \ref{AL:3D}. The method is programmed in Python, integrated into the FFT homogenization code FFTMAD \cite{Lucarini2019a,Lucarini2019c}. FFTW is used for Fourier transforms (through the \emph{pyfftw} project), and the \emph{scipy.sparse} functions are the core for the resolution of iterative linear problems.\\

\begin{algorithm}[H]\label{AL:3D}
\SetAlgoLined
\caption{FFT based elastodynamic algorithm --2D and 3D--}
\KwData{
Microstructure: $\mathbb{C}(\mathbf{x}), \rho(\mathbf{x})$ and definition of $\Gamma$: $\mathbf{x}_p \in \Gamma$\\
Prescribed fields on $\Gamma$:  $\mathbf{U}(t), \dot{\mathbf{U}}(t), \ddot{\mathbf{U}}(t)$, \\
Initial conditions: $\mathbf{u}_0(\mathbf{x}),\dot{\mathbf{u}}_0(\mathbf{x})$ ,\\
Time discretization [$t_0, n \Delta t, ..., t_f $] and space discretization $x = [0, x_m, ..., L]$, for $n$ and $m$ $\in \mathbb{Z}$; \\
Newmark constants, $\beta=0.25$, $\gamma=0.5$}
\KwResult{$\mathbf{u}(\mathbf{x},t), \dot{\mathbf{u}}(\mathbf{x},t), \ddot{\mathbf{u}}( \mathbf{x},t )$, $\mathbf{F}(t)$}	
{\bf Preprocessing:} 
Form function $\hat{\mathcal{A}}()$ eq. \eqref{eq:Fourier_ND_A} and compute $\ddot{\mathbf{u}}_0(\mathbf{x})=\mathcal{F}^{-1}\left(
i\boldsymbol{\xi} \mathcal{F} 
\left\lbrace 
\mathbb{C}(\mathbf{x}) : \mathcal{F}^{-1} 
\left\lbrace \dfrac{i}{2} 
\left( \hat{\mathbf{u}}_0 \otimes \boldsymbol{\xi} +\boldsymbol{\xi}  \otimes \hat{\mathbf{u}}_0 \right)
 \right\rbrace 
 \right\rbrace \right)$;\\
Obtain matrix $\mathbf{G}^{-1}$:
\begin{itemize}
    \item Solve $N_p$ linear problems in Eq.\eqref{eq:solve_green} for each point $p$ in $\Gamma$ and each d.o.f. using CG with preconditioner $\hat{\mathcal{M}}$  (Eq. \ref{eq:precond})
    \item Form matrix $\mathbf{g}(p,q)$ for every $p,q=1,N_p$ from the solution (Eq. \ref{eq:g})
    \item Form matrix $\mathbf{G}$ (Eq. \ref{eq:G})
    \item Invert matrix $\mathbf{G}$ using a direct $L,U$ solver. The factorization in matrices $\mathbf{L,U}$ is stored.
\end{itemize}

\While{$t \leq t_f$}{
$t = t_0 + \Delta t $\;
$\hat{\mathbf{b}}_{n} = \mathcal{F}\left\lbrace\rho(\mathbf{x}) \left(\mathbf{u}_{n-1} + \Delta t \dot{\mathbf{u}}_{n-1} + \Delta t^2 (0.5-\beta) \ddot{\mathbf{u}}_{n-1} \right)\right\rbrace$\;
$\mathbf{u}_{b_n} = \mathcal{F}^{-1}\left\lbrace \hat{\mathcal{A}}^{-1}(\hat{\mathbf{b}}_n) \right\rbrace$; $\hat{\mathcal{A}}^{-1}(\cdot)$ using CG ($\hat{\mathcal{A}}, \hat{\mathbf{b}}_n$)  with preconditioner $\hat{\mathcal{M}}$ (Eq. \ref{eq:precond})\;
Form $\mathbf{V}(t)$ using $\mathbf{u}_{b_n}$ and $\mathbf{U}(t)$ (Eq. \ref{eq:Vt})\;
$\mathbf{F} (t) = \frac{1}{\beta \Delta t^2}\mathbf{G}^{-1} \cdot \mathbf{V}(t)$\ using the factorization of $\mathbf{G}$\;
$\mathbf{u}_{f} =\beta \Delta t^2  \mathcal{F}^{-1}\left\lbrace \hat{\mathcal{A}}^{-1}(\hat{\mathbf{f}}_n) \right\rbrace$; $\hat{\mathcal{A}}^{-1}(\cdot)$ using CG ($\hat{\mathcal{A}}, \hat{\mathbf{f}}_n$)  with preconditioner $\hat{\mathcal{M}}$ (Eq. \ref{eq:precond})\;
$\mathbf{u}_n(\mathbf{x}) =\mathbf{u}_{f}(\mathbf{x})+\mathbf{u}_{b_n}(\mathbf{x})$\;
$\ddot{\mathbf{u}}_n(\mathbf{x}) = \frac{1}{\beta \Delta t^2} \left( \mathbf{u}_{n}-\mathbf{u}_{n-1}- \Delta t \dot{u}_{n-1} - \Delta t ^2 (0.5-\beta)\ddot{\mathbf{u}}_{n-1} \right)$\;
$\dot{\mathbf{u}}_n(\mathbf{x}) = \dot{\mathbf{u}}_{n-1} + \Delta t (1-\gamma)\ddot{\mathbf{u}}_{n-1} + \Delta t \gamma \ddot{\mathbf{u}}_{n}$\;
}
\end{algorithm}

\subsection{Explicit Integration Scheme}\label{sec:Explicit_integration}
An explicit algorithm based on FFT is developed to solve the wave propagation problem, in order to compare its accuracy and efficiency with the implicit scheme and FE. The method chosen is the explicit central difference scheme, which is the particularization of the Newmark-$\beta$ method for $\beta=0$ and $\gamma=1/2$.

Given a constant time increment $\Delta t$, the displacement and velocity fields in current time are obtained using $\beta=0$ and $\gamma=1/2$ in Eq.\eqref{eq:Newmark}, leading to
\begin{equation}
\label{eq:Central_diff_rule}
\begin{split}
    &{u}_{n} = {u}_{n-1} + \Delta t \dot{u}_{n-1} + \dfrac{{\Delta t}^2}{2} \ddot{u}_{n-1}\\
	&\dot{u}_{n}=\dot{u}_{n-1} + \frac{1}{2} \Delta t \ddot{u}_{n-1} + \frac{\Delta t }{2} \ddot{u}_{n}
\end{split}
\end{equation}
Note that now, the displacement at current time step $n$ is obtained explicitly from velocities and acceleration from previous time step $n-1$.
With respect to stability, the explicit scheme is conditionally stable and the maximum time increment for stability is given by the so-called Courant–Friedrichs–Lewy condition, usually called the $CFL$ time. $CFL$ time increment is defined from the highest eigenvalue $\omega_{max}$ of the problem.
\begin{equation}
    \label{eq:Courant_stable}
    \Delta t \leq \dfrac{2}{\omega_{max}}.
\end{equation}

\subsubsection{One dimensional case}
The problem to be solved is the one described in section \ref{subsec:one_dim} and consists in the propagation of an elastic wave in a one-dimensional domain with a prescribed condition on a point. Introducing the expressions of $u_n$ and $\ddot{u}_n$, given in Eq. \eqref{eq:Central_diff_rule}, in the conservation of linear momentum of the problem (Eq. \ref{eq:1D_consv_lm}) leads to

\begin{equation}
    \label{eq:Explicit_linear_momentum_1D}
    \frac{\mathrm{d}}{\mathrm{d}x}\left[E(x)\frac{\mathrm{d}}{\mathrm{d}x} \underbrace{\left(
    {u}_{n-1} + \Delta t \dot{u}_{n-1} + \dfrac{{\Delta t}^2}{2} \ddot{u}_{n-1}\right)}_{u_n}
        \right]+f_n(x) = 2
\frac{\rho(x)}{\Delta t} \left[\dot{u}_{n} - \dot{u}_{n-1} - \Delta t \frac{1}{2}\ddot{u}_{n-1} \right]
\end{equation}
where the velocity term $\dot{u}_{n}$ is the unknown that must be determined. 
Grouping terms (unknown on the left and known on the right), the equation can be written as:

\begin{equation}
    \label{eq:Explicit_linear_momentum_arranged_1D}
   \dot{u}_{n}=
\left[\dot{u}_{n-1} + \frac{1}{2} \Delta t \ddot{u}_{n-1} \right] + 
 \frac{\Delta t}{2\rho(x)}
\frac{\mathrm{d}}{\mathrm{d}x}\left[E(x)\frac{\mathrm{d}}{\mathrm{d}x} u_n(x) \right]+\frac{\Delta t}{2\rho(x)} f_n(x) = \dot{u}_b(x) +\frac{\Delta t}{2\rho(x)} f_n(x)  
\end{equation}

In Eq. \eqref{eq:Explicit_linear_momentum_arranged_1D} the right hand side depends only on fields at the previous time step $n-1$, with the exception of the force $f_n$. The term dependent on fields at $t_{n-1}$, $\dot{u}_b$, and can be computed applying the differential operator in Fourier space as
\begin{equation}\label{eq:b_explicit}
\dot{u}_b(x)=\left[\dot{u}_{n-1} + \frac{1}{2} \Delta t \ddot{u}_{n-1} \right] + 
 \frac{\Delta t}{2\rho(x)} \mathcal{F}^{-1}\left\{ 
 i\xi \mathcal{F}\left\{ E(x) \mathcal{F}^{-1} [i\xi \widehat{u}_n] \right\}
 \right\}.
\end{equation}
The value of the body force $f_n(x)$ can be obtained directly form the prescribed velocities on $\Gamma$, here the point $x_0$, 
$$\dot{u}_n(x=x_0)=\dot{U}(t)=\dot{u}_b(x_0)+\frac{\Delta t}{2\rho(x)} f_n(x)$$
leading to
\begin{equation}
F_n=\frac{2\rho(x)}{\Delta t}(\dot{U}(t)-\dot{u}_b(x_0))
\end{equation}
It is interesting to note that, contrary to the implicit algorithm (Eq. \ref{eq:L-determination}), there is no need to pre-compute the effect of the body force to obtain its value.
\subsubsection{Three dimensional case}
The extension of the 1D algorithm to three dimensions is straightforward and only the resulting equations will be provided.

The velocity field at time $n$ is obtained by introducing the time integration expressions in Eq. \eqref{eq:Central_diff_rule} into the conservation of linear momentum (Eq. \ref{eq:ND_consv_lm}).
\begin{equation}
    \label{eq:Explicit_linear_momentum_arranged_3D}
   \dot{\mathbf{u}}_{n}=
\left[\dot{\mathbf{u}}_{n-1} + \frac{1}{2} \Delta t \ddot{\mathbf{u}}_{n-1} \right] + 
 \frac{\Delta t}{2\rho(\mathbf{x})} \nabla\cdot \left[
 \mathbb{C}(\mathbf{x}):\nabla^s\mathbf{u}_n
 \right]
+\frac{\Delta t}{2\rho(\mathbf{x})} \mathbf{f}_n(\mathbf{x}) = \dot{\mathbf{u}}_b(\mathbf{x}) +\frac{\Delta t}{2\rho(\mathbf{x})} \mathbf{f}_n(\mathbf{x})  
\end{equation}
where $\dot{\mathbf{u}}_b(\mathbf{x})$ is computed similarly to Eq. \ref{eq:b_explicit} by differentiation in Fourier space as
$$\dot{u}_{b\ (j)}(\mathbf{x}) = \left[\dot{u}_{n-1 \ (j)} + \frac{1}{2} \Delta t \ddot{u}_{n-1 \ (j) } \right] + \frac{\Delta t}{2\rho(\mathbf{x})} 
\mathcal{F}^{-1} \left\lbrace 
i\xi_k \mathcal{F} \left\lbrace C_{jklm}(\mathbf{x}) \mathcal{F}^{-1} \left\lbrace \dfrac{i}{2} \left( \hat{u}_{l} \xi_m + \hat{u}_{m} \xi_l  \right) \right\rbrace \right\rbrace  \right\rbrace
$$

The value of the force field $\mathbf{f}_n(\mathbf{x})$ is set using the condition of the velocity $\dot{\mathbf{U}}(t)$in $\Gamma$ as
\begin{equation}
\mathbf{f}_n(x)=
\left \{ \begin{array}{cc}
0 & \mathbf{x} \not\in \Gamma \\
\frac{2\rho(\mathbf{x})}{\Delta t}(\dot{\mathbf{U}}(t)-\dot{\mathbf{u}}_b(\mathbf{x})) & \mathbf{x} \in \Gamma
\end{array}
\right.
\end{equation}

\subsubsection{Stable time increment}
For the stability of the solution, the maximum $\Delta t$ to be used in a simulation must be computed based on the eigenvalues of the discrete problem. In the case of a homogeneous medium ($\mathbb{C}(\mathbf{x})=\mathbb{C}_0$ and $\rho(\mathbf{x})=\rho$), the eigenvalues $\omega$ of the problem are obtained from the eigenvalues of the acoustic tensor. Following the work of Segurado and Lebensohn \cite{Segurado2021}, the eigenvalues of the problem for a given spatial frequency vector $\boldsymbol{\xi}$ are the three eigenvalues of the matrix $\mathbf{Q}$
\begin{equation}
    \label{eq:Eigenvalues_1D_homogeneo}
    Q_{ik}(\boldsymbol{\xi})=\frac{1}{\rho_0} [C_{ijkl}\xi_j\xi_l].
\end{equation}
The three largest eigenvalues $\omega_{1,2,3}$ are obtained substituting in Eq. \eqref{eq:Eigenvalues_1D_homogeneo} the largest  frequency vector of the discretization used (Eq. \ref{frequencies}), $\xi_i=\frac{N\pi}{L_i}$, and correspond to the longitudinal and the two transverse modes. Finally, stability is given by
\begin{equation}
    \label{eq:Courant_1D_homogeneo}
    \Delta t \leq \dfrac{2}{\mathrm{max}(\omega_1,\omega_2,\omega_3)}.
\end{equation}
In the case of a heterogeneous medium, the eigenvalues do not have a closed form expression and can be found by solving the discrete eigenvalue problem, as proposed in \cite{Segurado2021}.

\section{Validation and accuracy of the approach}
In this section, the developed numerical framework will be validated by comparing it with analytical solutions and the results obtained with the finite element method. All simulations have been performed using the same computer, a desktop computer with 6 cores i5-9500@3GHz and 8 GB RAM memory. The simulation times refer to this system.

\subsection{One dimension}

First, the proposed approach will be used to study the propagation of elastic waves in a one-dimensional medium. Two different problems will be analyzed: (1) wave propagation in a homogeneous medium and (2) wave propagation in a layered medium. In both cases, analytical solutions can be obtained and are described in Appendix \ref{anex:1Dsol}.

The simulated example consists of a 1D periodic domain $\Omega=\{0\leq x \leq L\}$ with $L = 2$m. A prescribed bell-shaped pulse, $U(t)$, is set at point $x_0=0$,  given by Eq. \eqref{eq-prescribed_displacement} 
\begin{equation}
\label{eq-prescribed_displacement}
U(t) = 
\begin{cases}
\dfrac{A\left( t(\pi / \omega-t)\right)^{\alpha}}{\left((\pi/\omega)^2 / 4\right)^{\alpha}} & \text{for} \ \  0 \leq t \leq \pi / \omega\\
0 & \text{ for the rest}
\end{cases}
\end{equation}
where the parameters that define the shape of the pulse are  $A=0.001$m, $\alpha = 4$ and $\omega = 5\pi c_{0} /L \ $s$^{-1}$ with $c_0$ the wave velocity in the one-dimensional medium.

The material considered in the simulation of the homgeneous domain is aluminum, while in case (2), two different configurations of materials were considered to have different contrast of elastic properties: aluminum-iron and aluminum-uranium. The elastic constants of the three materials are given in Table \ref{tab: 1D-mat-parameters}.  The wave velocity for longitudinal 1D waves $c_0=\sqrt{E/\rho}$ is also provided in the table. 

\begin{table}[h]
	\renewcommand*{\arraystretch}{1.5}
	\centering
	\caption{Parameters of metallic materials. Data from \cite{Meyers}}
	\label{tab: 1D-mat-parameters}      
	\begin{tabular}{lcccccc}
		\hline
		Material & $E$ (GPa)  & $\rho$ $(\mathrm{kg / m^3})$ & $\lambda$ (GPa) &$\mu$ (GPa) & $\nu$ & $c_0$ in 1D (m/s) \\\hline
		Aluminium & 70.3 & 2700.0 & 58.2 & 26.1 & 0.345 & 5102.6\\
		Iron & 211.4 & 7850.0 & 115.7 & 81.6 & 0.293 & 5189.4 \\
		Uranium & 172.0 & 18950.0 & 99.2 & 66.1 & 0.3 & 3012.7
		\\\hline
	\end{tabular}
\end{table}

The simulation time was set as $3.919 \cdot 10^{-4} \mathrm{s^{-1}}$, which is approximately the time the wave needs to travel the entire length of the aluminum bar. 

The homogeneous problem is solved using four different approaches, the implicit and explicit FFT schemes introduced in previous section and implemented in FFTMAD, and implicit and explicit finite element solvers.  
The finite element solutions are obtained using Abaqus \cite{Abaqus2013}. The implicit solution in FE is obtained with the same implicit Newmark solver and a unidimensional mesh of linear truss elements with the same number of elements as the points used in the FFT approach. Periodic boundary conditions 
$$u(x=0)=u(x=L)$$ are applied using multipoint constraints. The prescribed pulse is introduced as a time-dependent Dirichlet boundary condition at the point $x=0$. The explicit simulations are also performed in Abaqus, using in this case central differences with an identical integration scheme as the one proposed in Section 3 for the FFT approach. To avoid instabilities produced by the use of multipoint constraints in explicit FE \cite{SADABA2019434}, the periodic boundary conditions are prescribed in this case applying the same pulse at the two model external nodes.

The 1D bar is discretized with eight different number of voxels, ranging from $N$ = 256 to $N$=32768 voxels. The problem is solved also using different time increments. The reference time increment was the one given by Courant's condition, computed using the standard definition in explicit FE, $CFL=\Delta x/c_0$, and which varies with discretization. The time increments used correspond to $\Delta t / CFL = \left[ 0.5, \ 1, \ 2, \ 3, \ 4, \ 8, \ 10, \ 15, \ 20, \ 50 \right]$. The 80 feasible combinations were simulated with the implicit FFT method presented, while only the problem with $N$=6561 was used with the implicit FE ($\Delta t / CFL = 0.5, 1, 2, 5, 10, 15$),  explicit FE ($\Delta t / CFL = 0.1, 0.7, 1$)  and  explicit FFT ($\Delta t / CFL = 0.1, 0.5, 1$) solvers for comparative purposes. 

\begin{figure}[h]
	\captionsetup[subfigure]{labelformat=empty}
	\centering
	\subfloat{
		\label{fig11}
		\includegraphics[width=0.33\textwidth,trim=0 0cm 0 0cm,clip]{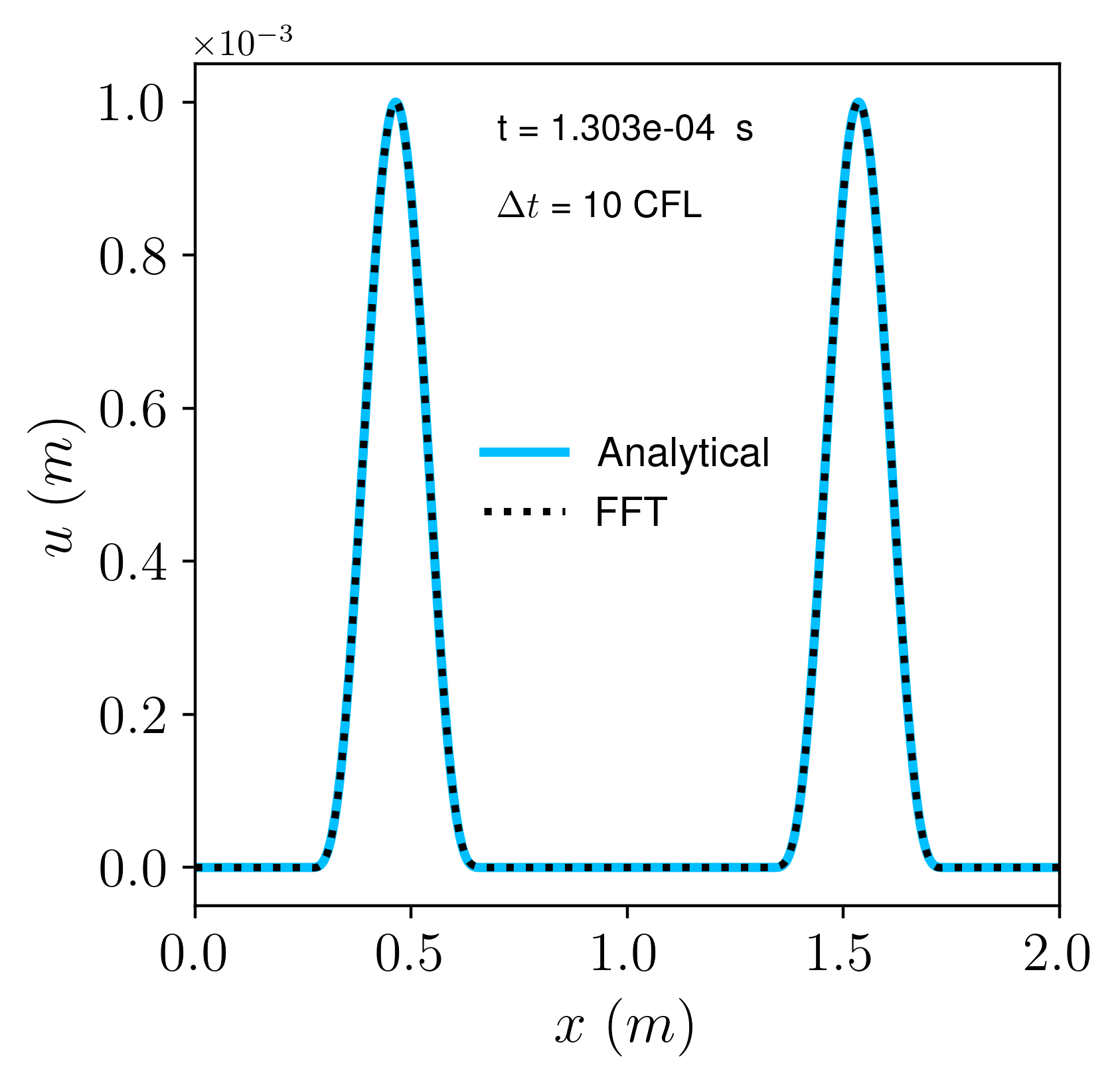}}
	\subfloat{
		\label{fig12}
		\includegraphics[width=0.33\textwidth,trim=0 0cm 0 0cm,clip]{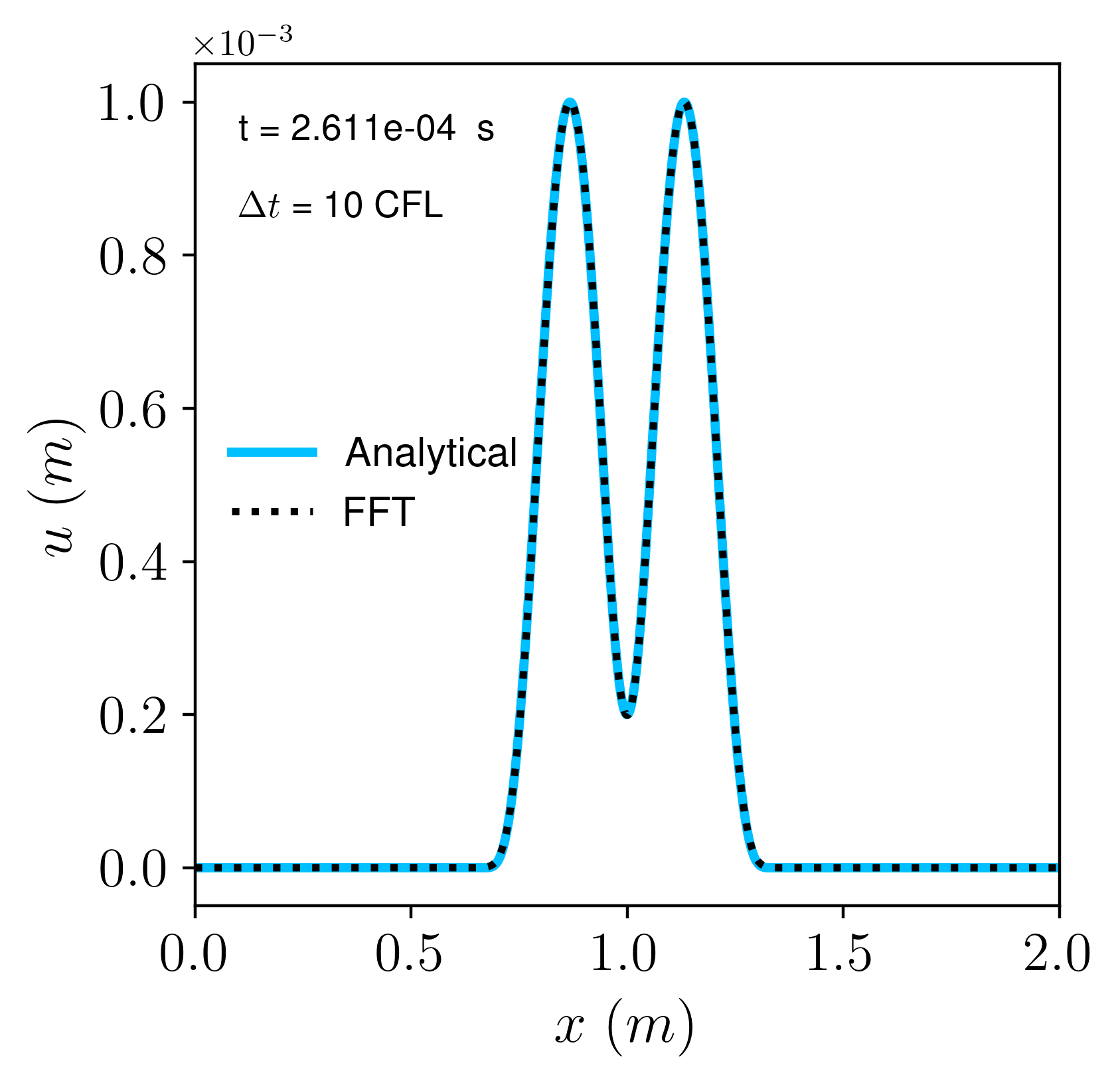}}
	\subfloat{
		\label{fig13}
		\includegraphics[width=0.33\textwidth,trim=0 0cm 0 0cm,clip]{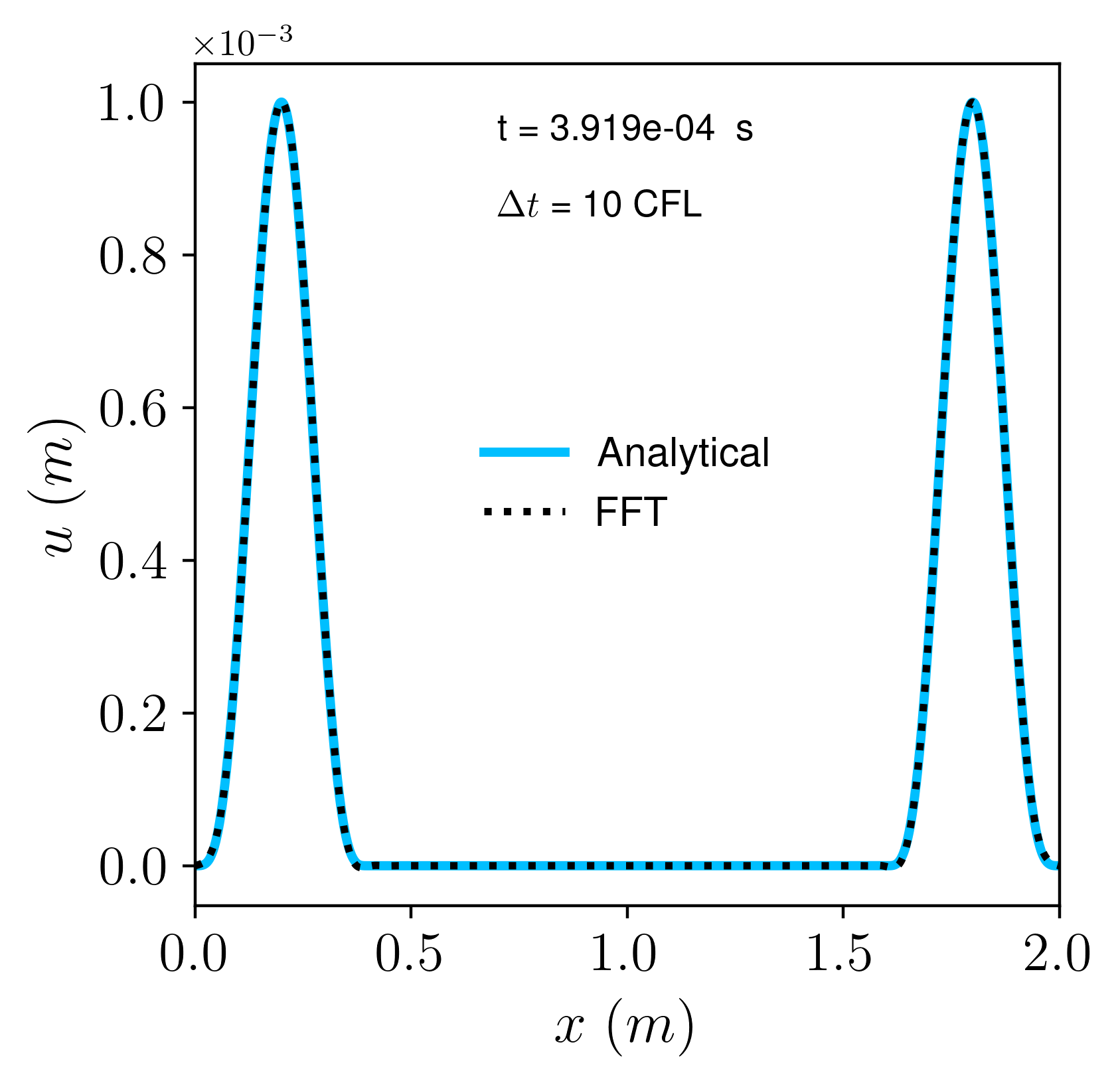}}
	\caption{1D wave propagation (displacement field) along the Al homogeneous domain discretized with 6561 voxels. }
	\label{fig_1D-homogeneous}       
\end{figure}

To illustrate the solution of the problem and the numerical response obtained, a fixed discretization of $N=6561$ voxels is analyzed first. The displacement field in the homogeneous material given by the analytical solution (Annex \ref{anex:1Dsol}) is represented in Fig. \ref{fig_1D-homogeneous} for three different times. The graphs also include the solution obtained by the implicit FFT method proposed here, using a time increment of $\Delta t=10CFL$. The solutions of explicit FFT and FE are very similar and are not represented in the graph for the sake of clarity. Two pulses traveling in opposite directions are obtained, a consequence of the periodic boundary conditions. The perturbation is applied on $x=0$ therefore, one pulse propagates from that point to the right. The other pulse propagates from that point to the left, appearing then in the right corner, which is the periodic point of $x=0$, and moving to the left. In the figure, it can be seen that qualitatively the numerical response is indistinguishable from the analytical solution.  To quantify the errors, in table \ref{tab: 1D-propagation} the $L_2$-norm of the difference between the analytical solution $u(x)$ and the numerical ones $u^h(x)$,
\begin{equation}\label{eq:normL2}
error =\frac{\| u(x)-u^h(x)\|_{L_2}}{\|u(x)\|_{L_2}}.
\end{equation} 
is presented. It can be observed that using implicit FFT, the relative error increases slightly with increasing time step but is always behind $10^{-3}$. The error using explicit FE was always well above the implicit schemes, around one order of magnitude larger. In the case of the explicit FFT, the error was similar to the implicit version. Finally, the implicit FE error was also greater than that of FFT (around a factor of 2). 

\begin{table}[h]
	\renewcommand*{\arraystretch}{1.5}
	\centering
	\caption{Accuracy and performance of simulation results for Al homogeneous domain discretized with 6561 voxels. 1 core computation}
	\label{tab: 1D-propagation}      
	\begin{tabular}{ccccc}
		\hline
		$\Delta$t/CFL & $t$ (s) \tiny{ABQ}  & $\mathrm{e_r}$ \tiny{ABQ} & $t$ (s) \tiny{FFT} &$\mathrm{e_r}$ \tiny{FFT} 
		\\\hline
		1 & - & - & 7.23  & $\mathrm{2.59 \cdot 10^{-3}}$ \\
		2 & - & - & 4.26  & $\mathrm{2.48 \cdot 10^{-3}}$ \\
		5 & - & - & 1.66  & $\mathrm{2.09 \cdot 10^{-3}}$ \\
		10 & 184 & $\mathrm{7.48 \cdot 10^{-3}}$ & 0.92  & $\mathrm{5.36 \cdot 10^{-3}}$ 
		\\\hline
		\multicolumn{5}{c}{Explicit scheme} \\
		1 & .4 & $\mathrm{1.20 \cdot 10^{-1}}$ & 4.88 & $\mathrm{2.64 \cdot 10^{-3}}$
		\\\hline
	\end{tabular}
\end{table}

The second one dimensional case analyzed is a layered medium. The periodic domain is made, in the first example, of aluminum from $x=0$ to $x=0.6$mm, iron from 0.6 to 1.2 mm, and again aluminum between 1.2 mm and 2 mm. In the second example iron is replaced by uranium, achieving a higher property contrast at the materials interface. Again, the solution at three different times is shown in Fig. \ref{fig_1D-layered medium} together with the analytical solution, provided in Appendix \ref{anex:1Dsol}. The pulse, when moving to a different medium, is partially reflected and transmitted, and the numerical methods are able to accurately reproduce this effect, providing a response almost coincident with the analytical solution (Fig.  \ref{fig_1D-layered medium}). Differences between numerical approaches and analytical solutions are quantified using Eq. \eqref{eq:normL2} and represented in Table \ref{tab: 1D-propagation_hetero}. Conclusions are similar to the homogeneous case, errors are very small (about $10^{-3}$) and similar for both FE and FFT. 

\begin{figure}[h]
	\captionsetup[subfigure]{labelformat=empty}
	\centering
	\subfloat{
		\label{a}
		\includegraphics[width=0.33\textwidth,trim=0 0cm 0 0cm,clip]{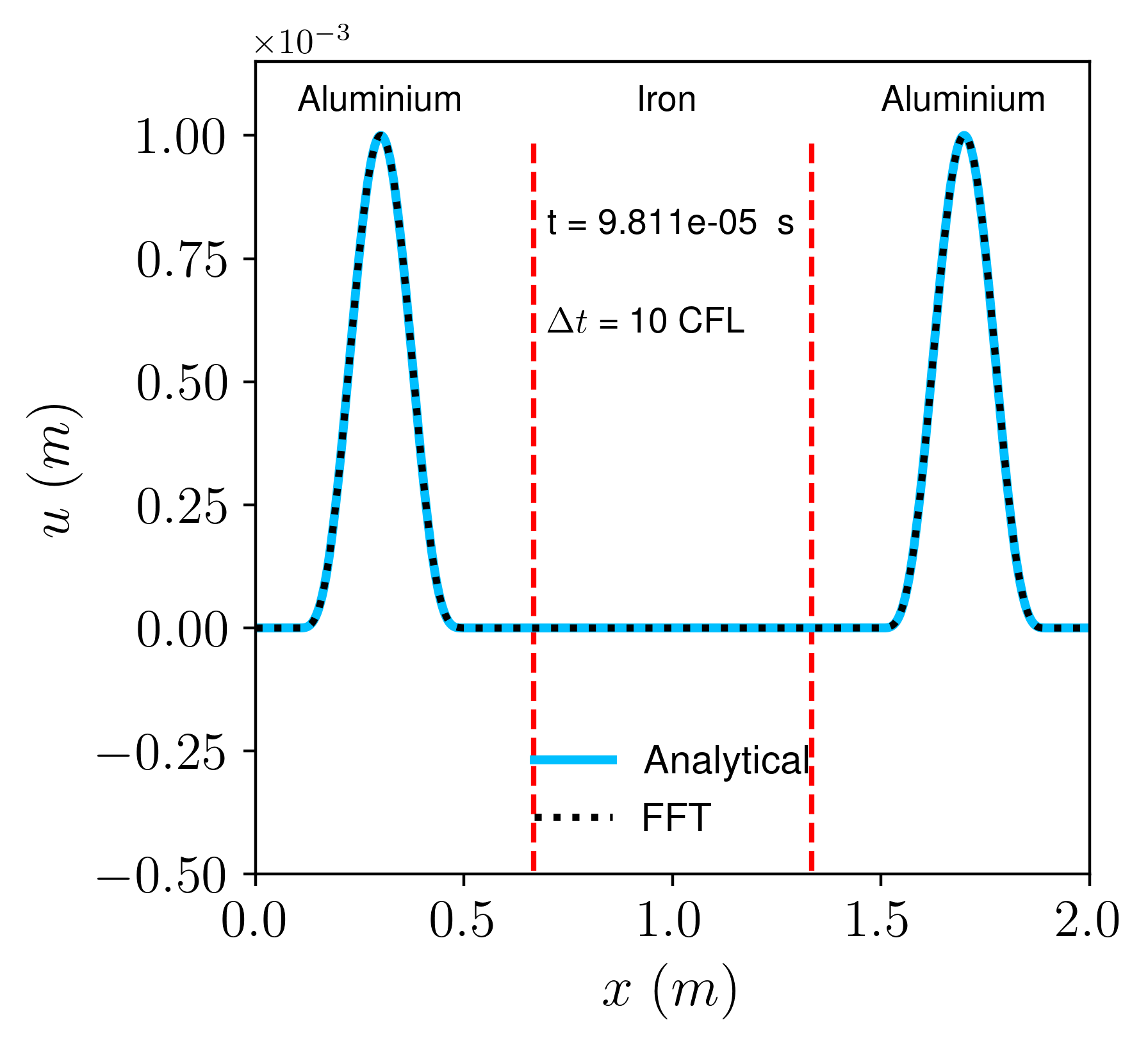}}
	\subfloat{
		\label{fig12a}
		\includegraphics[width=0.33\textwidth,trim=0 0cm 0 0cm,clip]{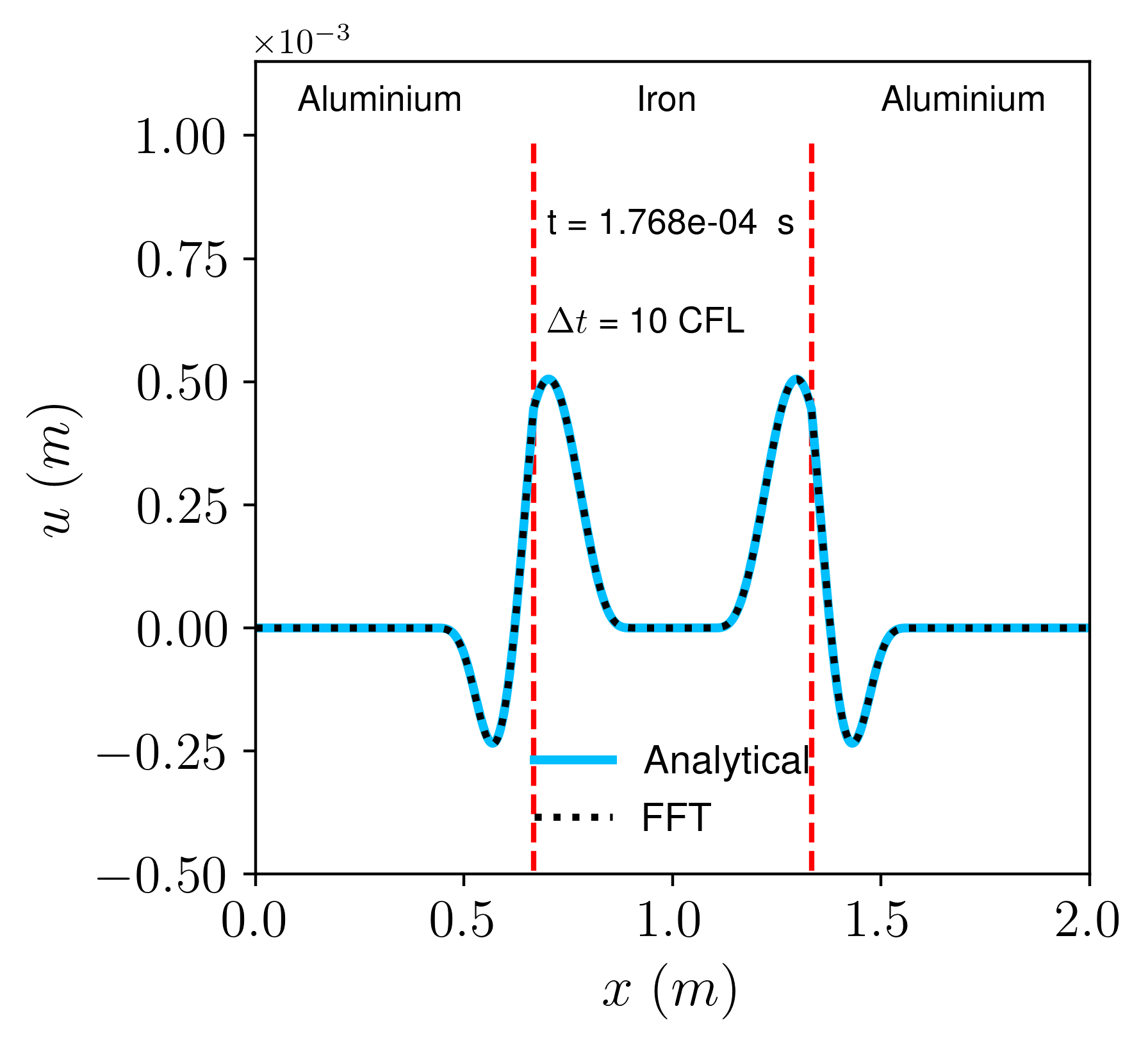}}
	\subfloat{
		\label{fig13a}
		\includegraphics[width=0.33\textwidth,trim=0 0cm 0 0cm,clip]{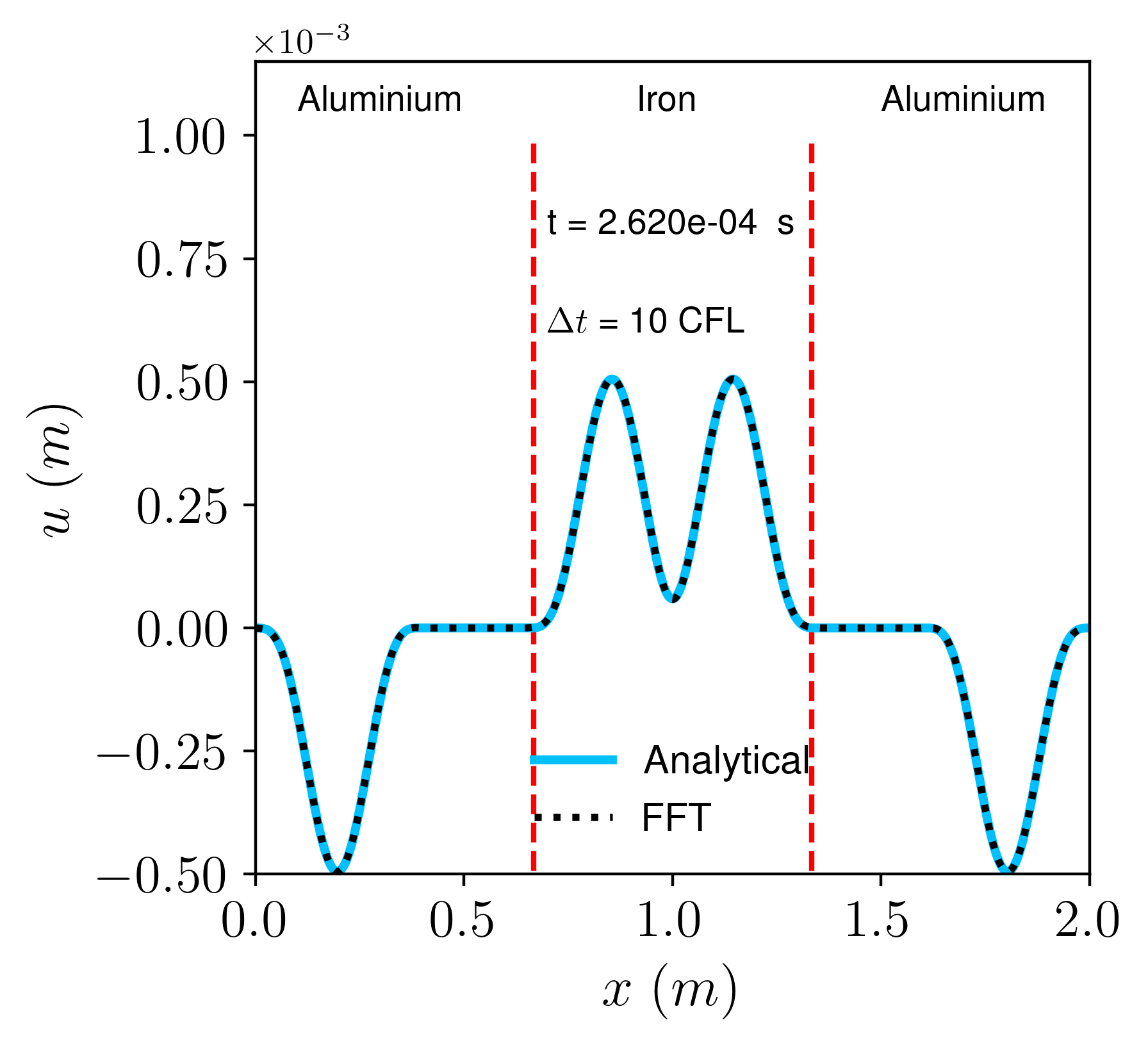}} \\
	\subfloat{
		\label{a}
		\includegraphics[width=0.33\textwidth,trim=0 0cm 0 0cm,clip]{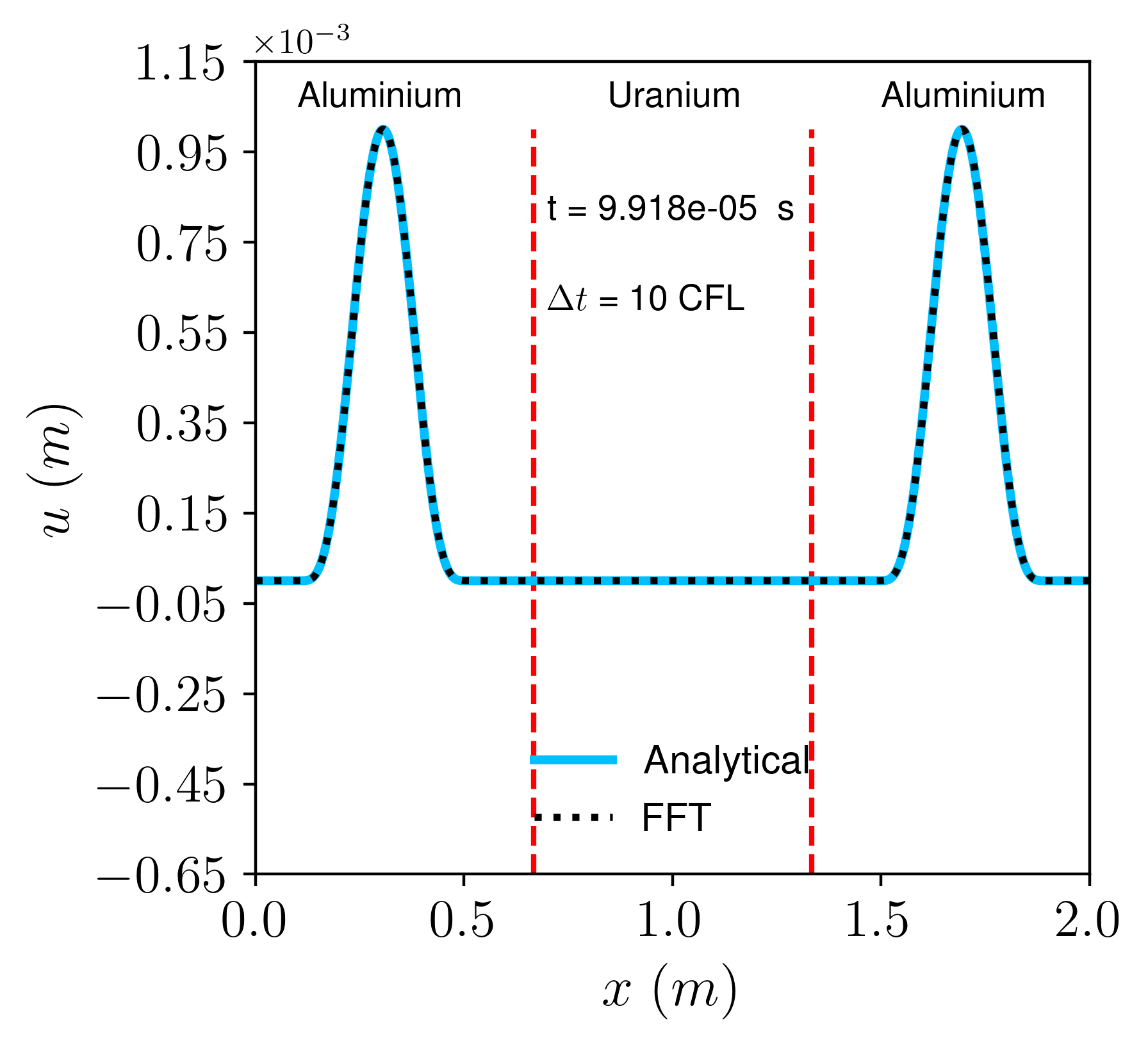}}
	\subfloat{
		\label{fig12a}
		\includegraphics[width=0.33\textwidth,trim=0 0cm 0 0cm,clip]{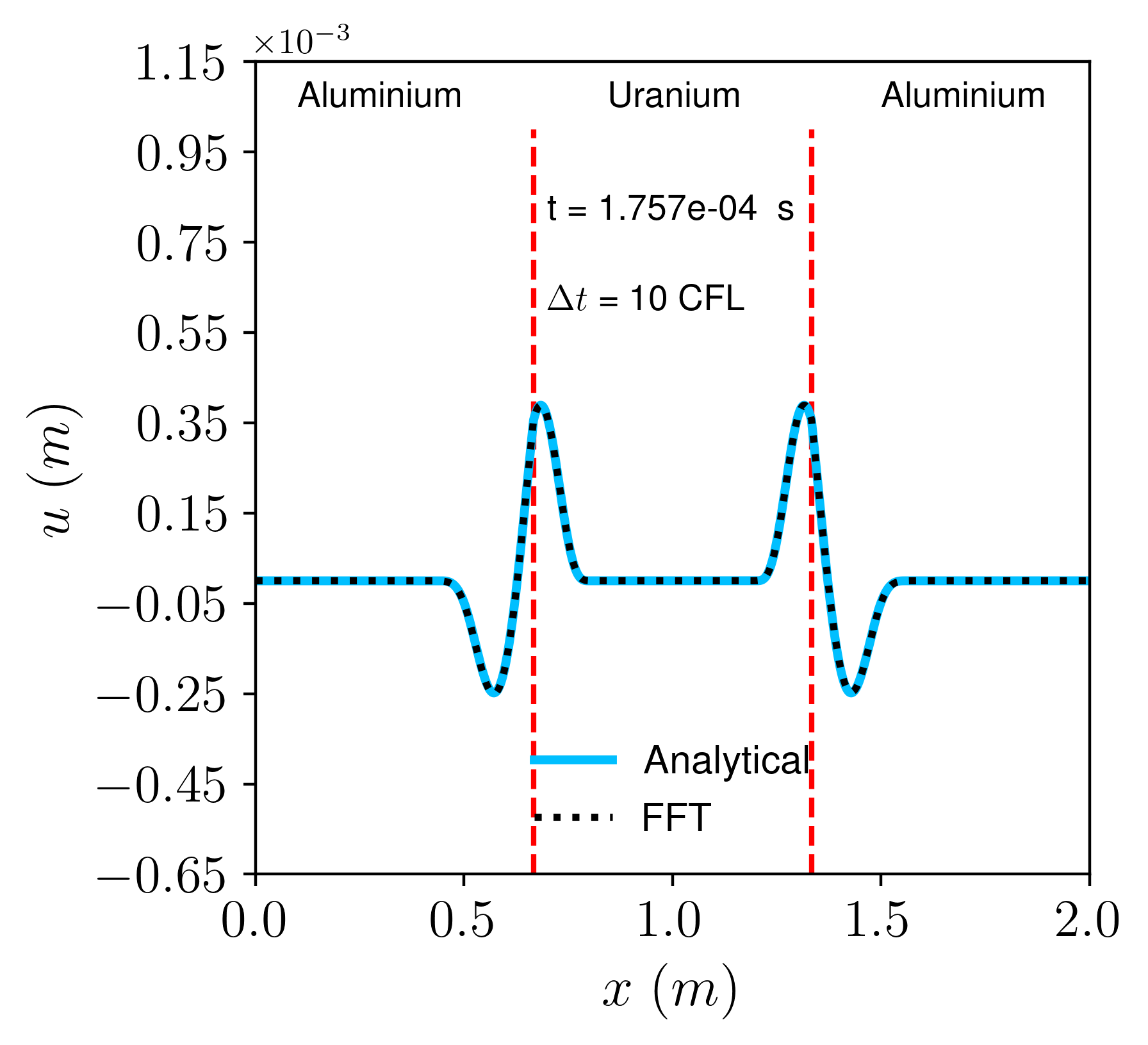}}
	\subfloat{
		\label{fig13a}
		\includegraphics[width=0.33\textwidth,trim=0 0cm 0 0cm,clip]{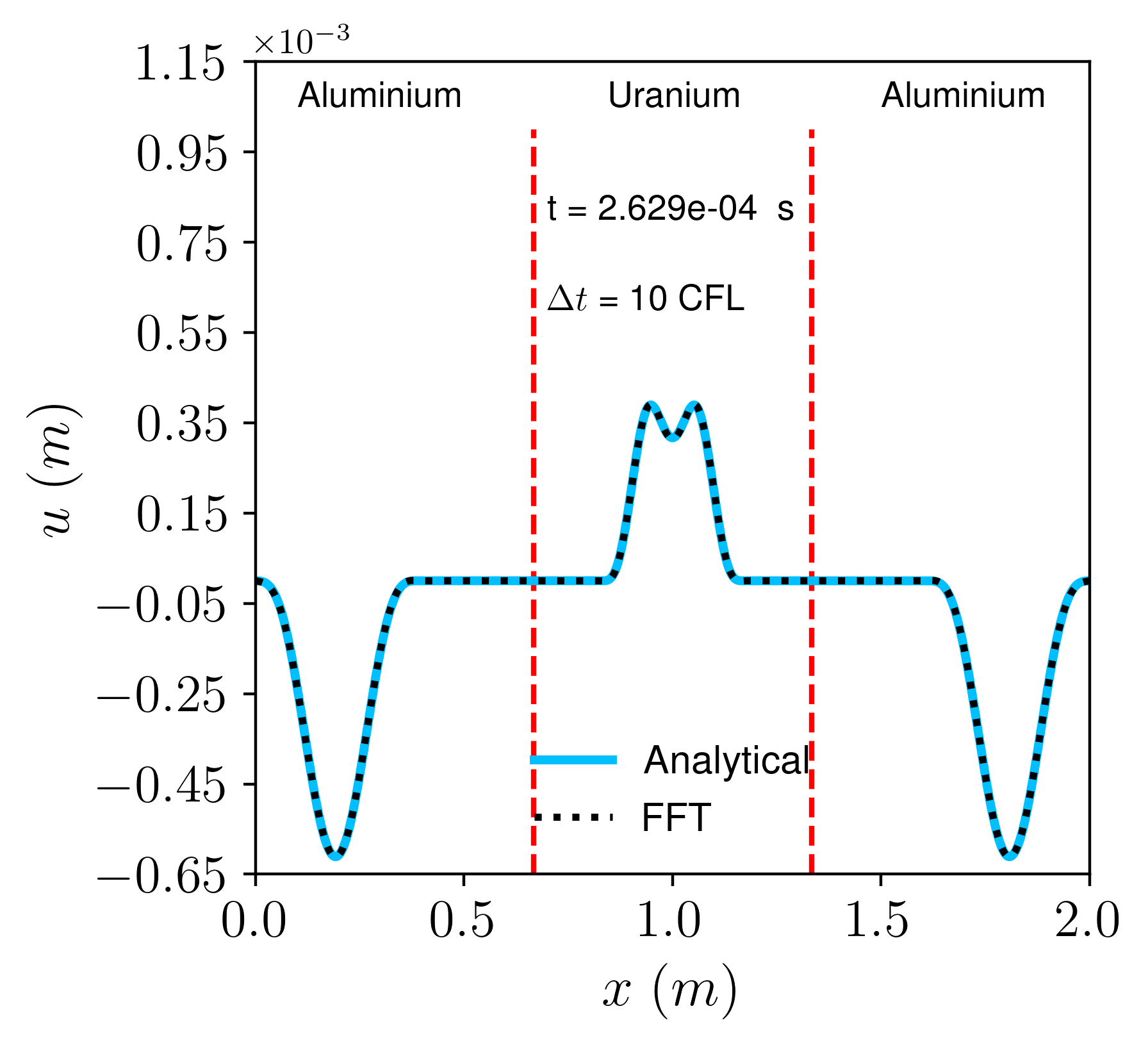}} \\
	\caption{1D wave propagation (displacement field) along the layered domain. N = 6561 voxels}
	\label{fig_1D-layered medium}    
\end{figure}

\begin{table}[h]
	\renewcommand*{\arraystretch}{1.5}
	\centering
	\caption{Accuracy and performance of simulation results for the layered medium. 1 core computation}
	\label{tab: 1D-propagation_hetero}      
	\begin{tabular}{ccccccccc}
	     & \multicolumn{4}{ c }{Aluminium-Iron} & \multicolumn{4}{ c }{Aluminium-Uranium}
		\\\hline
		$\Delta$t/CFL & $t$ (s) \tiny{ABQ}  & $\mathrm{e_r}$ \tiny{ABQ} & $t$ (s) \tiny{FFT} &$\mathrm{e_r}$ \tiny{FFT} & $t$ (s) \tiny{ABQ}  & $\mathrm{e_r}$ \tiny{ABQ} & $t$ (s) \tiny{FFT} &$\mathrm{e_r}$ \tiny{FFT} 
		\\\hline
		1 & - & - & 33.37  & $\mathrm{3.46 \cdot 10^{-3}}$ & - & - & 44.28 & $\mathrm{4.12 \cdot 10^{-3}}$ \\
		2 & - & - & 16.68  & $\mathrm{3.40 \cdot 10^{-3}}$ & - & - & 23.22 & $\mathrm{4.00 \cdot 10^{-3}}$ \\
		5 & - & - & 7.12  & $\mathrm{3.03 \cdot 10^{-3}}$  & - & - & 10.02 & $\mathrm{3.53 \cdot 10^{-3}}$\\
		10 & 194 & $\mathrm{8.50 \cdot 10^{-3}}$  & 4.01  & $\mathrm{3.50 \cdot 10^{-3}}$ & 181 & $\mathrm{8.99 \cdot 10^{-3}}$& 5.61 & $\mathrm{3.54 \cdot 10^{-3}}$ 
		\\\hline
		\multicolumn{9}{c}{Explicit scheme} \\
		1 & 1 & $\mathrm{6.24 \cdot 10^{-2}}$ & 7.69 & $\mathrm{4.21 \cdot 10^{-3}}$ & 1 & $\mathrm{6.37 \cdot 10^{-2}}$ & 7.15 & $\mathrm{4.51 \cdot 10^{-3}}$
		\\\hline
	\end{tabular}
\end{table}

\subsubsection*{Numerical performance}
The calculation times for the two cases studied, with a discretization of $N$=6561 voxels, are represented in tables \ref{tab: 1D-propagation} and \ref{tab: 1D-propagation_hetero} respectively. First, it can be observed that the stability of the implicit solver allows one to obtain a solution as good as that for $\Delta t=CFL$ using $\Delta t=10CFL$ with almost an order of magnitude lower computational cost. The deviation from an inverse linear scaling with time increment size is small; e.g., the ratio of simulation times for $CFL$ and $\Delta t=10CFL$ for the heterogeneous simulation is 8.32 instead of 10, and its origin can be found in iterative linear solver. To solve $u$ at time $n$, the solution $u_{n-1}$ is used as the initial guess in the conjugate gradient method, and this guess is closer to the solution for shorter time steps, so the number of iterations per time step decreases slightly. 

Regarding the comparison of the computational cost of the different schemes, the time spent with the implicit FFT approach to solve this relatively small problem compared to FE with the same integrator, discretization, and time step was remarkably shorter. The FFT-based solution was 200 times faster for the homogeneous case and 50 times faster for the heterogeneous. This increase in performance will be even more clear in large three-dimensional problems, as will be presented in Section 3.3, due to the order $n \log n$ of the computational cost in FFT-based methods. Moreover, it is important to note that the FFT-based method needs larger times to solve the problem when the contrast between phases increases. This behavior is expected since the preconditioner that is used to solve the corresponding linear systems is based on the homogeneous solution. Regarding the comparison with explicit schemes, explicit FE is for this relatively small problem the fastest approach, but paying the prize of having much less accuracy. Explicit FFT is much more accurate than explicit FE  , but is also computationally more demanding.

\subsubsection*{Order of accuracy}
\begin{figure}[H]
    \centering
		\includegraphics[width=0.48\textwidth]{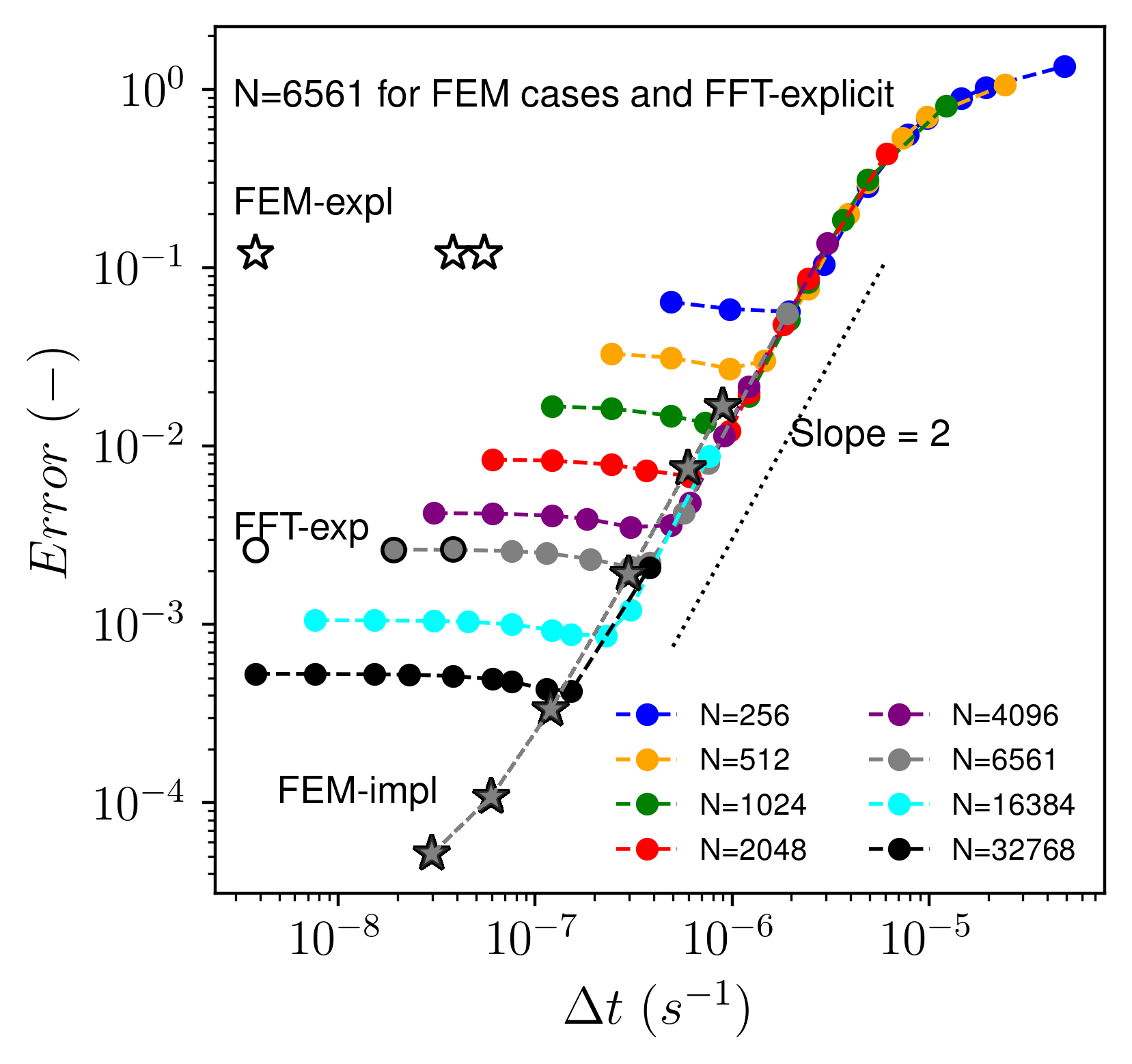} 
		\includegraphics[width=0.48\textwidth]{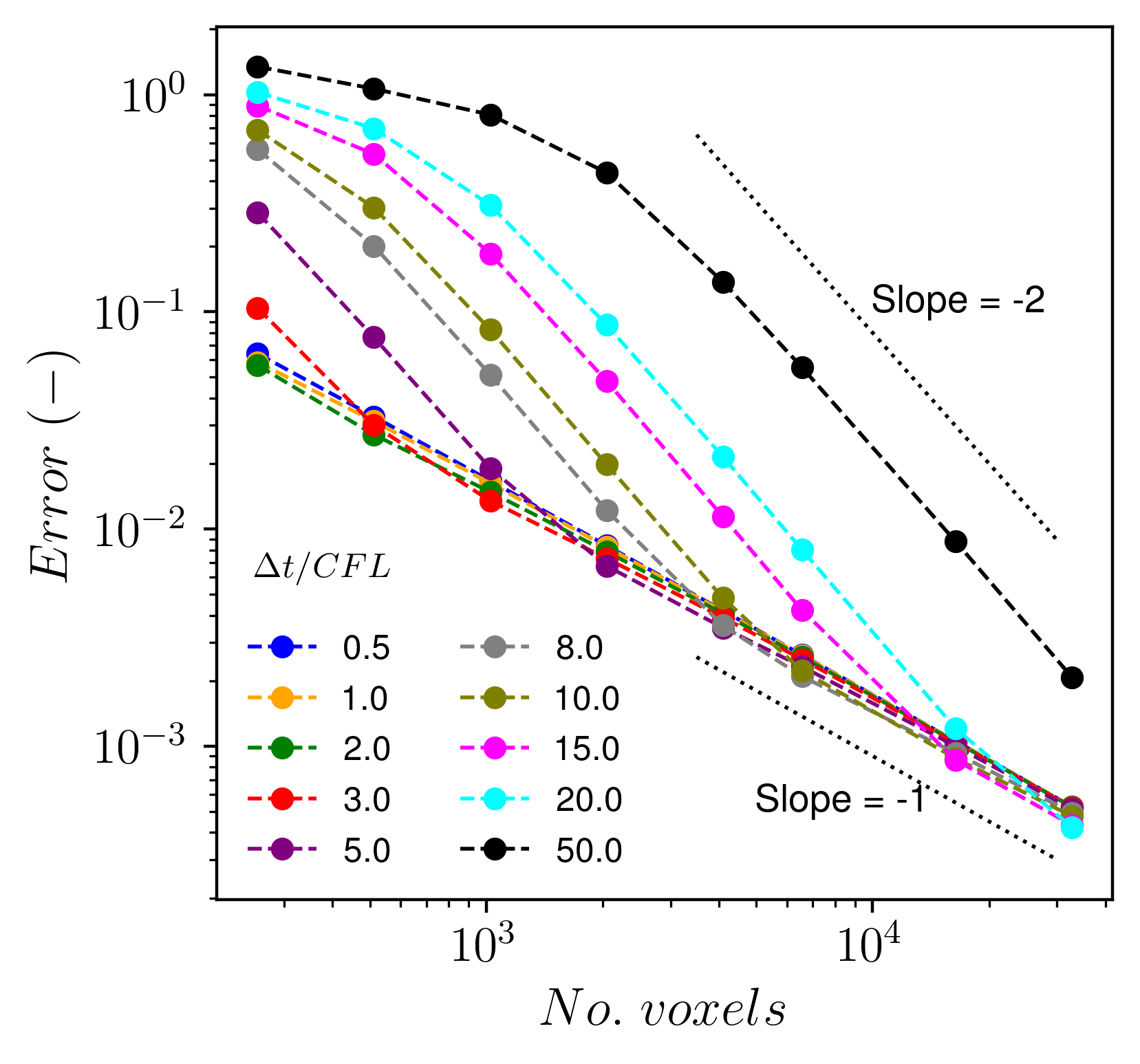}
	\caption{Order of accuracy in time (left) and space (right) of the presented method for wave propagation in 1D homogeneous elastic medium. The time convergence is compared with the FE   -implicit, -explicit and FFT-explicit method.}
	\label{fig:1D_convergence}    
\end{figure}

To extend the analysis of the methods in terms of precision, the results of homogeneous material are analyzed for all combinations of discretization of time and space. The results are presented in Fig.\ref{fig:1D_convergence}, showing  error vs time step (left) and spatial  discretization (right). The proposed method shows an order of accuracy in time of two, the same as the implicit FE solver, for time steps greater than the $CFL$ condition (second marker of each curve). However, the results also show that there is a critical time step that depends on the spatial discretization, below which the error remains constant and does not decrease as the time step decreases. This kind of behavior is common in implicit integrators \cite{Marfurt84}, and a similar trend can be inferred from the two lowest points of the FE  -implicit results ($\Delta t / CFL = 0.5, \ 1$). However, the implicit FE solver critical time is shorter than the FFT one, so it is capable of achieving more accurate results for time steps on the order of $CFL$ number. Looking at the results of the explicit schemes, it can be observed that the error keeps also constant when decreasing the time step below the $CFL$ number (larger time increment). It is interesting to note that the accuracy of explicit FFT for time steps below $CFL$ is similar to the implicit version and one order of magnitude higher in the FE solver. 
Regarding the order of convergence in space, it is observed that the FFT-implicit method converges linearly for times steps similar to the $CFL$ condition and presents a quadratic convergence for times steps larger than ten times the $CFL$ condition. Therefore, from these results, it is clear that the proposed implicit FFT method is optimal for problems which involve fine spatial discretization in space because it allows the use of large time steps preserving accuracy and efficiency. As will be presented in the forthcoming examples, wave propagation in heterogeneous microstructures that require fine discretizations are ideal systems for this approach.

\subsection{Three dimensional problem} \label{sec:3D-problem}
In this section, we consider the propagation of plane waves in a periodic 3D medium, represented by the domain $\Omega$. A uniform displacement (Eq. \eqref{eq-prescribed_displacement-nd}) is applied to one of the planes at the cell boundary $\Gamma=\{ \mathbf{x}  | \ x_3= 0 \}$. Two cases are considered, longitudinal waves, when the applied displacement is normal to $\Gamma$ and shear waves (distortional) when the prescribed displacement is parallel to that plane.

Note that, because $\Gamma$ corresponds to a plane and due to periodicity, the solution is equivalent to a 1D case under plane strain. Therefore, the analytical solution to the problem is simply given by the 1D-wave solution (see Appendix \ref{anex:1Dsol}) in which the velocity constant $c_0$ has now a different definition depending on whether the wave is longitudinal or distortional (Eq. \ref{eq-wave_velocity-nd}). 
\begin{equation}
\label{eq-wave_velocity-nd}
\mathbf c_0 = 
\begin{cases}
\sqrt{\dfrac{\lambda + 2\mu}{\rho}};  & \text{for longitudinal waves}\\
\sqrt{\dfrac{\mu}{\rho}};  & \text{for shear waves} \\
\end{cases}
\end{equation}

The periodic domain used $\Omega$, is defined as $\Omega =\{ 0\leq x_1 \leq L_1 \ \times \ 0\leq x_2 \leq L_2 \ \times \ 0\leq x_3 \leq L_3 \}$ with $L_1 = 0.1 \ m$, $L_2 = 1 \ m$,  $L_3 = 2 \ m$ and is discretized with $N_1 = 1$, $N_2 = 27$ and $N_3 = 2187$ points.  The bell-shaped pulse is introduced prescribing a displacement on $\Gamma$ given by Eq. \eqref{eq-prescribed_displacement-nd}.
\begin{equation}
\label{eq-prescribed_displacement-nd}
\mathbf{u}(x_3=0,t) = 
\begin{cases}
u_1 = 0; \  u_2 = 0;  \ u_3 = U (t);  & \text{for longitudinal waves}\\
u_1 = 0; \  u_2 = U (t); \  u_3 = 0;  & \text{for shear waves} \\
\text{being $U(t)$ defined as in the 1D problem, Eq. \eqref{eq-prescribed_displacement}}
\end{cases}
\end{equation}
The simulation final time was set equal to the time required by a longitudinal wave to travel the length of $L_3 = 2 \ m$. As in the one-dimensional problem, two cases are considered, a homogeneous material and a layered medium. The materials used are the same as in the one-dimensional case, Al and Fe. A time step of time step of $\Delta t=10CFL$ is considered, and therefore only the implicit solvers are considered in this analysis.

In the case of FE simulations, trilinear 3D continuum elements with reduced integration (C3D8R in Abaqus) are used. The mesh had the same number of elements as the number of voxels in the FFT-based model. For the implicit FE solver, periodic boundary conditions along the three axes are introduced using multipoint constraints and the prescribed displacement is introduced as an applied displacement on the boundary $x_3=0$.

\begin{figure}
	\captionsetup[subfigure]{labelformat=empty}
	\centering
	\subfloat{
		\includegraphics[width=0.4\textwidth,trim=0 0cm 0 0cm,clip]{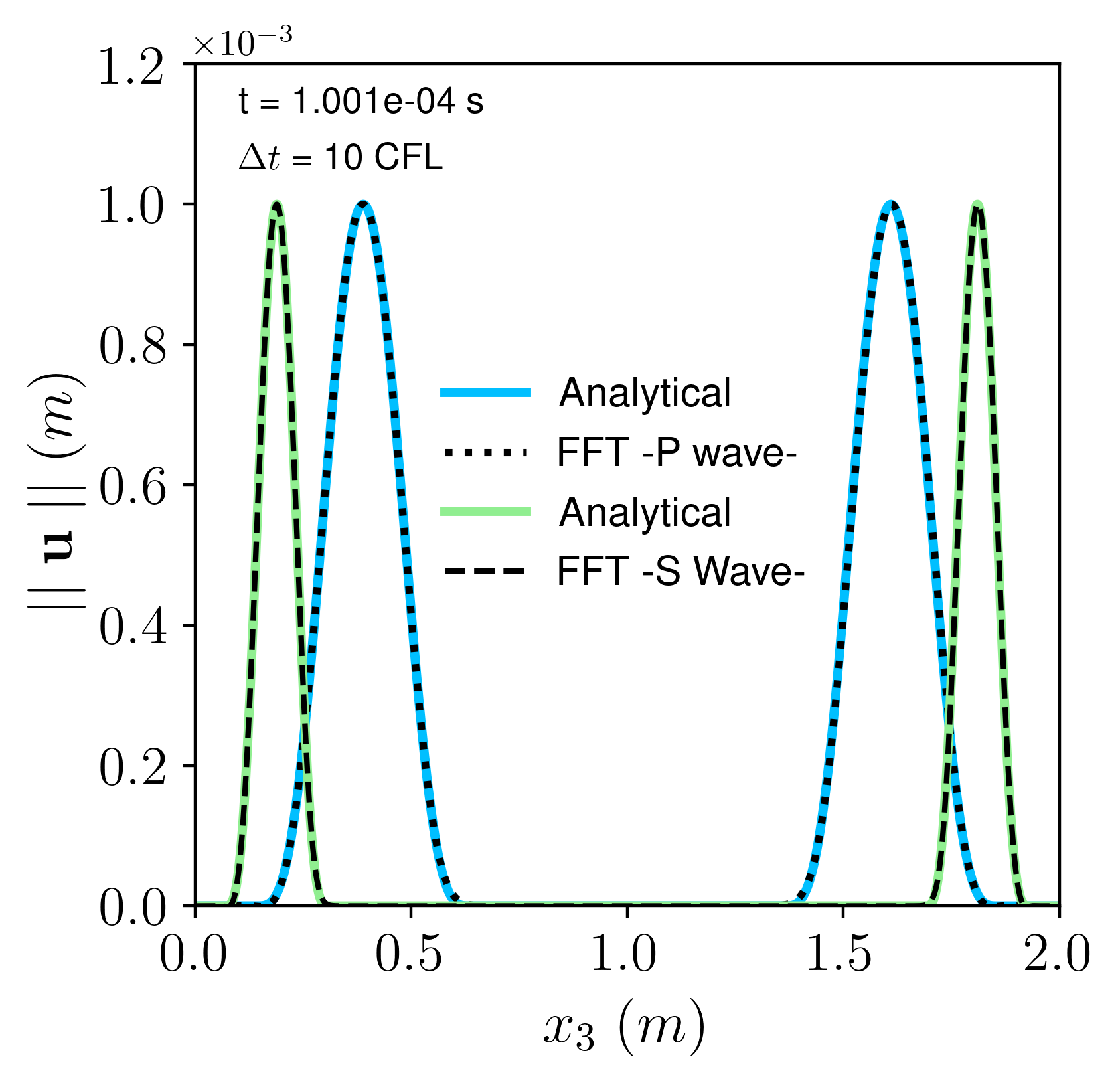}}
	\subfloat{
		\includegraphics[width=0.42\textwidth,trim=0 0cm 0 0cm,clip]{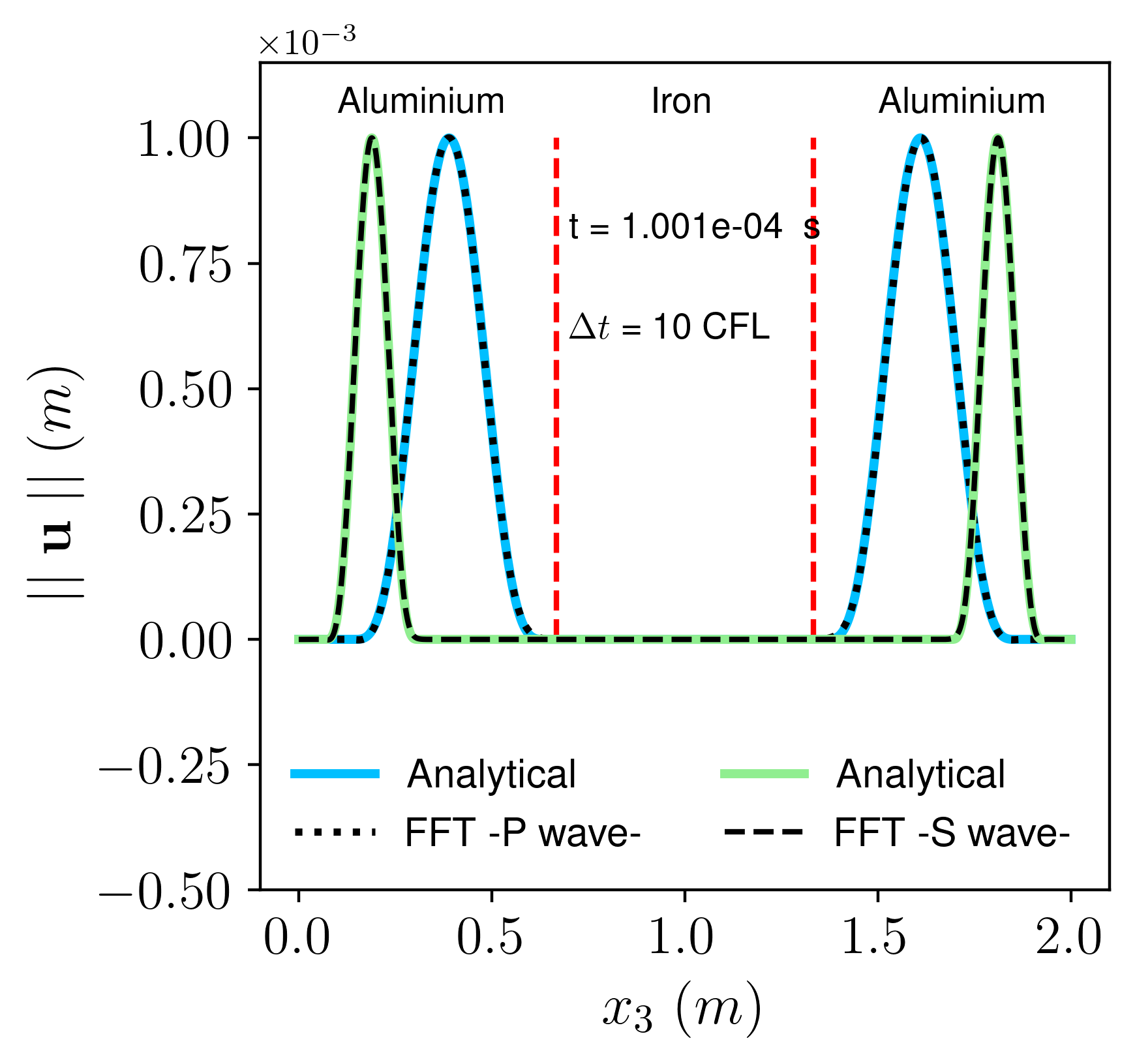}}\\
	\subfloat{
		\includegraphics[width=0.4\textwidth,trim=0 0cm 0 0cm,clip]{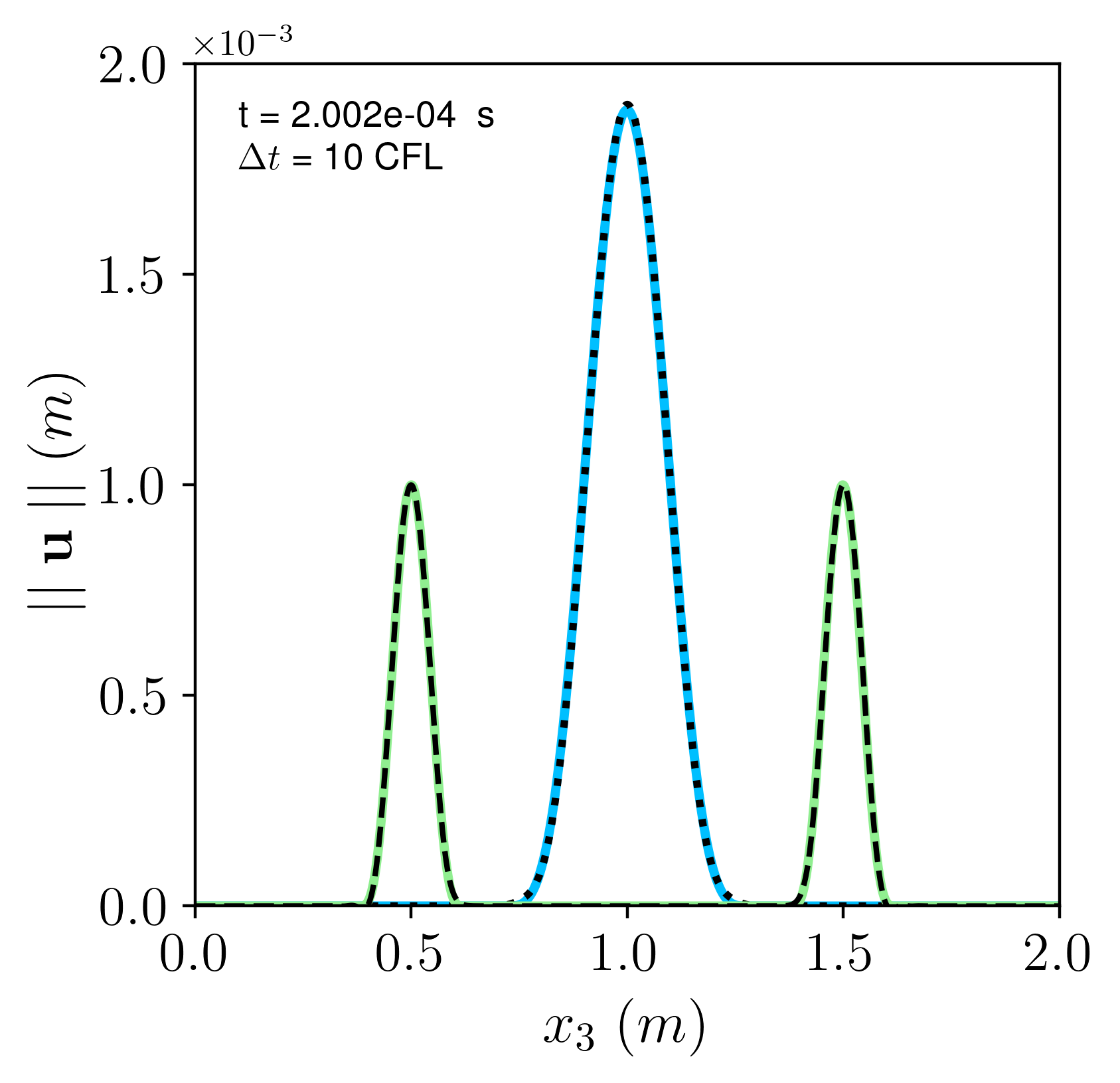}}
	\subfloat{
		\includegraphics[width=0.42\textwidth,trim=0 0cm 0 0cm,clip]{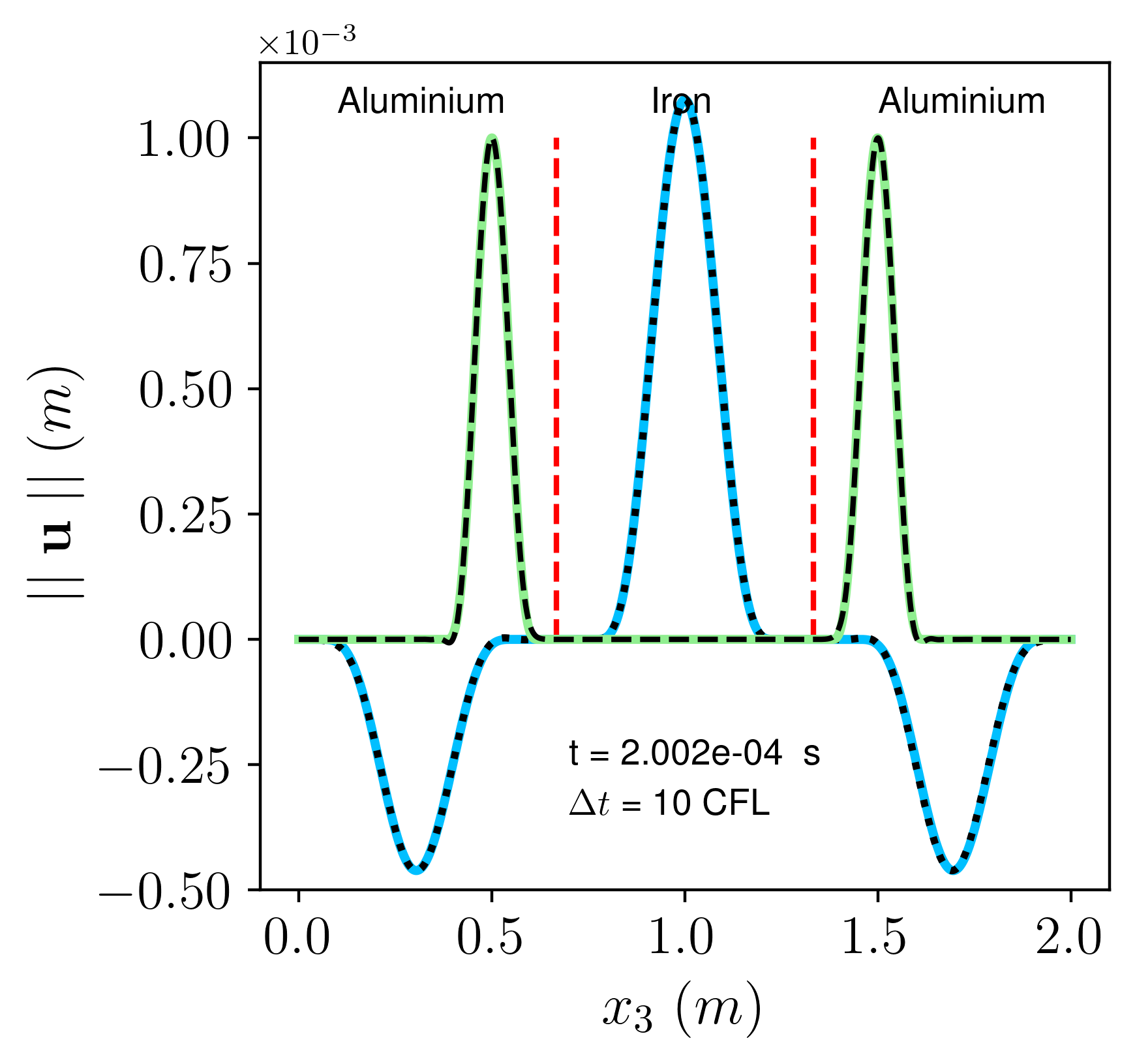}}\\
	\subfloat[(a) Homogeneous medium]{
		\includegraphics[width=0.4\textwidth,trim=0 0cm 0 0cm,clip]{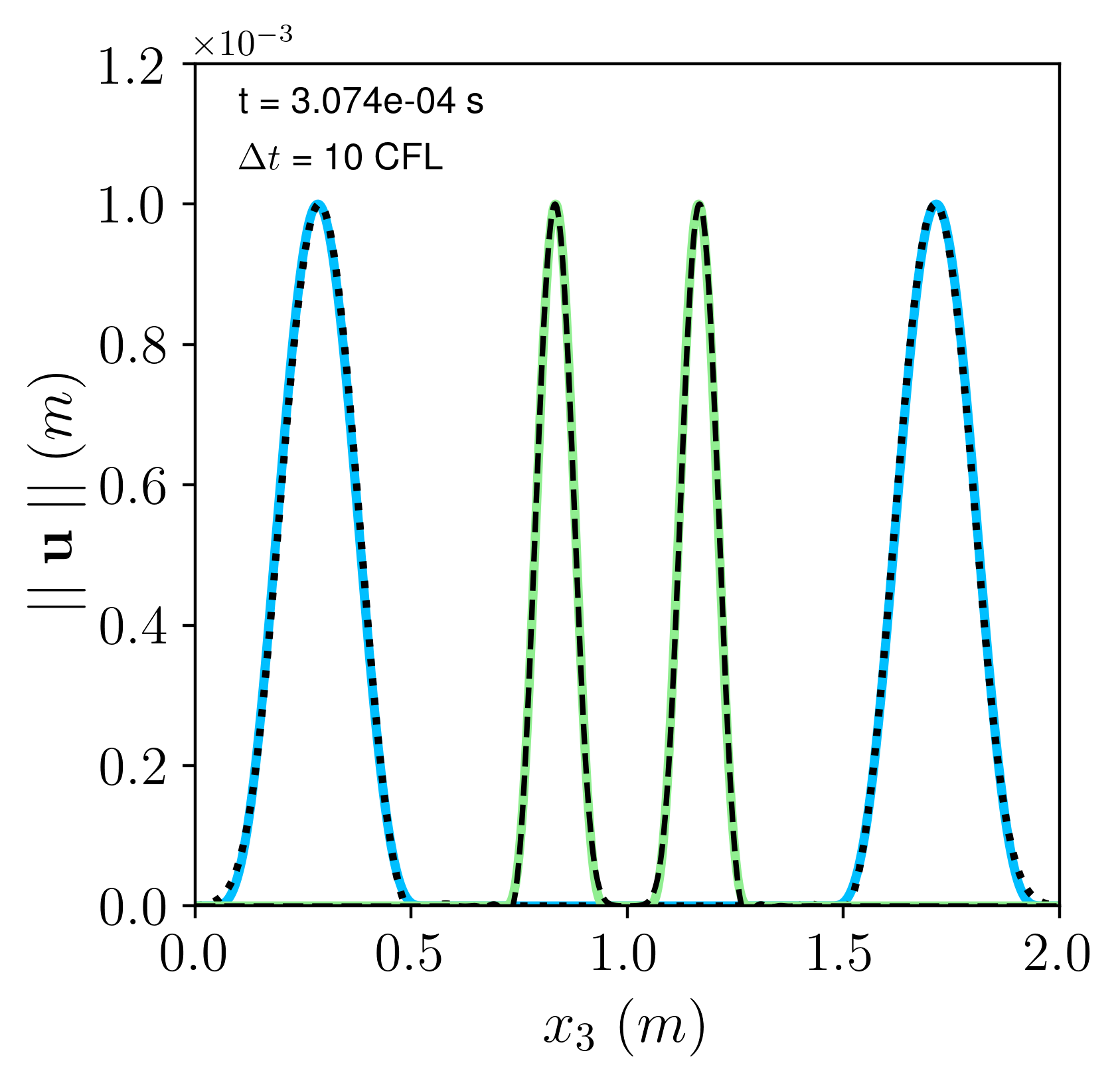}}
	\subfloat[(b) Heterogeneous medium]{
		\includegraphics[width=0.42\textwidth,trim=0 0cm 0 0cm,clip]{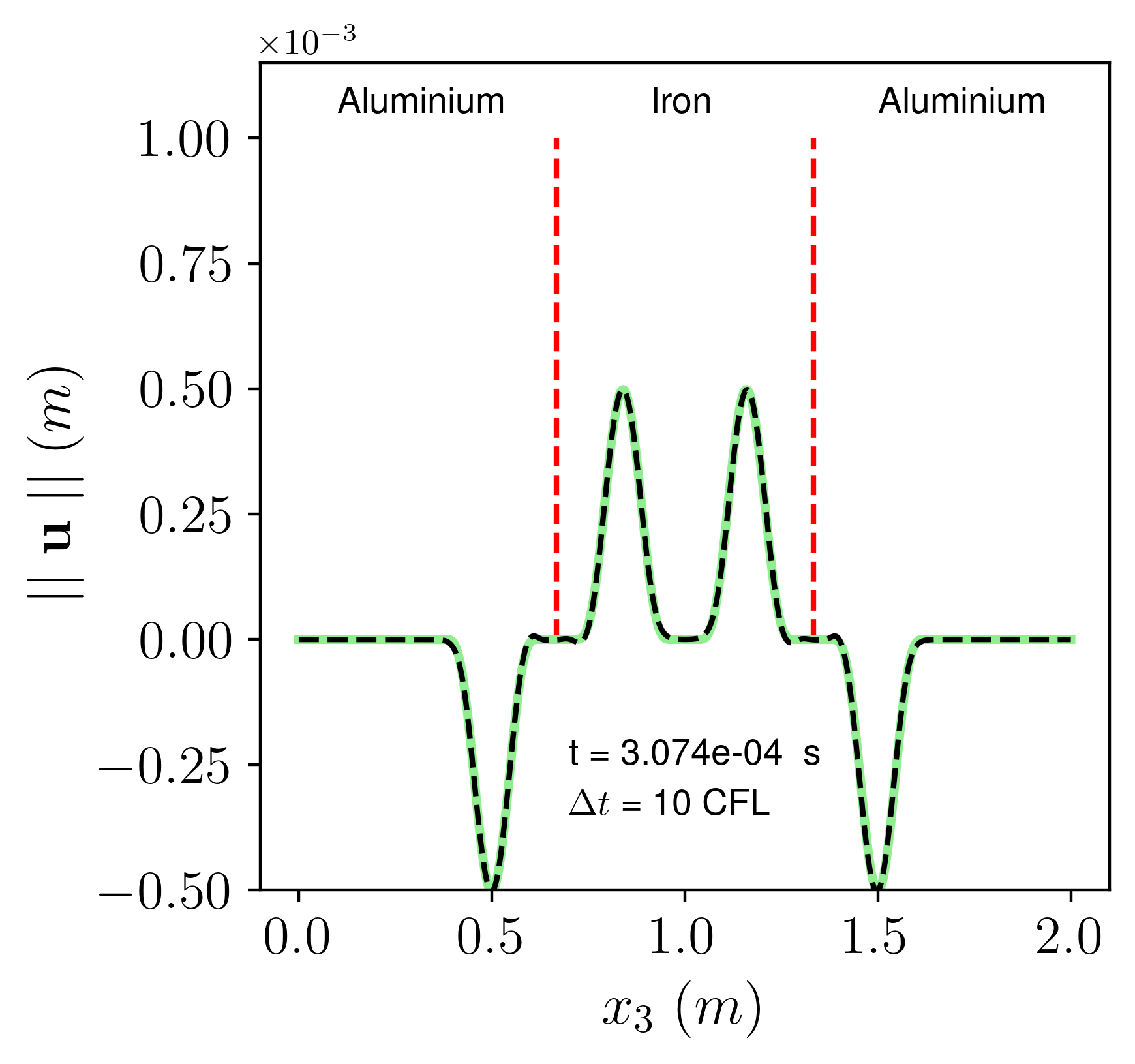}}
	\caption{1D-like wave propagation (displacement field) along a 3D Al homogeneous domain and Al-Fe layered medium.}
	\label{fig-nd_simu}       
\end{figure}

The displacement along the axis $x_3$ for different times using a is represented in Fig. \ref{fig-nd_simu} for both cases, together with the analytical solution. It is important to note that, for simplicity, the analytical solution was only computed up to the interaction of any wave with the starting or end point of the simulated domain. 

Conclusions are almost the same as for the one dimensional case. Qualitatively, analytical solution is indistinguishable from the implicit FFT result, regardless the reflections and transmissions through the heterogeneity. To quantify the accuracy, the  $L_2$-norm of the difference in the displacement field with the analytical solution has been computed, and the results are provided in table \ref{tab:nD-propagation_1} and table \ref{tab:nD-propagation_2} for the homogeneous and heterogeneous medium.

It can be seen that the errors are below $5\cdot 10^{-2}$. The difference with the analytical solution was equal to or smaller than the implicit FE results.

Another interesting result of the simulation is the wave group velocity, which can be obtained from the numerical simulations as the difference between the position of the wave at two different times divided by the time difference. The results, presented in Tables \ref{tab:nD-propagation_1} and \ref{tab:nD-propagation_2}, capture the theoretical wave speed for both modes, given in Eq. (\ref{eq-wave_velocity-nd}), with an error for the largest time increment of around $10^{-2}$.

Simulation times for both cases (homogeneous and heterogeneous medium) using FE and FFT have been summarized in tables \ref{tab:nD-propagation_1} and \ref{tab:nD-propagation_2}. Note that in the case of FFT-based simulations, the simulation time includes the preprocessing step in which the Green's functions for the nodes in $\Gamma$ are computed. In the preprocessing time of the 3D simulations, the time needed for computing the effect of each point in the manifold $\Gamma$ was 1.5 the time to perform the simulation of one time step. This is because for each point in $\Gamma$ the linear problem in Eq. \ref{eq:solve_green} has to be solved three times, one for each direction, while in a standard time step, the system is solved twice, once to obtain the value of $\mathbf{F}$ on $\Gamma$ and another to include the effect of the obtained eigenforces. In addition to this, a matrix of size $3m\times 3m$ with $m$ the number of points in $\Gamma$ have to be factorized. This last time is usually negligible compared to forming the system. As an example of the preprocessing time, in the layered case with $N_1 = 1$, $N_2 = 27$ and $N_3 = 2187$ points, the preprocessing took 45s and the rest of the simulation took  211s. The ratio of these times corresponds to 0.21, near to the time estimated of 27*3/2 knowing that the simulation took 215 increments. Furthermore, the the preprocessing of a RVE can be stored and reused for any other dynamic simulation with the same RVE.

Including preprocessing times, in the homogeneous simulation, the FFT-based approach took 14 times less than FE. In the case of the heterogeneous medium, this ratio was reduced and the FFT-based model was more than 2.5 times faster.  
If the preprocessing time is subtracted from the total time, the FFT-based simulation time becomes approximately 1/3 of the FE time. Although these times illustrate well the benefit of the proposed approach, they are just obtained for a particular model size which is indeed quite small. For a deeper analysis of the numerical performance, a more detailed analysis involving large 3D problems (millions of voxels) will be made in the next section.

\begin{table}[H]
	\renewcommand*{\arraystretch}{1.5}
	\centering
	\caption{Accuracy and performance of the FFT-based simulations for the 3D homogeneous medium. 6-core computation. $\mathrm{t_{pre}}$ refers to the preprocessing time to compute the Green's functions.}
	\label{tab:nD-propagation_1}      
	\begin{tabular}{cccccccc}
		\hline
		$\Delta$t/CFL & $t$ (s) \tiny{ABQ}  & $\mathrm{e_r}$ \tiny{ABQ} & $t_{tot}$ (s) \tiny{FFT} & $t_{pre}$ (s) \tiny{FFT} &$\mathrm{e_r}$ \tiny{FFT} & $c_{FFT}$ m/s & $c_{real}$ m/s
		\\\hline
		10 (P wave) & 552 & $\mathrm{4.04} \cdot 10^{-2}$ & 40.0 & 4 & $\mathrm{2.47} \cdot 10^{-2}$ & 6398.4 & 6395.5\\
		10 (S wave) & 481 & $\mathrm{3.34} \cdot 10^{-2}$ & 41.2 & 4 & $\mathrm{2.80} \cdot 10^{-2}$ & 3105.4 & 3111.1
		\\\hline
	\end{tabular}
\end{table}

\begin{table}[h]
	\renewcommand*{\arraystretch}{1.5}
	\centering
	\caption{Accuracy and performance of the FFT-based simulations for the 3D layered medium. 6-core computation. $\mathrm{t_{pre}}$ refers to the preprocessing time to compute the Green's functions. }
	\label{tab:nD-propagation_2}      
	\begin{tabular}{cccccc}
		\hline
		$\Delta$t/CFL & $t$ (s) \tiny{ABQ}  & $\mathrm{e_r}$ \tiny{ABQ} & $t_{tot}$ (s) \tiny{FFT} & $t_{pre}$ (s) \tiny{FFT} &$\mathrm{e_r}$ \tiny{FFT}
		\\\hline
		10 shear & 637 & $\mathrm{5.29} \cdot 10^{-2}$ & 256.5 & 45 & $\mathrm{2.73} \cdot 10^{-2}$ \\
		10 lon & 680 & $\mathrm{2.96} \cdot 10^{-2}$ & 280.1 & 44 & $\mathrm{1.70} \cdot 10^{-2}$
		\\\hline
	\end{tabular}
\end{table}

\subsection{Numerical performance}

With the aim of checking the performance of the FFT scheme developed in large 3D problems, the propagation of a longitudinal wave in an Al-Fe-Al layered medium is simulated using a variable number of elements. The problem is solved using the proposed implicit FFT approach as well as with the explicit version of FFT and both implicit and explicit FE. In all cases, implicit simulations are performed using Newmark-$\beta$ and explicit ones using central differences, as explained in Section 3. The dimensions of the prismatic bar used for the simulation are 0.1 $\mathrm{x}$ 0.1 $\mathrm{x}$ 2 m and are discretized with $N_1 = 5$, $N_2= 5$, $N_3 = \left[243, \ 729, \ 2187, \ 6561, \ 19683, \ 59049, \ 177147 \right]$ voxels. The reason for choosing this problem is twofold (1) is simple, so an analytical solution is available to compute the error of each simulation, and (2) it symmetry allows us to use standard Dirichlet and Neumann boundary conditions to impose the periodicity, avoiding the use of multipoint constraints which lead to instabilities in explicit FE \cite{SADABA2019434}.These conditions correspond to displacements in the $x$ and $y$ directions on planes $\mathrm{X=0}, \ \mathrm{X_1=L_1}, \ \mathrm{X_2=0}, \ \mathrm{X_3=L_3}$. 

The elements used in explicit FE are the same as in the implicit simulations, reduced order 8 node quadrilateral elements. Regarding the time step, it is important to note that in the implicit simulations, the time increment was set equal to the $CFL$ condition of the coarsest mesh, with N=5,5,243 voxels. With this setting, the time steps for the finest models with $\mathrm{N_3}=59049$ and $\mathrm{N_3}=177147$ have been done with steps of 243 and 729 times, respectively, the CFL condition for these discretizations.   

\begin{figure}[H]
    \centering
		\includegraphics[width=0.45\textwidth]{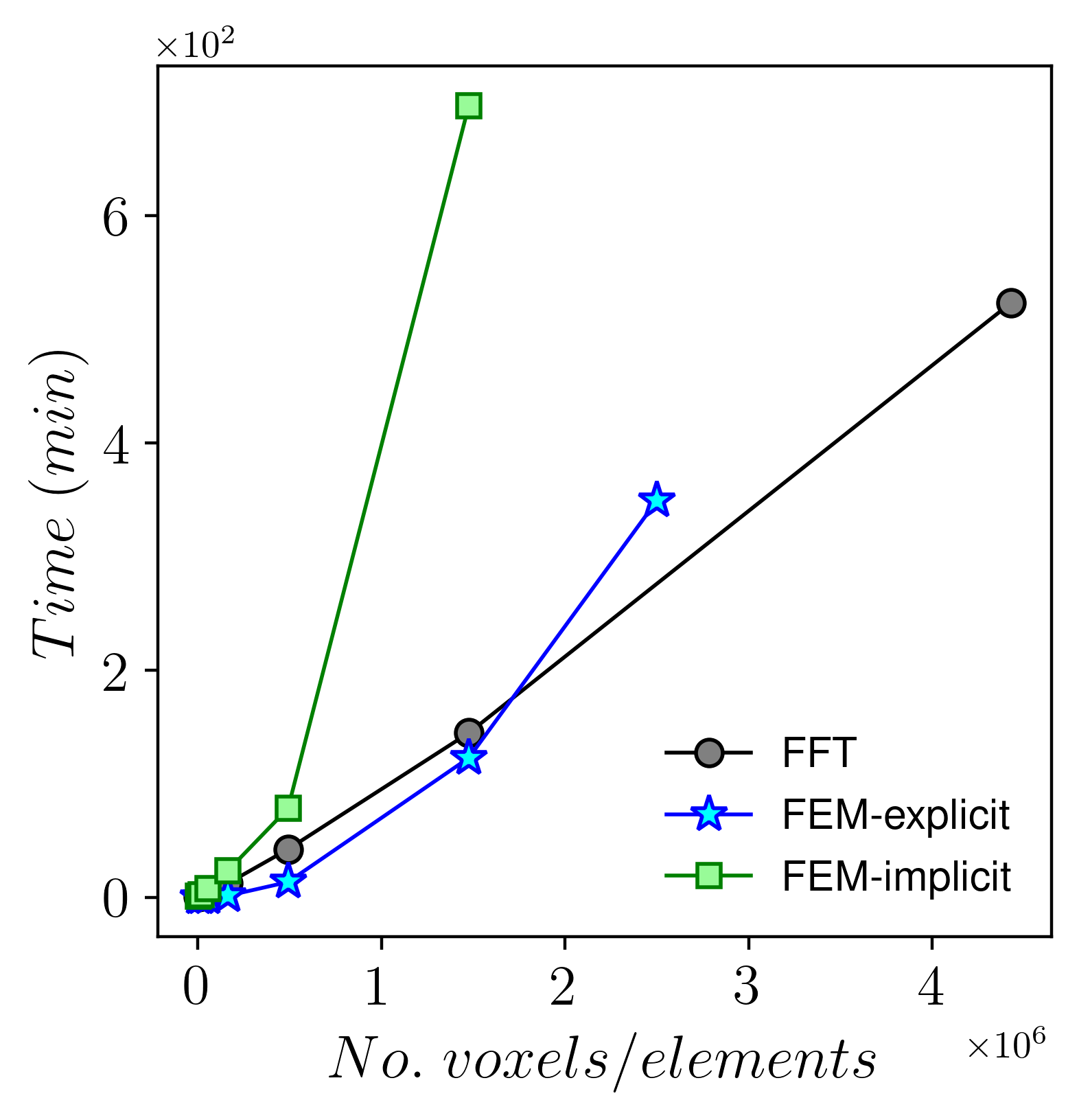} 
	\caption{Performance comparison among the FFT method, FE   explicit and FE   implicit considering the computational cost. The time increment in FE  -explicit is the one given by the $CFL$ condition, while in FFT and FE  -implicit is constant and set as the $CFL$ condition for the smallest problem.}
	\label{fig:3D_performance_analysis}    
\end{figure}

Fig. \ref{fig:3D_performance_analysis} represents the simulation times obtained as function of the number of elements/voxels. It can be observed that the implicit FFT method becomes the fastest approach for problems with more than $2 \cdot 10^{6}$ voxels. For these sizes, even explicit FE  ---which is considerably less accurate but faster--- becomes less competitive due to the small time step needed for such fine models. The order of growth of the cost of implicit FE and implicit FFT methods with the number of elements is very different, and the first one becomes orders of magnitude slower than the proposed method for a large number of voxels. Note that the conditions used for both methods are almost identical: same number of integration points and nodes than voxels, same boundary conditions, same linear solver (conjugate gradient) with same tolerance (10$^{-6}$). It has also to be noted that the efficiency of the proposed FFT approach has still place for improvement. For example, current simulations are performed using an odd number of voxels but the algorithms can be adapted for an even number of points which would improve the performance using the number of elements powers of 2.

Regarding the precision in these simulations, the error in the FFT solution was lower than for the explicit and implicit FE methods for all cases except for the two coarsest meshes, where the time step was close to $CFL$, reaching the point where the error decay loses quadratic convergence. To give some examples, for the problem with 1476225 voxels, the error of the FFT method was $9.40 \cdot 10^{-3}$, while it was equal to $1.02 \cdot 10^{-2}$ and $1.38 \cdot 10^{-2}$ for explicit and implicit FE solvers.

\section{Results}
In this section, three different problems will be solved to show the capabilities of the proposed method. The cases do not have full analytical solutions and the validation of the results is done by analyzing some physical outputs such as the resulting wave velocities, the amplitude attenuation, etc.

\subsection{Circular and spherical waves}

\subsubsection{Circular waves in a homogeneous plate}
This example simulates the deformation of a plate caused by the application of an excitation at its center. The material is Al. The domain is a square plate of dimensions $0.001\times1\times1$m, discretized in $1\times1025\times1025$ voxels, resulting in a model with more than 1 million voxels. Due to the periodicity and the use of a single voxel in $x_1$, this 3D case is equivalent to a 2D plane strain simulation. The time increment is $\Delta t = 10 \ CFL = 1.9 \cdot 10^{-6} s$ and the final time was $T = 2.86 \cdot 10^{-4} \ s$. The displacement is imposed in the center of the plate, in a region with radius $R=2.5$mm, suing the pulse of amplitude 1mm defined in Eq. (\ref{eq-prescribed_displacement}). The direction of the pulse is in $x_1$, perpendicular to the plate section, causing transverse waves.

\begin{figure}[H]
    \centering
		\includegraphics[width=0.45\textwidth]{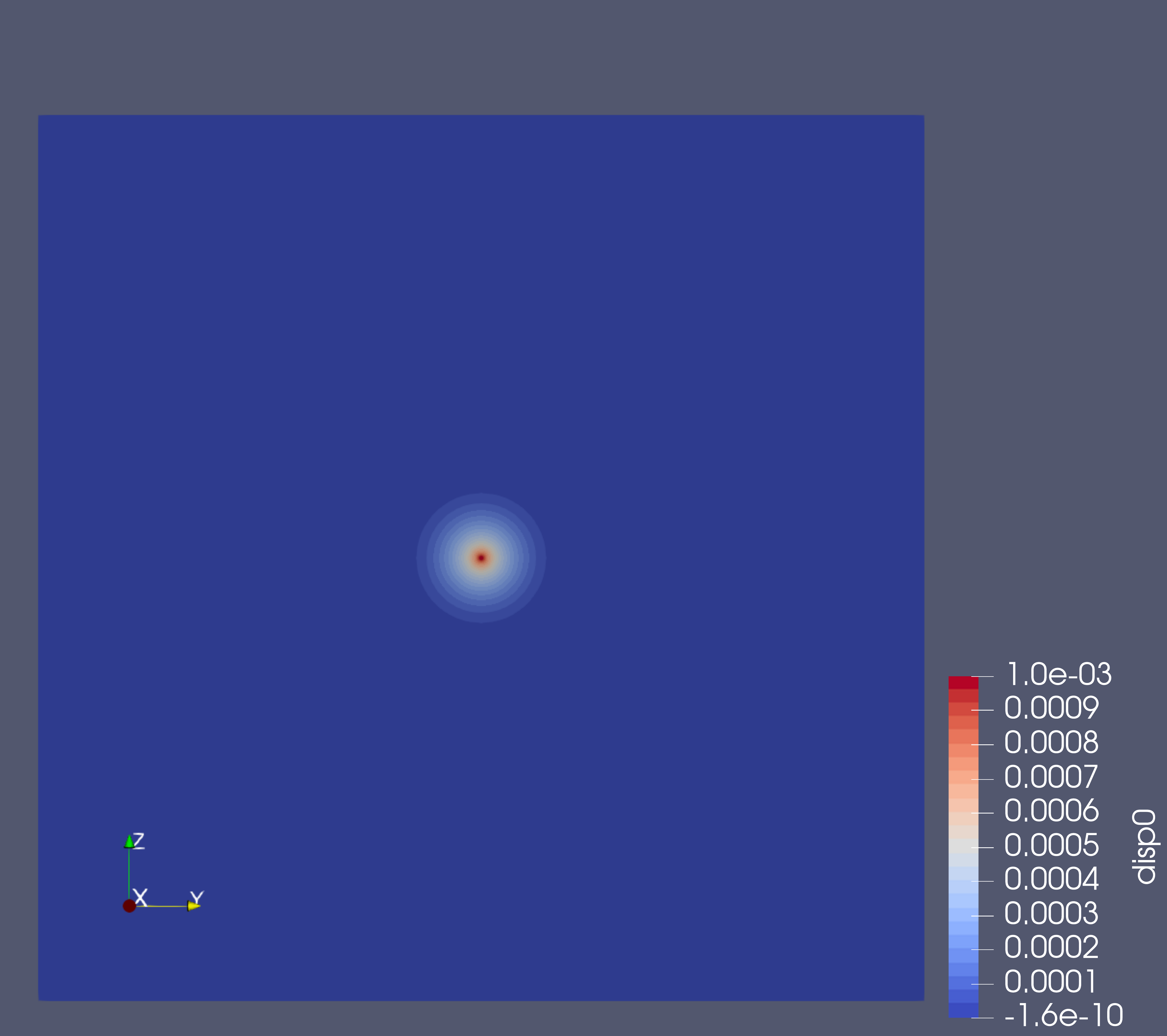} 
		\includegraphics[width=0.45\textwidth]{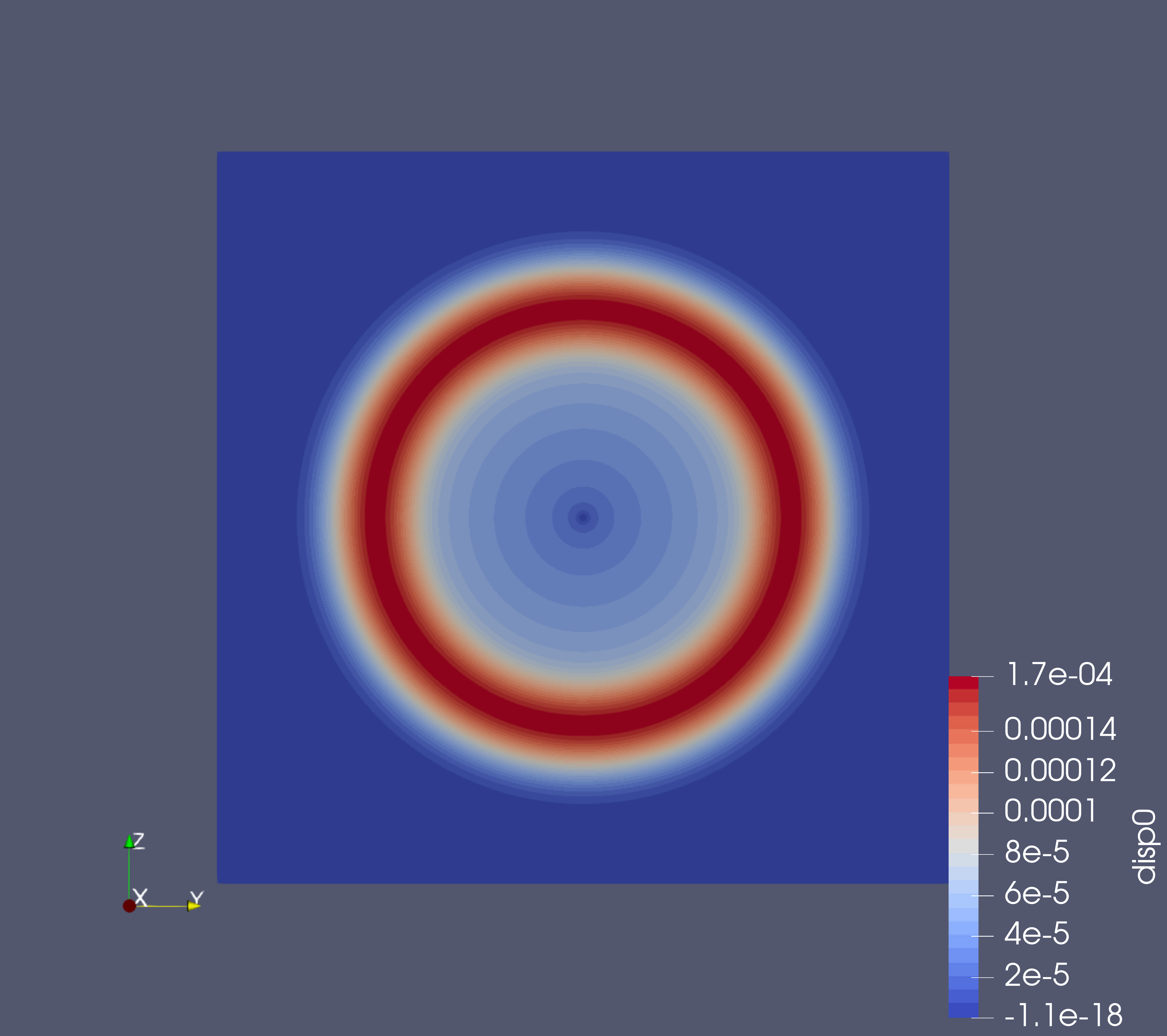}
		\includegraphics[width=0.45\textwidth]{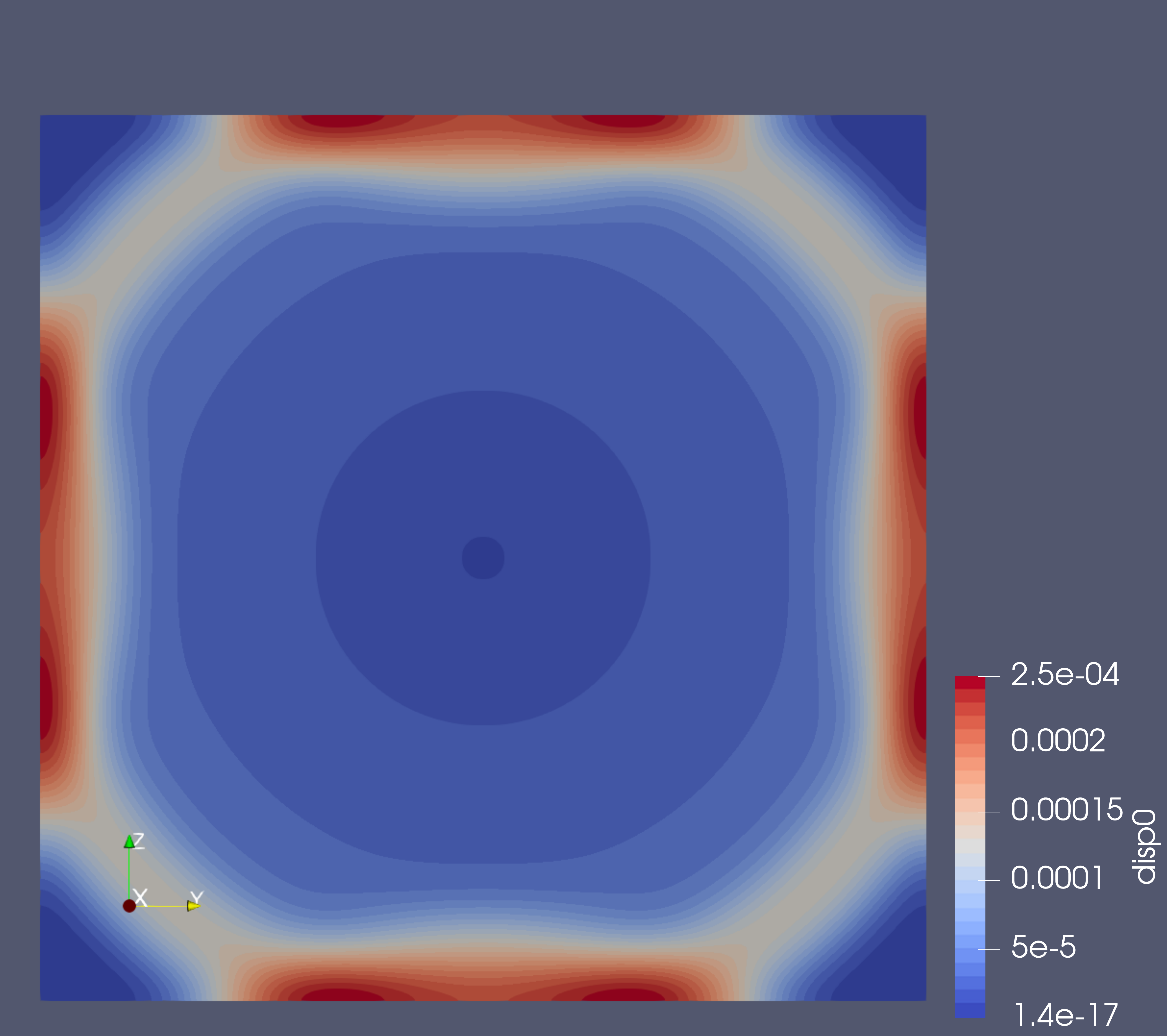} 
		\includegraphics[width=0.45\textwidth]{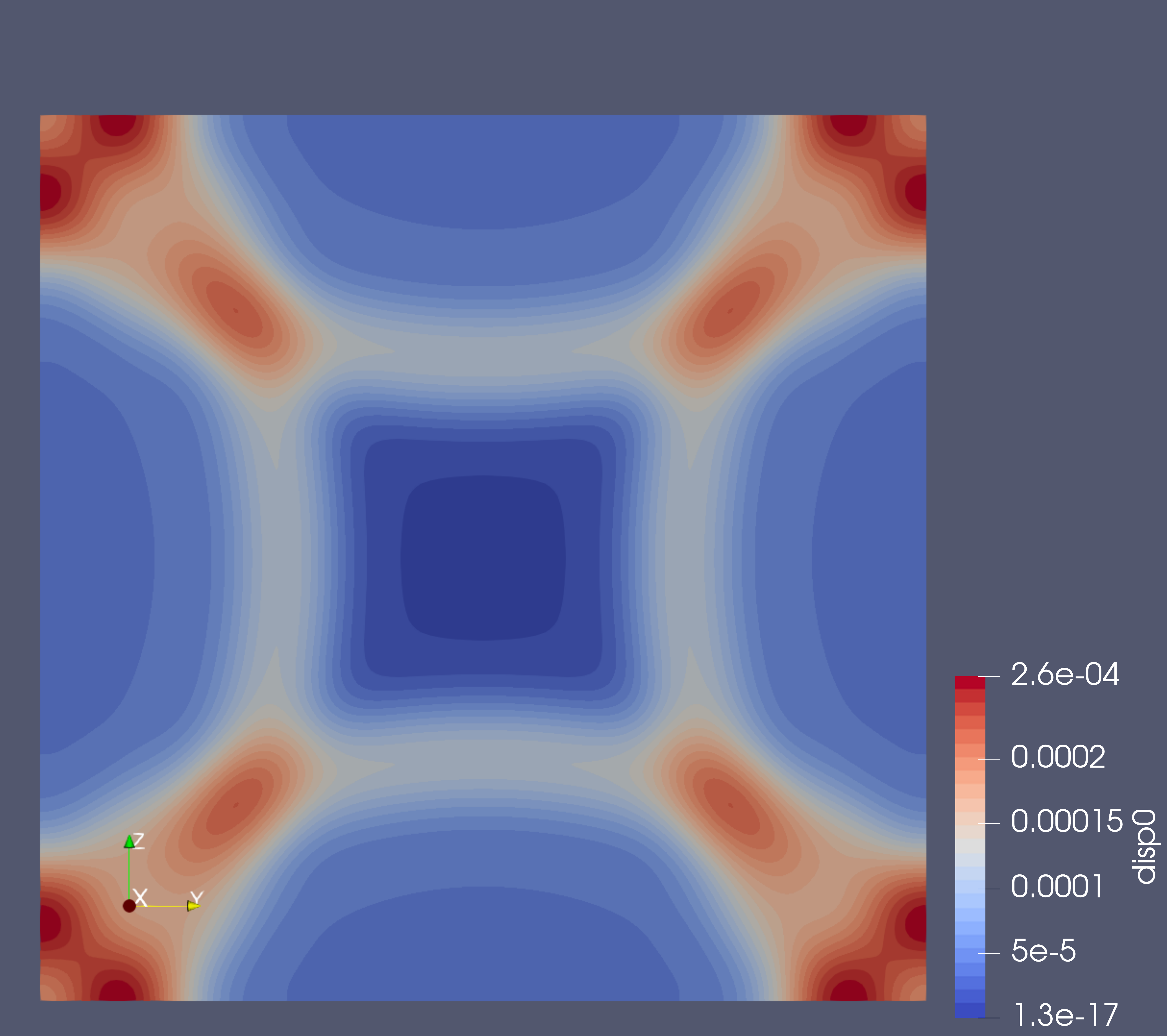}
	\caption{Transverse displacement for circular waves in plane strain traveling in Al, $t= 1.91 \cdot 10^{-5}$s, $ 1.115 \cdot 10^{-4}$s, $ 2.3 \cdot 10^{-4}$s, $ 2.86 \cdot 10^{-4}$s. Disp0 refers to the displacement in the $x$ direction.}
	\label{fig:circular_homo}    
\end{figure}

\begin{figure}[H]
    \centering
		\includegraphics[width=0.5\textwidth]{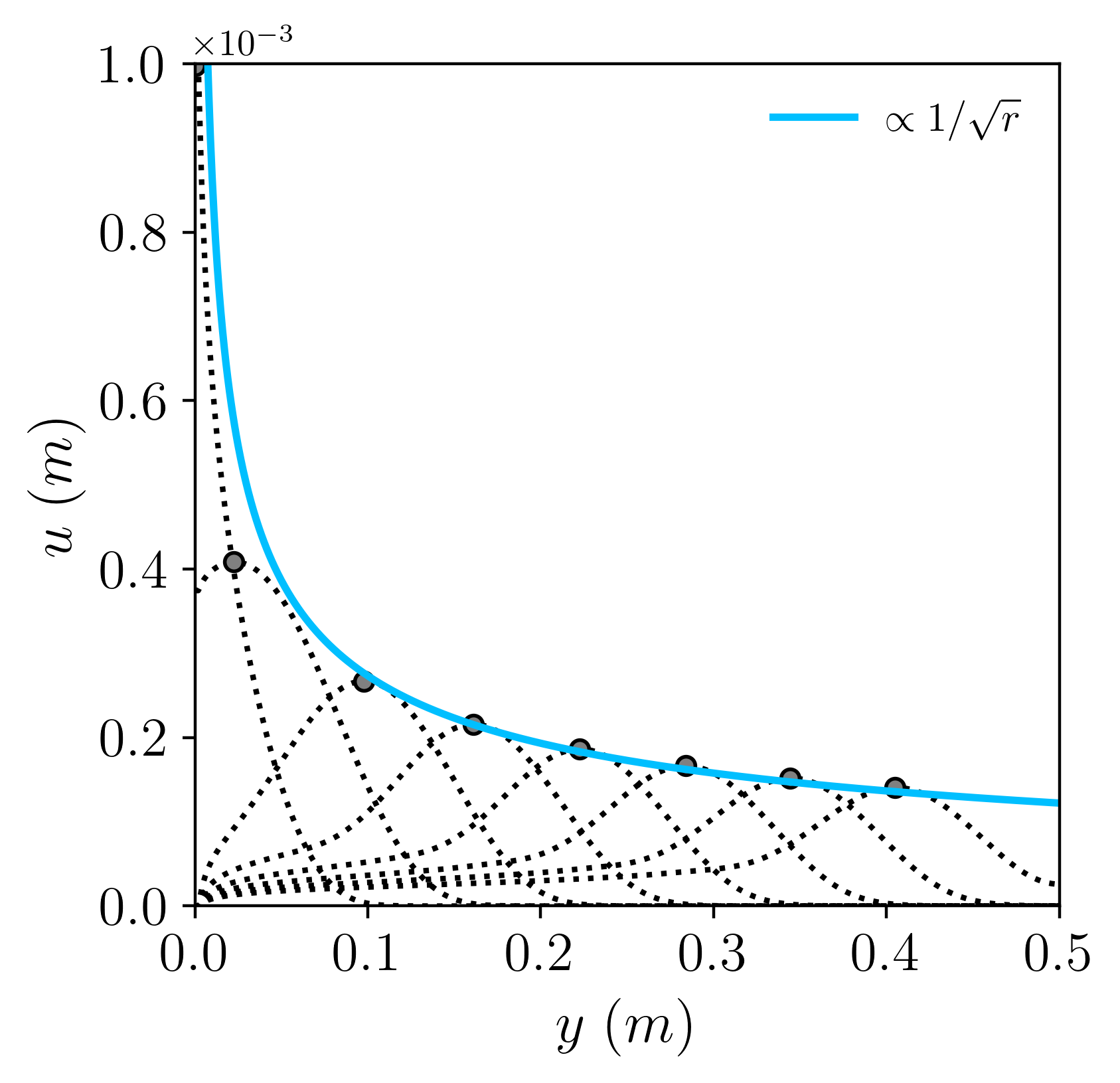} 
	\caption{Circular-wave amplitude evolution, in the homogeneous case, along the line $z = 0$ for the positive $y$ values. All times plotted are previous to the wave-front - plate-edges interaction.  The point $y=0$ corresponds to the center of the plate. $r$ value in the legend refers to the distance to the origin.}
	\label{fig:Circular_wave_decay}    
\end{figure}

The contour plots of the transverse displacement at three different times are represented in Fig. \ref{fig:circular_homo}, where the colors represent the displacement in the perpendicular direction. The simulations show how the amplitude of the traveling wave decays with the distance to the center. Note that although the maximum displacement decays, the colors in Fig. \ref{fig:circular_homo} are rescaled to the maximum displacement at each time for a clearer representation. To assess the accuracy of the simulation, we calculated the reduction in wave amplitude with time and compared it with the analytical solution. The amplitude decay is represented in Fig. \ref{fig:Circular_wave_decay}. Denoting $r$ as the distance to the center of the plate, it can be observed that the envelope of the front waves adjusts perfectly with the theoretical trend $1/\sqrt{r}$ for circular waves originating from a punctual perturbation \cite{Meyers}. 

\subsubsection{Circular waves in a framed plate}
\begin{figure}[H]
    \centering
		\includegraphics[width=0.45\textwidth]{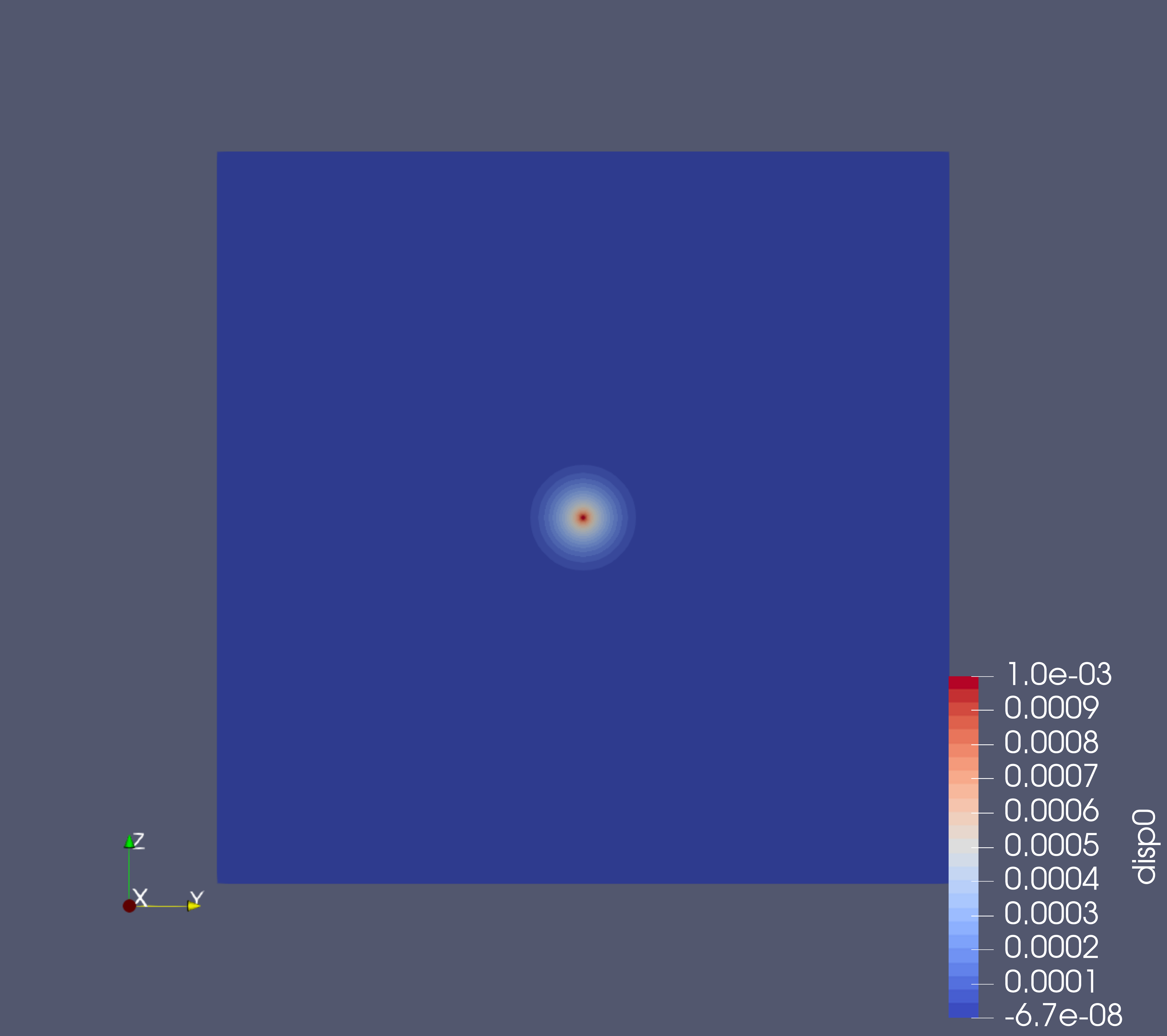} 
		\includegraphics[width=0.45\textwidth]{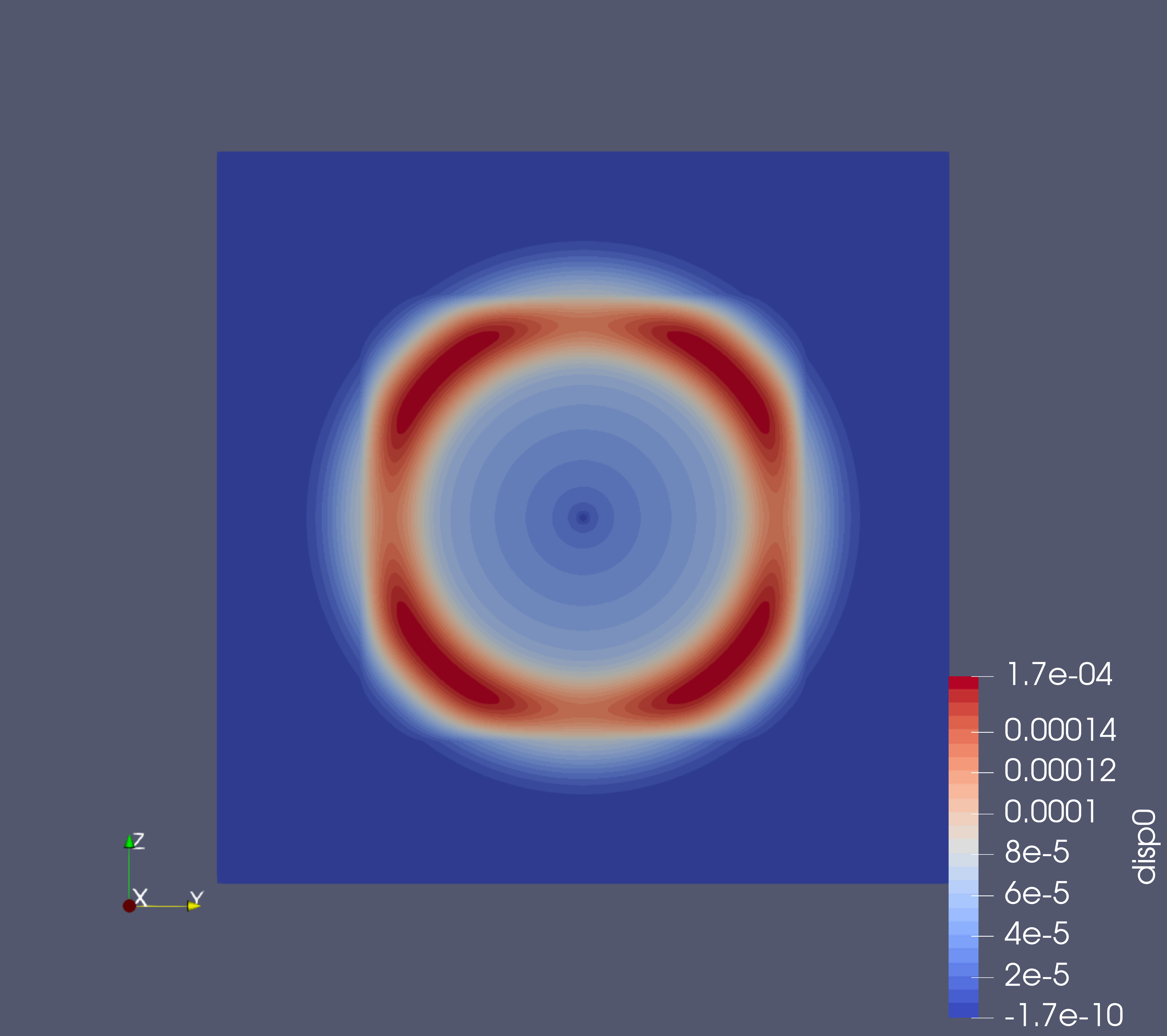}
		\includegraphics[width=0.45\textwidth]{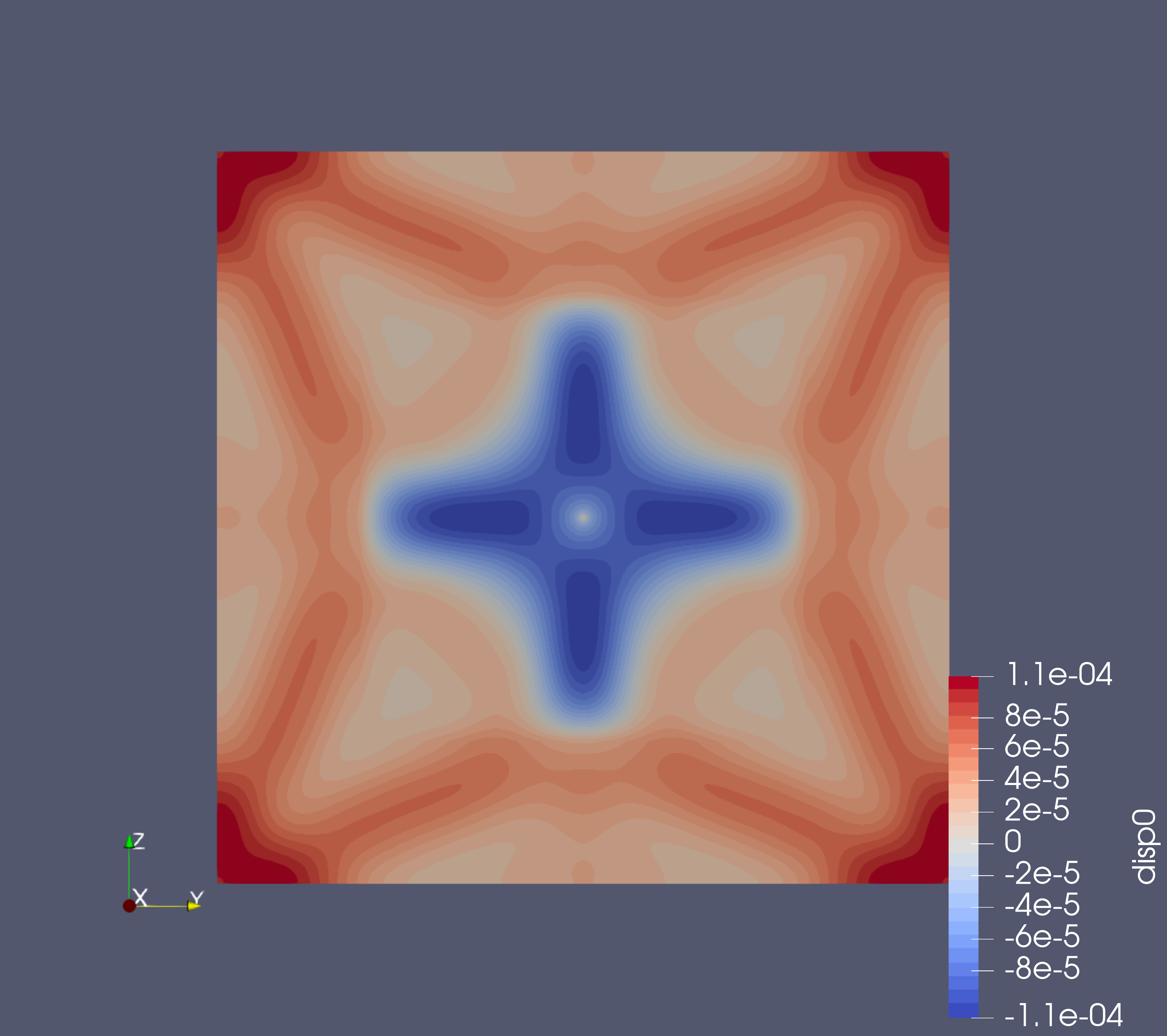} 
		\includegraphics[width=0.45\textwidth]{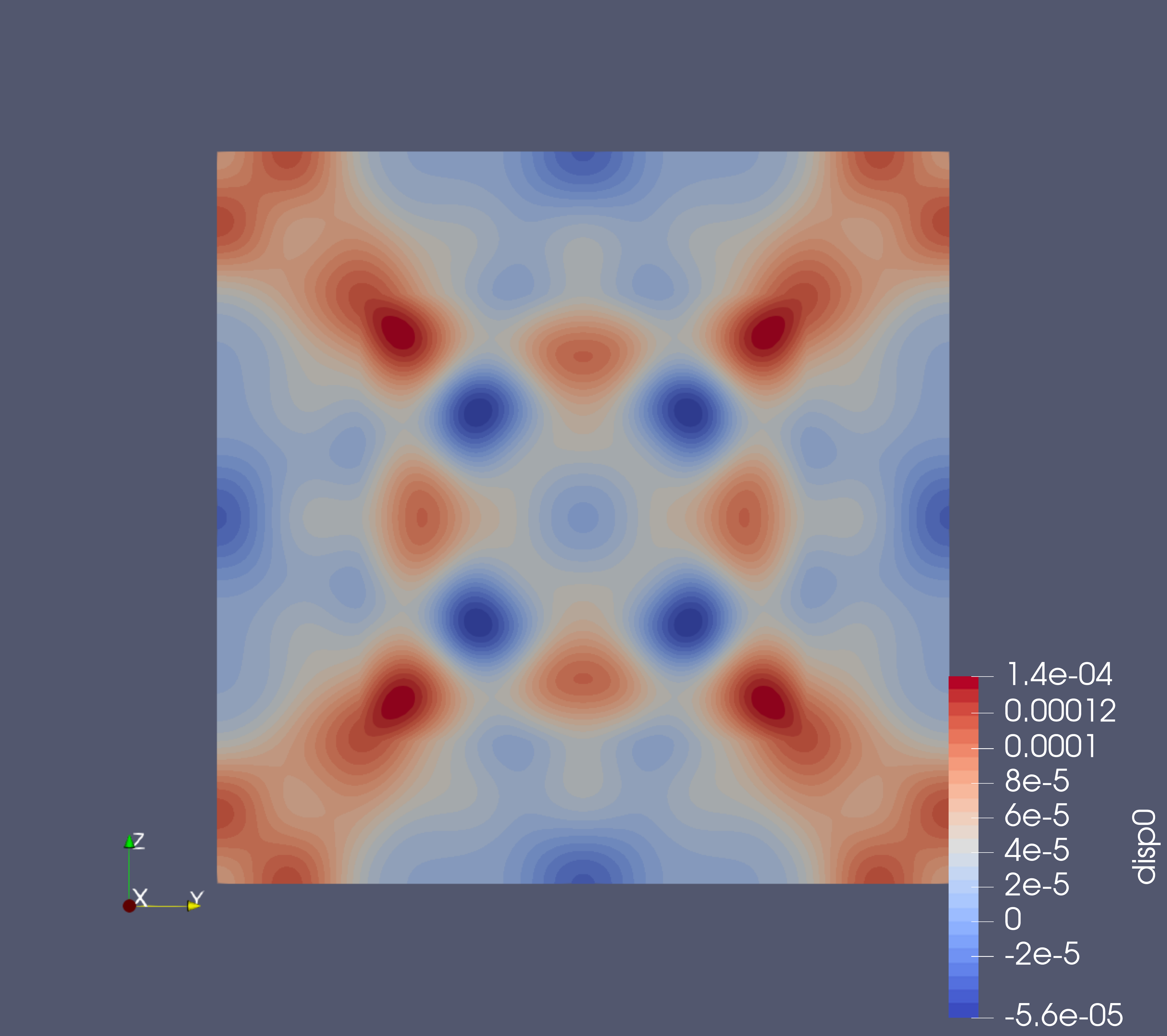}
	\caption{Transverse displacement for circular waves in plane strain traveling in Al framed in Fe, $t= 1.91 \cdot 10^{-5}$s, $ 1.115 \cdot 10^{-4}$s, $ 2.3 \cdot 10^{-4}$s, $ 2.86 \cdot 10^{-4}$s}
	\label{fig:circular_hetero}    
\end{figure}

This case is a modification of the previous one using a heterogeneous medium. The size and shape of the domain are the same, a square plate of dimensions $0.001\times1\times1$m, but the material occupying an internal square of dimensions $0.39\times0.39$m is made of Al, while the external frame is made of Fe. The properties of both materials are given in  table \ref{tab: 1D-mat-parameters}. The simulation conditions are the same as in the homogeneous case, including time discretization, grid of $1\times1025\times1025$ voxels, shape of the pulse applied, and area where applied. The results of this simulation are represented in Fig. \ref{fig:circular_hetero}, which shows the contour plots of the displacement in the direction perpendicular to the plate at three different times.
 
The solutions of both cases are identical until the wave reaches the frame. In terms of computational cost, the homogeneous problem requires only evaluating the preconditioner at each time step, and the total simulation time was 10 minutes on the desktop computer. The heterogeneous problem required approximately 10 CG iterations per time step for a tolerance $10^{-8}$, and the total time was 1 hour and 50 minutes using the same computer. These numbers are remarkable for a three-dimensional problem with a million points and a standard desktop computer. As an illustration, we tried to solve the same problem in FE with one element per voxel and periodic boundary conditions, but the required memory exceeded the computed RAM of 8GB. In order to use the same computer for clear time comparisons, periodicity was eliminated in the FE model to reduce the memory needs. In this case, the total simulation time was 3 days. The comparison shows that the FFT-based model was $\approx \times$ 36 faster than the implicit FE simulation, even after eliminating the limiting periodicity condition for FE.
 
\subsubsection{Spherical waves in a homogeneous solid}
This simulation corresponds to the propagation of waves in a homogeneous solid medium caused by the application of an excitation in its center. The material is Al. The periodic domain is a cube of dimensions $1\times1\times1$m, discretized in $129\times129\times129$ voxels, with around $2\cdot 10^6$ voxels. The time increment used is d$t=1.52 \cdot 10^{-6}$ s, above $CFL$, but sufficiently small to have a complete deformation map in 6 time steps before reaching the periodic boundary. The final time was $T=2.2 \cdot 10^{-4}$s. The displacement is imposed in the center of the volume, in a region with radius $R=1$cm, using the pulse of amplitude 1mm defined in Eq. (\ref{eq-prescribed_displacement}).  The direction of the pulse is forced in $x_1$, and the resulting wave has a mixed character, being longitudinal or transverse depending on the position of the point with respect to the origin, which defines the propagation direction for that point.

\begin{figure}[H]
    \centering
		\includegraphics[width=0.45\textwidth]{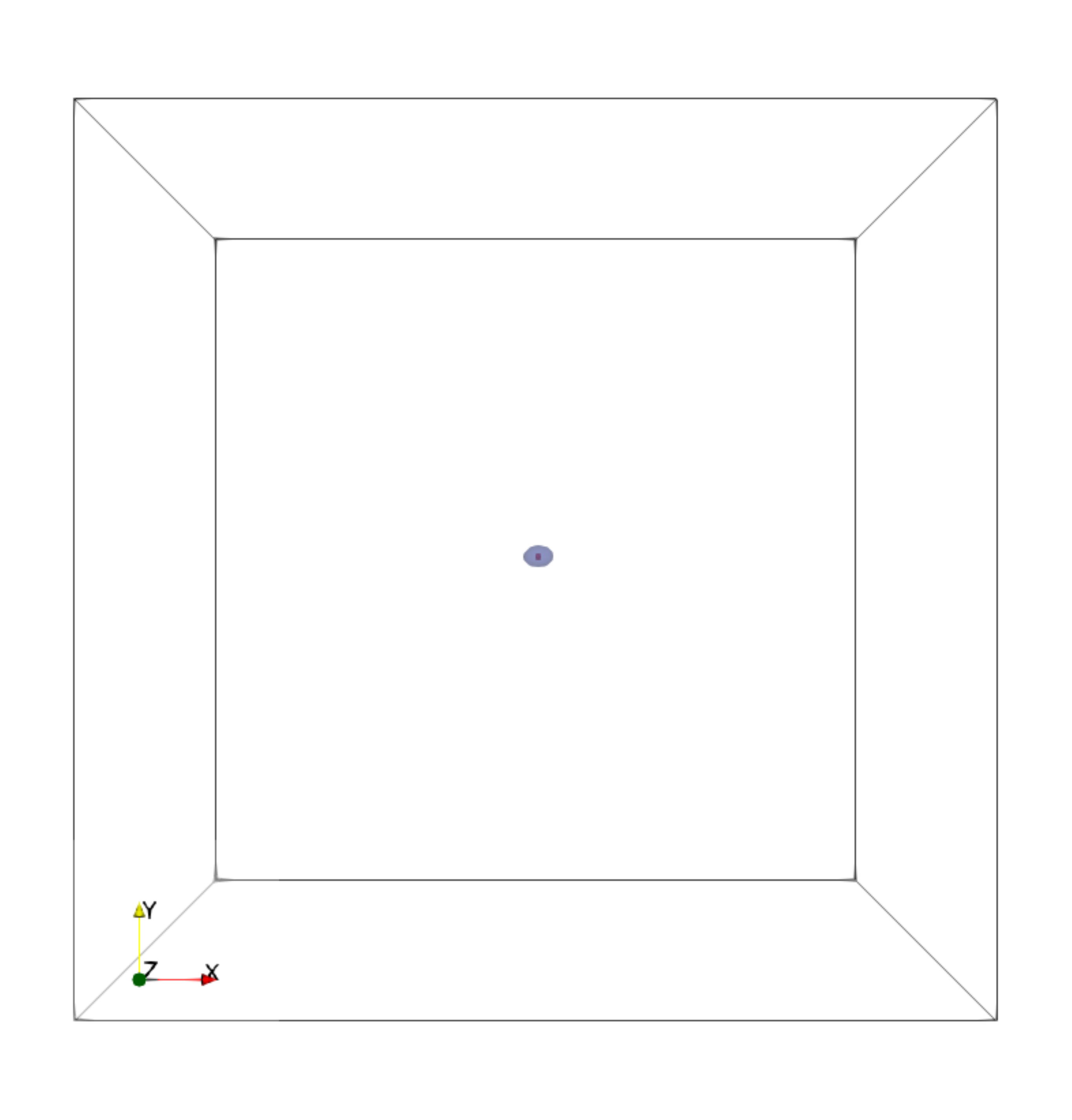} 
		\includegraphics[width=0.45\textwidth]{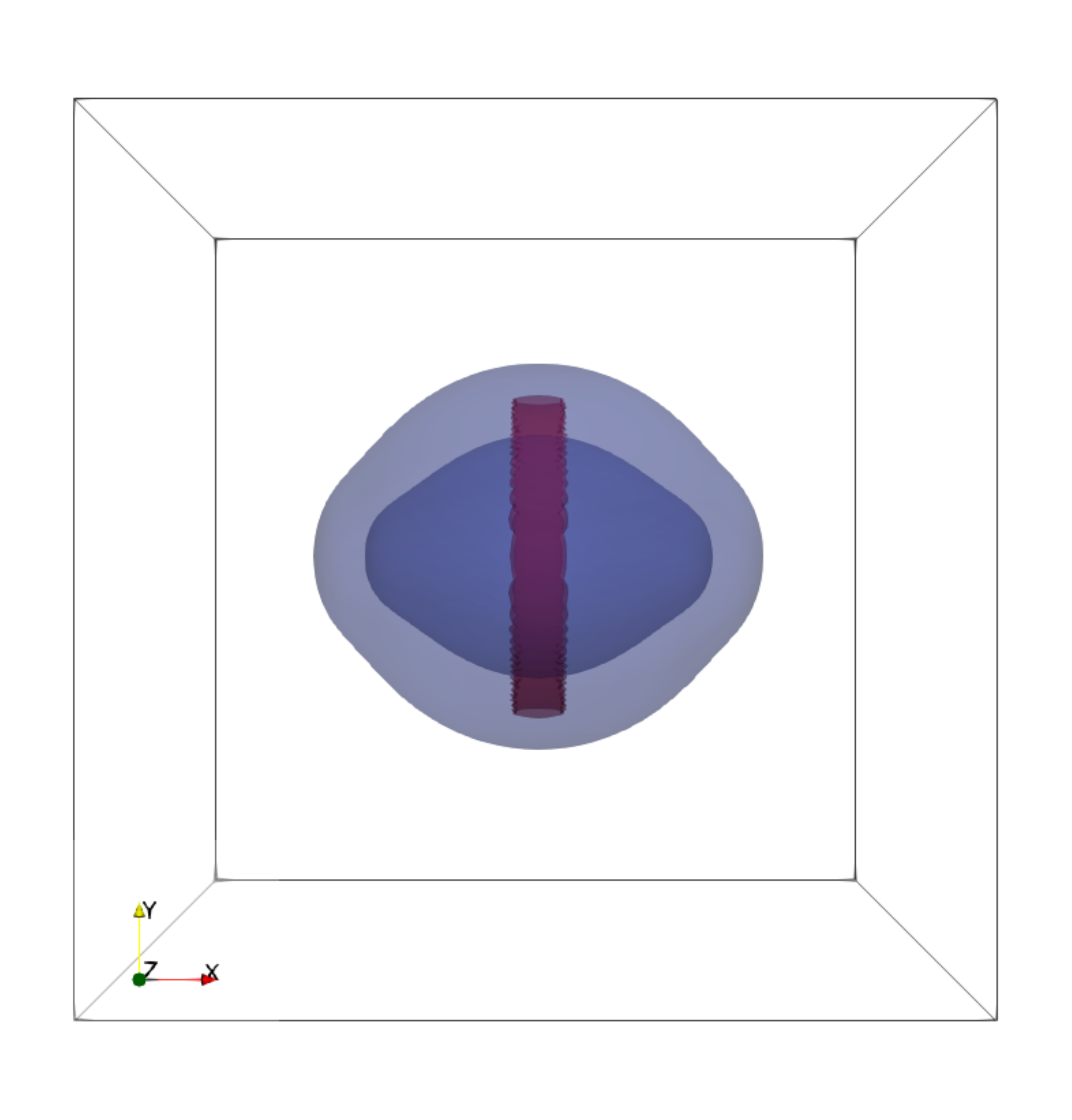}
		\includegraphics[width=0.45\textwidth]{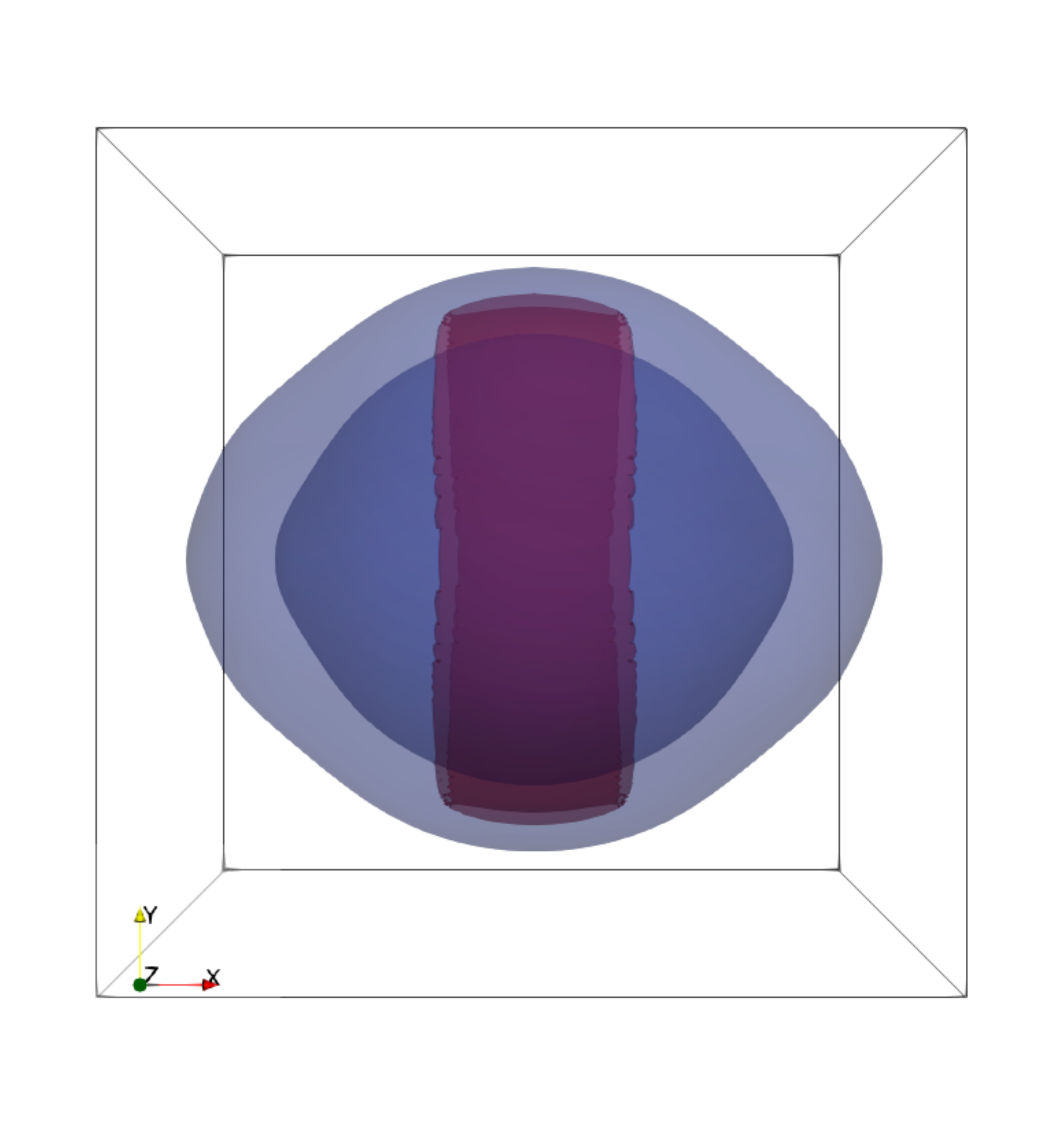} 
		\includegraphics[width=0.45\textwidth]{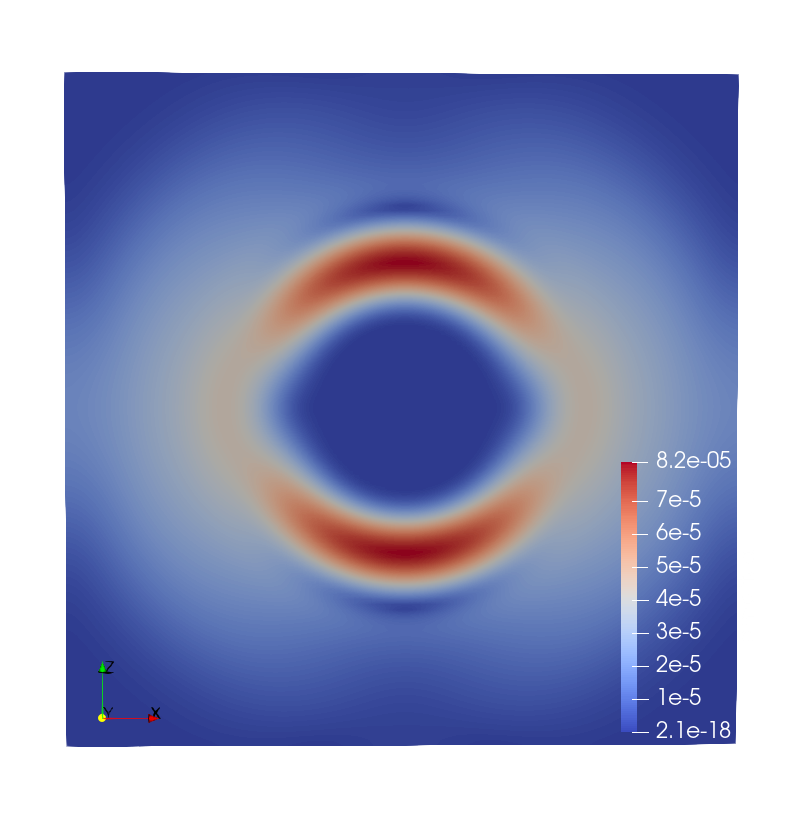}
	\caption{(a-c) Surfaces having displacements amplitudes of $0.5A_{max}$ (blue) and $0.9A_{max}$ (red) for $t=3.04 \cdot 10^{-6}$s ,$9.13 \cdot 10^{-6}$s ,$1.37 \cdot 10^{-5}$. (d) Contour plot of the amplitude value in the plane $x_1-x_3$ passing through the cell center.}
	\label{spherical_homo}  
\end{figure}

The results of the simulation are represented in Fig. \ref{spherical_homo}. Figs. a-c show the locus of the points in which the amplitude is equal to 0.5 and 0.9 the maximum amplitude ($A_{max}$) for three different times. It can be observed that the shape of the constant amplitude locus has cylindrical symmetry with respect to $x_1$, but, contrary to the 2D case, the point symmetry is lost. The locus of the points with 0.5 normalized amplitude is an elongated spheroid in the direction $x_1$, as the waves are longitudinal in this direction. The locus of the points with maximum amplitude (in Fig. \ref{spherical_homo}) is a ring with axis $x_1$ that corresponds to the points of the front more near the origin of the pulse, in which excitation is a transverse wave. A contour plot showing the amplitude values in $x_1-x_2$ is represented in Fig. \ref{spherical_homo}(d). 
The amplitude decay with the distance to the center, $r$, accurately follows the theoretical relation of $\approx 1/r$, and no numerical dissipation has been found. Moreover, the results do not show any spurious oscillations. It should be noted that, in this homogeneous simulation, the implicit integration has a closed expression in Fourier space (the preconditioner in Eq. \ref{eq:precond}) and therefore the solution of the system does not require the use of the conjugate gradient. The time needed for the entire simulation was only 20 minutes on the desktop computer. In FE, the memory requirements exceeded the computer's capabilities even after the periodicity was removed.

\subsection{Spherical waves in a framed solid}
\begin{figure}[H]
    \centering
		\includegraphics[width=0.45\textwidth]{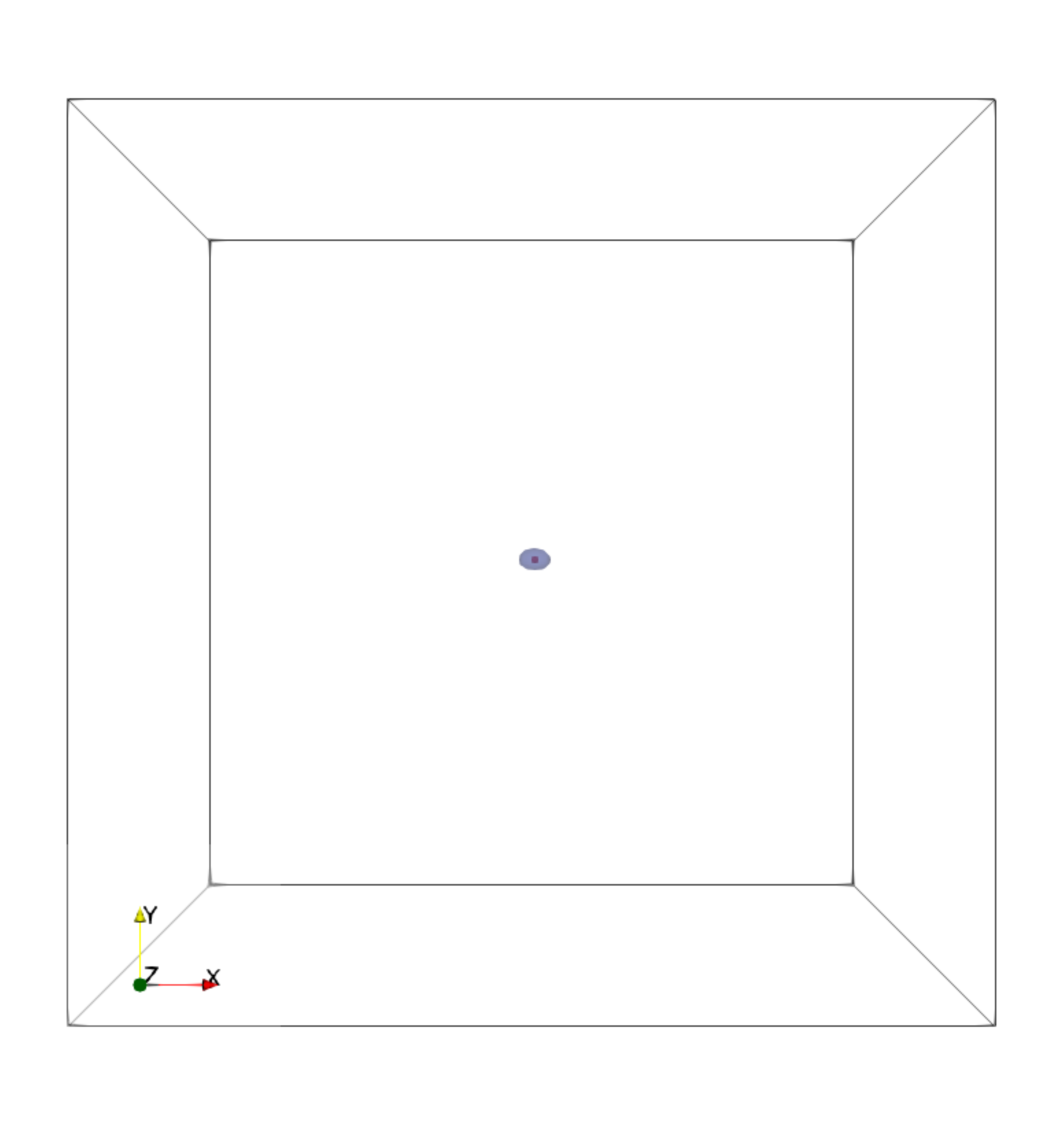}
		\includegraphics[width=0.45\textwidth]{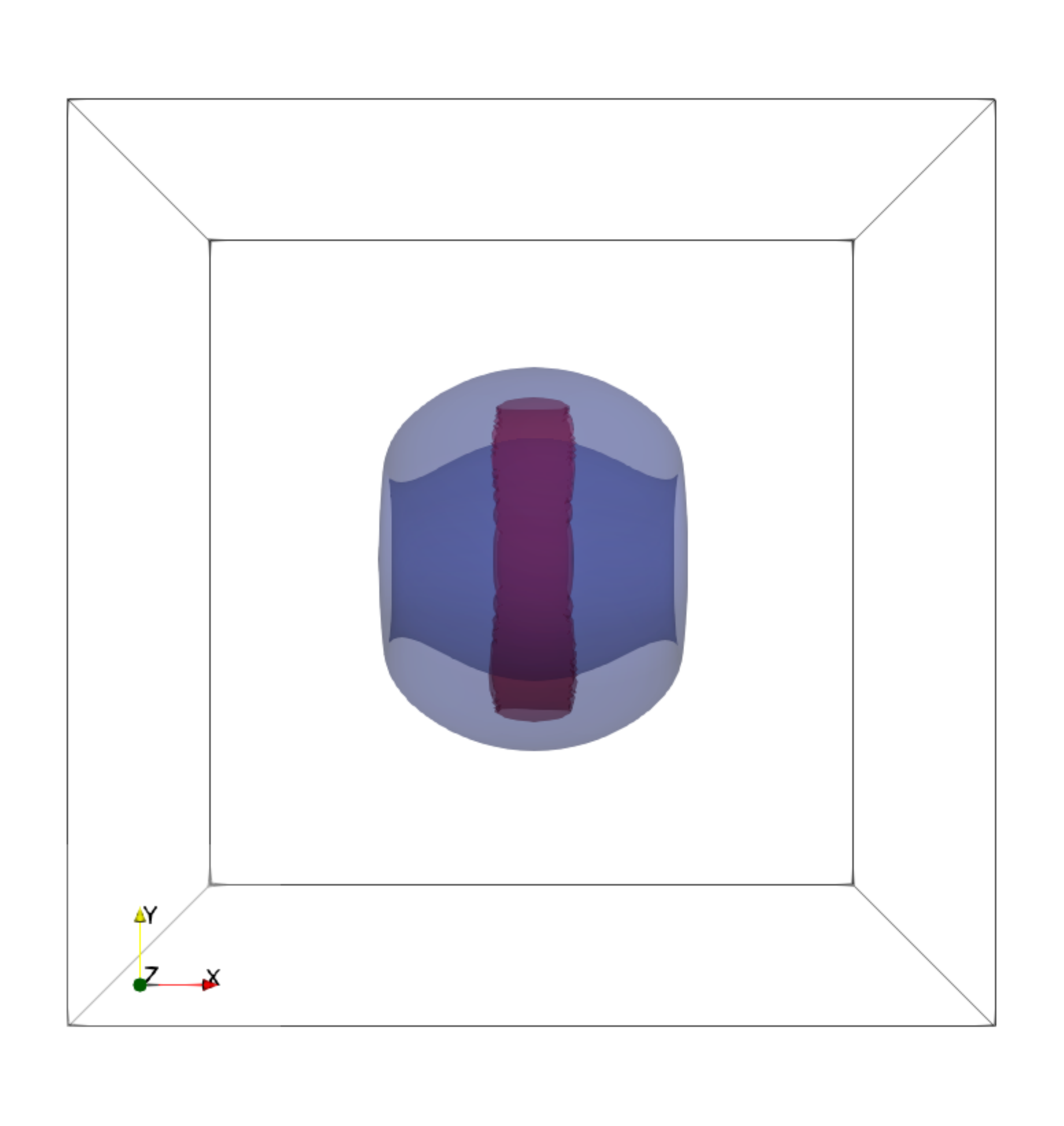}
		\includegraphics[width=0.45\textwidth]{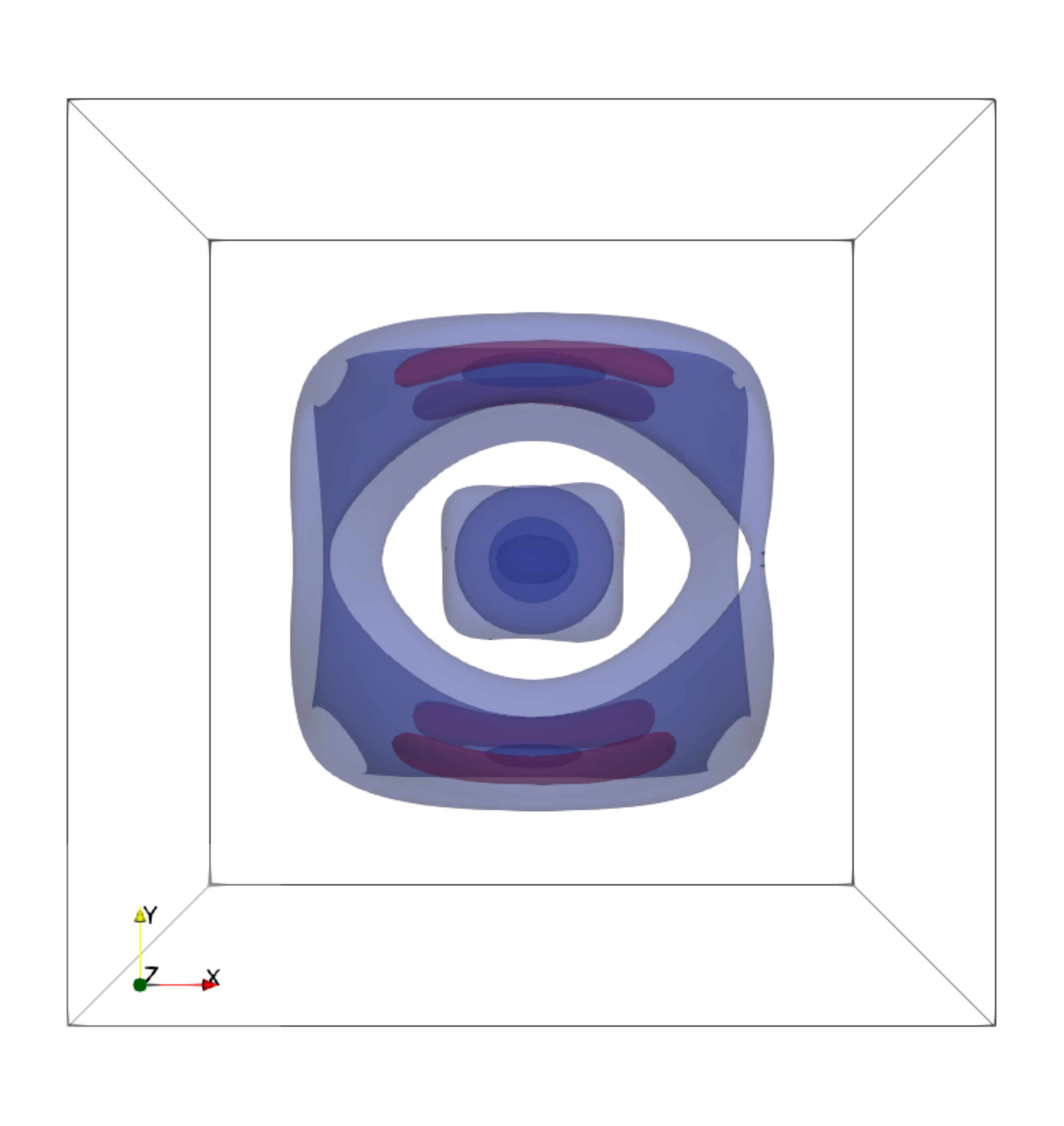}
		\includegraphics[width=0.45\textwidth]{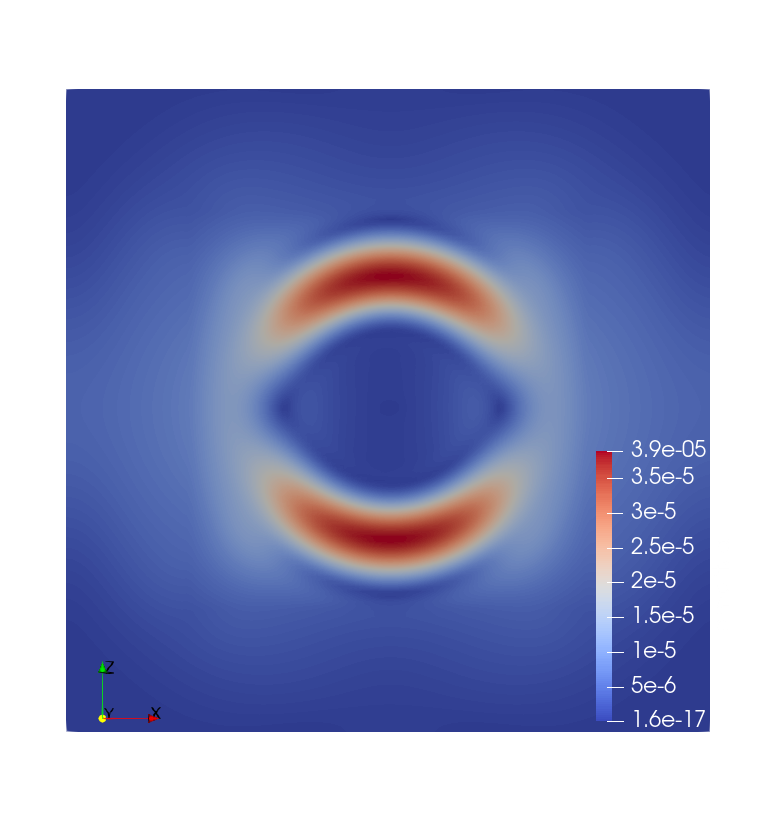}
	\caption{(a-c) Surfaces having displacements amplitudes of $0.5A_{max}$ (blue) and $0.9A_{max}$ (red) for $t=3.04 \cdot 10^{-6}$s ,$9.13 \cdot 10^{-6}$s ,$1.37 \cdot 10^{-5}$. (d) Contour plot of the amplitude value in the plane $x_1-x_3$ passing through the cell center.}\label{spherical_hetero}
\end{figure}

This next example corresponds to the propagation of spherical waves in a solid medium composed of two different materials. The inner cube, with length 0.62 m, is made of aluminum, while the outer frame with thickness of 0.194 m on each face is made of iron. The properties of both materials are the same as in the rest of the article. The applied perturbation,  spatial discretization, time increment size, and final time are the same as in the homogeneous case. The results of the simulation are represented in Fig. \ref{spherical_hetero}. 

Figs. \ref{spherical_hetero} (a)---(c) show the locus of the points in which the amplitude is equal to 0.5 and 0.9 the maximum amplitude for three different times. The shape and symmetries of the constant amplitude locus are the same as in the homogeneous case before the wave reaches the frame (Fig. \ref{spherical_hetero}(a)), and changes after reaching those points (Figs. \ref{spherical_hetero}(b-c)). Axial symmetry is also lost when the wave touches the frame because of the lack of spherical symmetry of the cubic unit cell. The contour plot of the displacement amplitude represented in Fig. \ref{spherical_hetero}(d) corresponds to the incipient contact of the wave with the frame, showing a slightly different shape than the one represented in Fig. \ref{spherical_homo}(d) for the homogeneous solid at the same time. 

The solution in this heterogeneous solid preserves the smoothness and accuracy of the solution in the homogeneous medium. Regarding the computational cost, the number of iterations of the conjugate gradient was less than 10 at each time step. The total time needed for the entire simulation was less than 3 hours on the same computer. Again, simulation with FE for time comparison was not possible on the same computer due to memory restrictions, even using standard Dirichlet boundary conditions.

\subsection{Wave propagation in polycrystal}

\begin{figure}[H]
    \centering
    \includegraphics[width = 0.9\textwidth]{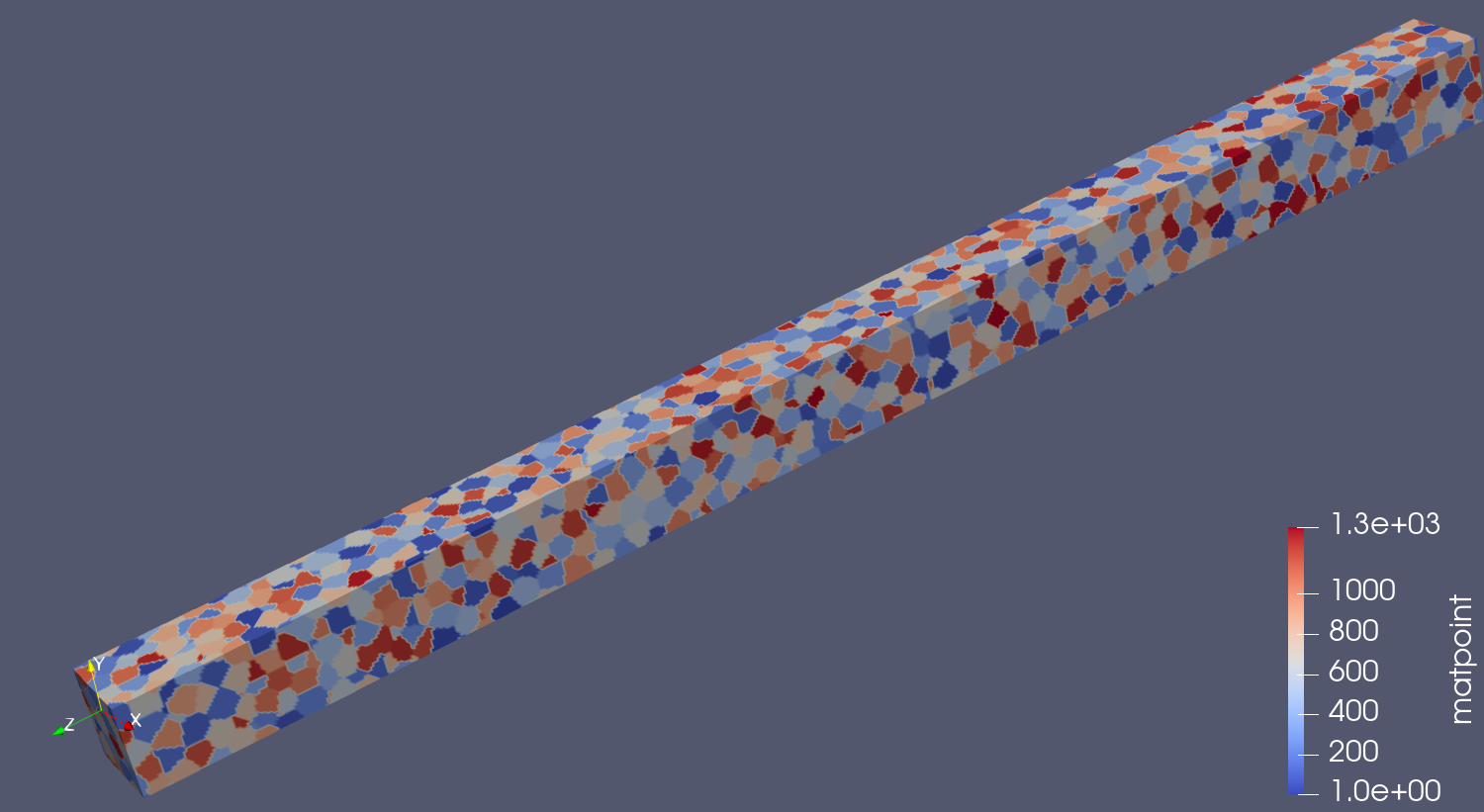}
    \caption{Microstructure of the polycrystal. $35 \times 35 \times 567$ voxels and $0.35 \times 0.35 \times 5.67$ mm. Mean grain diameter = 100 $\mu$m}
    \label{fig:polycrystal_microstructure}
\end{figure}

Finally, we study the propagation of the elastic wave in a polycrystalline RVE in which the grains are explicitly represented. The domain is a prismatic bar of dimensions 0.35 $\times$ 0.35 $\times$ 5.67mm, discretized in 35 $\times$ 35 $\times$ 567 voxels. The 3D periodic microstructure is represented in Figure \ref{fig:polycrystal_microstructure} and was generated using a weighted Voronoi tessellation to statistically represent a log-normal grain size distribution with a mean grain diameter of $\overline{d}_g = $ 100 $\mu$m and a standard deviation of 5 $\mu$m. The properties of the single crystal correspond to Ni, a very anisotropic crystal with a Zener ratio greater than 2 with elastic constants $C_{11} = 249 $GPa, $C_{12} = 155 $GPa, $C_{44} = 114 $GPa and density $\rho = 8908 \ kg/m^3$ . The crystallographic texture was adopted to be random and therefore the elastic stiffness tensor differ from grain to grain in a random fashion, resulting in a macroscopic isotropic elastic response.

As in Sections 3.1 and 3.2, the bell-shaped pulse $U(t)$ defined in Eq. \eqref{eq-prescribed_displacement} is used to prescribe the displacement boundary conditions $u_3 = U(t)$; $u_1 =u_2= 0$ in the plane $x_3 = 0$. The amplitude is set to $A = 5.67 \cdot 10^{-3}$mm while two different pulse lengths (by changing the parameter $\omega$) are used to analyze two situations: 1) a short pulse with respect to the grain size, approximately, 3 $\overline{d}_g$ ($\omega = \pi \sqrt{c_{11}/\rho} / (3\overline{d}_g )$) and 2) a long pulse with respect to the grain size, 20 $\overline{d}_g$ ($\omega = \pi \sqrt{c_{11}/\rho} / (20\overline{d}_g )$). The total time of the first simulation (3$\overline{d}_g$) was around 11 hours, with approximately 80\% of the time consumed in calculating the Green's functions. Reusing the $\mathbf{G}$ matrix built with Green's functions, the second simulation (20$\overline{d}_g$) only took 2 hours. 

\begin{figure}[H]
    \centering
	\subfloat[][Pulse length = 2-3 grain sizes]{
		\includegraphics[width=0.9\textwidth]{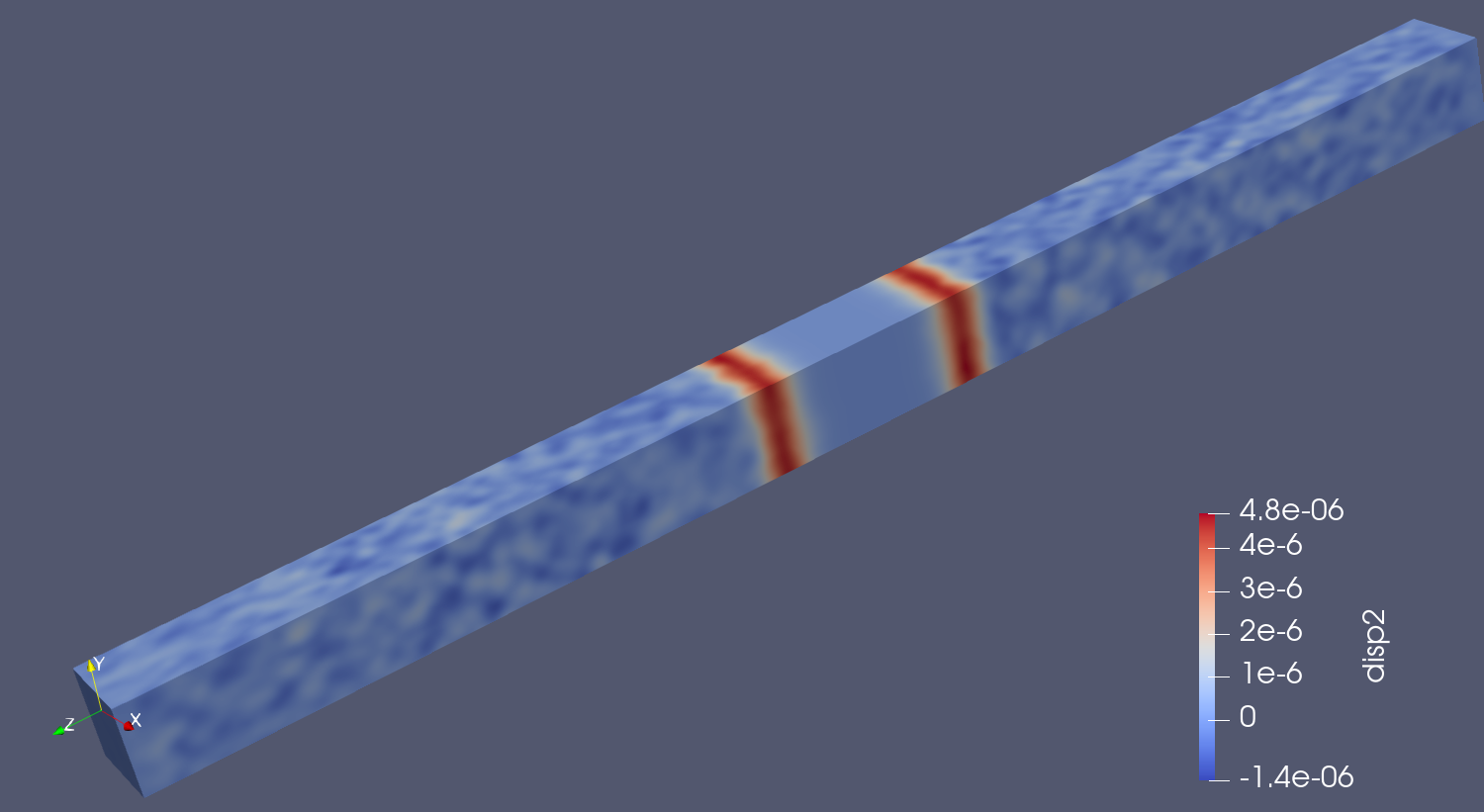}} \\
	\subfloat[][Pulse length = 20 grain sizes]{
		\includegraphics[width=0.9\textwidth]{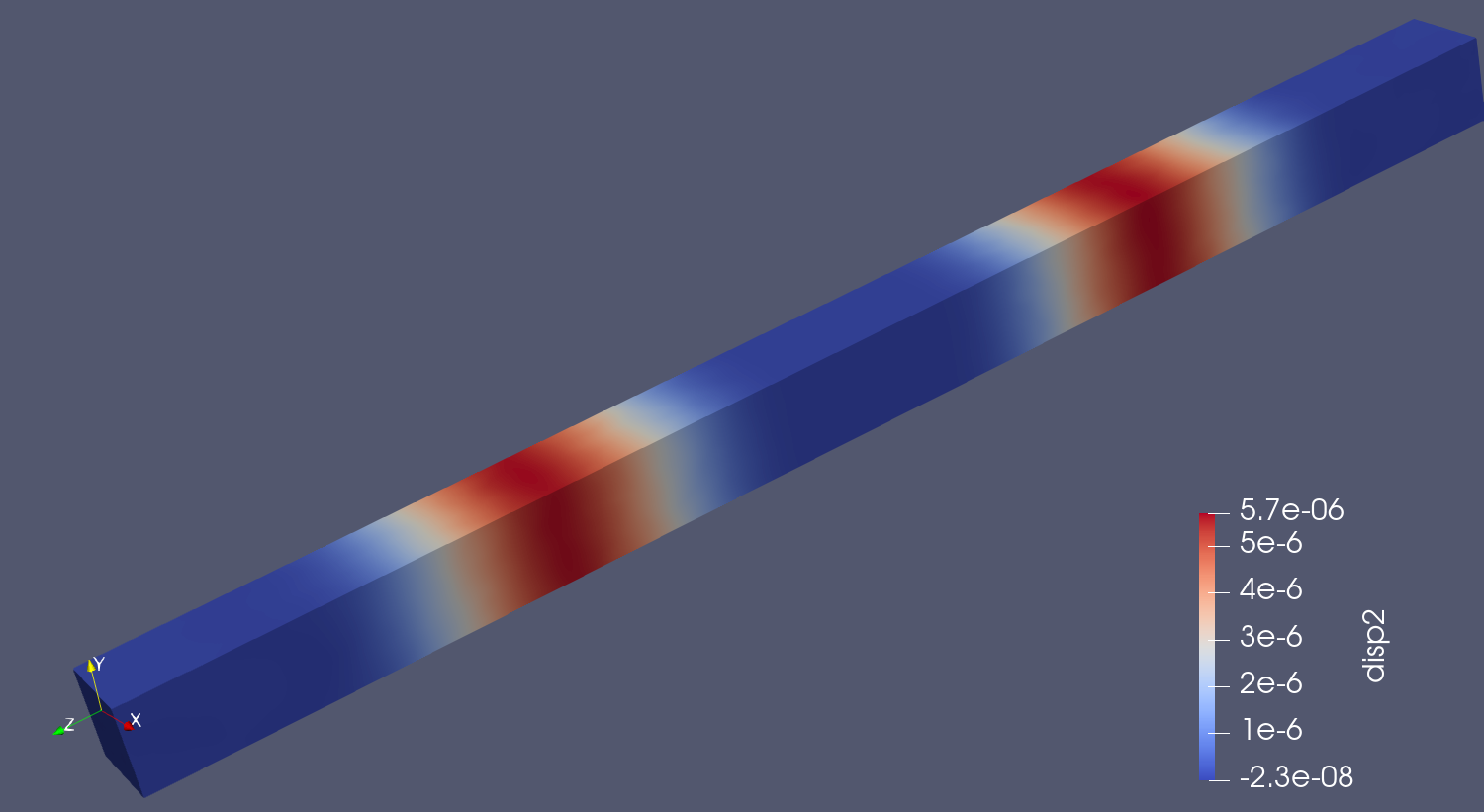}}
	\caption{$u_3$ displacement field at the same time for the same micro-structure but for two different pulse lengths.}
	\label{fig:polycrystal}    
\end{figure}

Figure \ref{fig:polycrystal} shows the contour plot of the field $u_3$ for the two cases for a time in which waves have propagated to approximately 1/4 of the domain length. In the first case, it can be seen that the grains promote the scattering of the wave, whereas in the second case the wave interaction with the micostructure is minimal and behaves similarly to a wave traveling in a homogeneous solid.

\begin{figure}[H]
    \centering
    \includegraphics[width = 0.45\textwidth]{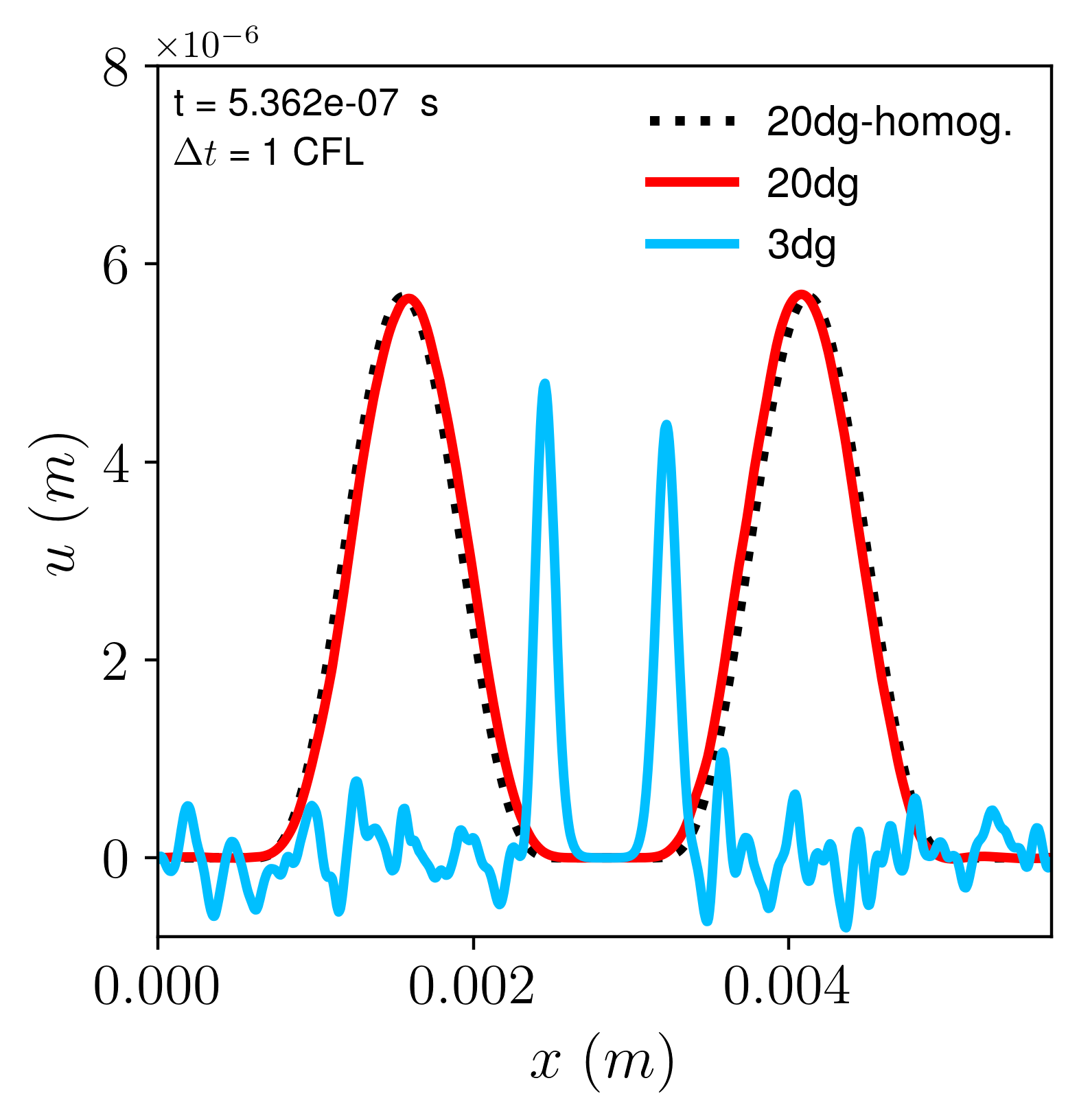}
    \includegraphics[width = 0.45\textwidth]{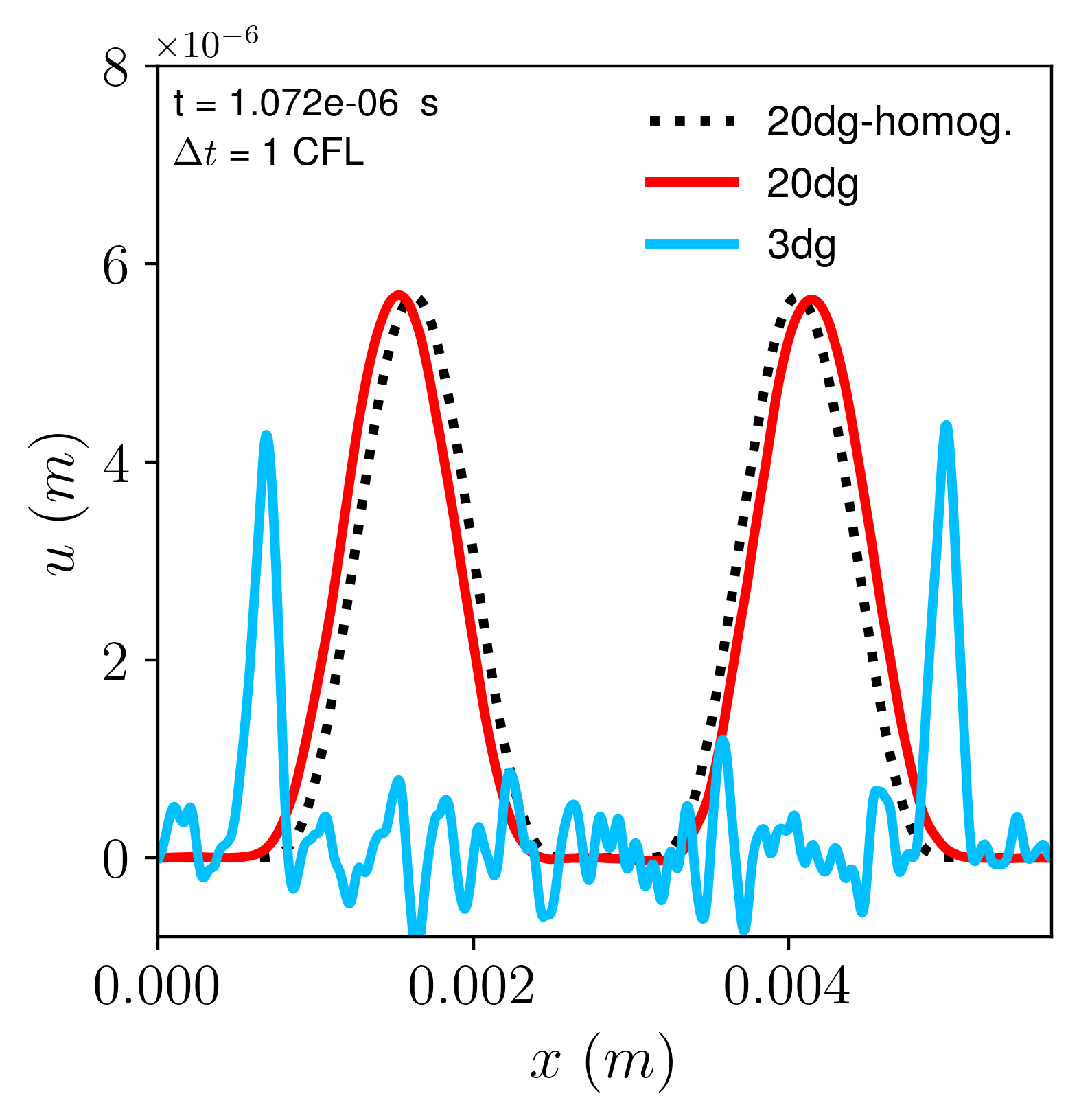}
    \caption{$u_3$ field along the polycrystalline bar length}
    \label{fig:polycrystal_plot}
\end{figure}

To further analyze the interaction of the waves with the single crystal grains, the longitudinal displacement along a line in the $x_3$ direction is represented in Fig. \ref{fig:polycrystal_plot} for the two wave lengths at different times. The results in Fig. \ref{fig:polycrystal_plot} also include the results for a homogeneous equivalent isotropic solid. The properties of the homogeneous material are obtained by elastic computational homogenization of the same domain, using the DBFFT approach \cite{Lucarini2019c}. Owing to the random texture of the polycrystalline unit cell, the homogenized material was isotropic with elastic constants $E=198$ GPa, $\nu=0.306$.
It can be observed that when the pulse is sufficiently larger than the characteristic length of the microstructure, the solid behaves similarly to a homogeneous material, as observed in similar simulations using FE \cite{vanpamel2017}. Nevertheless, even for the long wave, the maximum of the wave traveling in the polycrystal was slightly delayed with respect to the homogeneous material. This reduction in group velocity is in agreement with the results obtained by Segurado and Lebensohn \cite{Segurado2021} using a completely different approach based on computing the dispersion relations. The results in the case of a short pulse are very different. In this case, it can be observed that the displacement behind the wave front is extremely wavy (Fig. \ref{fig:polycrystal_plot}). This oscillation corresponds to the reflections and refractions of the traveling wave on the grain boundaries. As a result, although the total elastic energy is conserved in the full bar, the energy around the wave front is progressively reduced. This loss of energy will eventually lead to the disappearance of the traveling wave, as reproduced using massive parallel FE simulations \cite{vanpamel2017,Ryzy2018}. The technique proposed here will allow us to study the maximum propagation distance for a general microstructure in a very efficient manner compared to full-field simulations using explicit FE.

\pagebreak
\section{Conclusions}

A novel FFT-based algorithm has been described to simulate propagation of elastic waves in heterogeneous $d$-dimensional rectangular shape domains. The method allows prescribing the displacement as a function of time in a subregion of the domain, emulating the application of Dirichlet boundary conditions on an outer face of the domain. Time discretization is performed using an unconditionally stable implicit beta-Newmark approach. The implicit problem for obtaining the displacement at each time step is solved by transforming the equilibrium equations into Fourier space and solving the corresponding linear system with a preconditioned Krylov solver. For comparison purposes, an explicit version of the FFT model is also proposed, following a central difference integration scheme.

The main conclusions of this work are:
\begin{itemize}
    \item The framework is capable of reproducing the propagation of plane, circular, and spherical waves very accurately and without numerical noise.
    \item The Green's function method used to impose the prescribed displacement was totally equivalent to applying time-dependent Dirichlet boundary conditions in a boundary value problem using FE  .
    \item The use of Newmark implicit integration allows us to resolve large models, preserving the same accuracy and absence of noise,  using time steps orders of magnitude larger than Courant's $CFL$ condition. 
        
    \item The study of the order of accuracy show that the error reduction with time is quadratic up to a limiting time, near the Courant condition. Below that time increment, accuracy is not improved by time increment reduction. 
    \item The order of growth of the computational cost of the FFT implicit model is much lower than for implicit FE. The difference in time using the exact conditions in both solvers becomes order of magnitude smaller for sufficiently large models.
    \item The implicit FFT method using a fixed time step became also faster than the explicit FE, due to the progressive reduction of the time step in explicit FE. This improvement in efficiency is important considering that the accuracy in FFT was much better than in explicit FE.
    \item The use of very large phase contrasts reduces computational efficiency and introduces noise in the results
    \item In summary, the method is an ideal framework for studying the propagation of waves in heterogeneous media where very fine discretizations are needed. The present approach allows one to solve large problems with the same or better accuracy at a fraction of the cost of using FE.
\end{itemize}

We have identified several features and possible improvements, that will be covered in upcoming research:
\begin{itemize}

    \item The periodicity enforces the presence of two traveling waves. It will be interesting to suppress one of the waves for a more clear analysis of the results. The use of a buffer of a very compliant medium or a viscous material is a potential solution.
    \item The use of staggered grids or discrete finite difference differentiation rules could be explored to reduce noise in the presence of large mechanical contrast.
\end{itemize}
 These potential extensions and improvements will allow us to use this approach to study complex problems explicitly considering microstructure, including simulations of impact, resonant ultrasound spectroscopy (RUS), spallation, Hopkinson-bar tests, etc.
\clearpage

\section{Acknowledgments}
 Javier Segurado acknowledges the European Union's Horizon 2020 research and innovation program for the MOAMMM project, grant agreement No. 862015, of the H2020-EU.1.2.1. - FET Open Program and the Spanish Ministry of Science for the project ADSORBENT, Plan estatal de I+D+i-20019: PID2019-106759GB-I00.
 
Ricardo Lebensohn acknowledges support from Los Alamos National Laboratory's Laboratory-Directed Research \& Development (LDRD) Program.

Rafael Sancho-Cadenas gratefully acknowledges the support received for this work under grant PID2020-116440RA-I00 funded by MCIN/AEI/10.13039/501100011033.

\appendix

\section{One dimensional Green's function}\label{anex:green}
In this Appendix, the analytical expression of the Green's function of the spatial operator resulting from Newmark integration in the case of homogeneous materials is analyzed in order to show that it is nonsingular everywhere and to analyze its numerical representation in FFT.

In the one dimensional case, the linear operator of the equation of linear momentum conservation integrated using $\beta-$Newmark for a homogeneous medium is obtained by particularizing Eq. \eqref{eq:1D_consv_lm_Newmark_ord} to $E(x)=E$ and $\rho(x)=rho$
\begin{equation}
\mathcal{A}(u_n) =  \frac{\mathrm{d}^2 u_n}{\mathrm{d}x^2} - \frac{\rho}{\beta \Delta t^2 E }u_{n}.
\end{equation}
This equation corresponds to a 1D Helmholtz equation as
$$
\frac{\mathrm{d}^2 u}{\mathrm{d}x^2} - K^2 u = 0 
$$
with $K= \left(\frac{\rho}{\beta \Delta t^2 E }\right)^{1/2}$. The Green´s function of this linear operator in an infinite medium where $\lim_{x\rightarrow \pm \infty} u = 0$ can be easily found integrating the equation with a delta function as right hand side. The resulting  expression is
\begin{equation}
g(x,x') = -\frac{\exp(-K|x-x'|)}{2K}.
\label{eq:green1Db}
\end{equation}
It can be observed that the function is not singular when evaluated at $x=x'$, being its value $g(x,x)=-1/2K$. 

\begin{figure}
    \centering
    \includegraphics[width=0.48\textwidth]{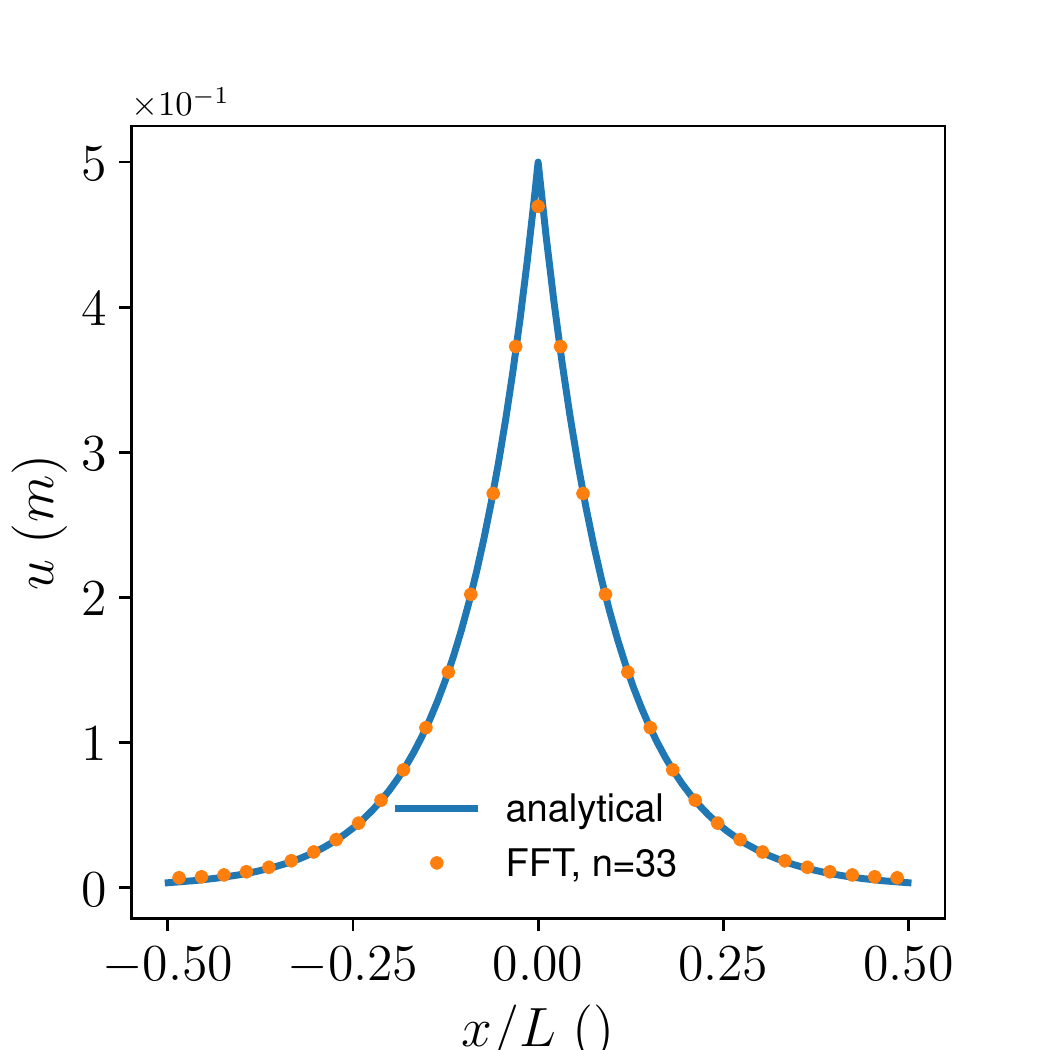}
    \caption{1D Green operator, numerical solution with N=33 and analytical solution}
    \label{fig:green1D}
\end{figure}

The FFT resolution of equation $\mathcal{A}(u_n)=\delta$, provides a numerical result of the Green's function in the case of a periodic domain, that is the one that is used in the simulations. To illustrate the accuracy of the numerical Green's function, in the figure \ref{fig:green1D} the analytical expression Eq. \eqref{eq:green1Db} is represented together with the FFT solution for a coarse grid size (33 voxels) in a sufficiently long domain to avoid the effect of periodicity. It can be observed that points of the numerical solution lie on the analytical one even for this coarse discretization.

In the case of a heterogeneous medium, the Green's function depends on the particular microstructure. However, if a homogeneous reference medium $E$, $\rho$ is defined, the heterogeneity can be accounted as a source field acting on that reference medium, and the Green's function is preserved. This is the usual approach in FFT-based homogenization approaches based on polarization.

\section{Three dimensional Green's function}
In this appendix, the analytical expression of the Green's function of the 3D operator for a homogeneous materials and its FFT representation are analyzed.

\begin{figure}
    \centering
    \includegraphics[width=0.48\textwidth]{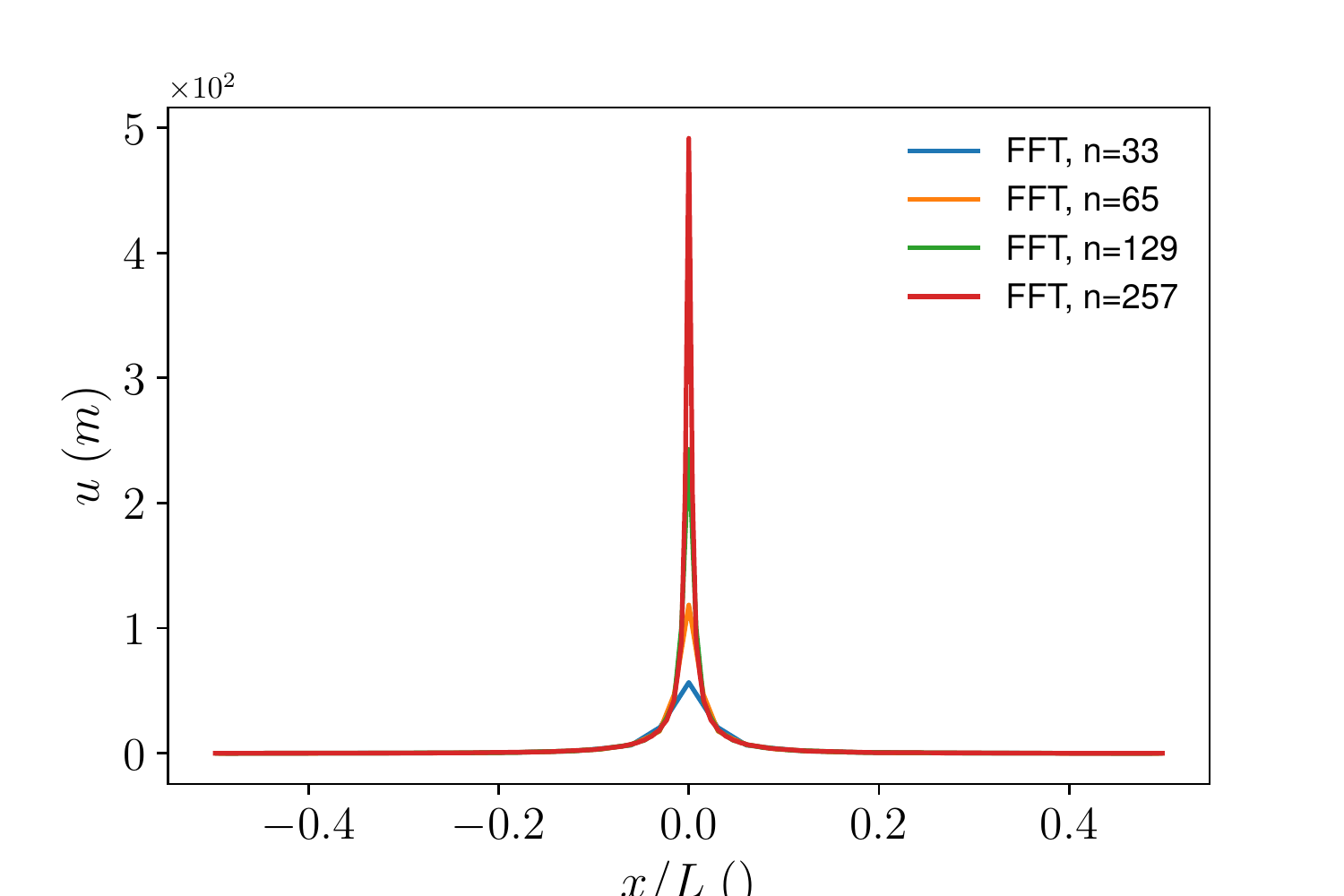} 
    \includegraphics[width=0.48\textwidth]{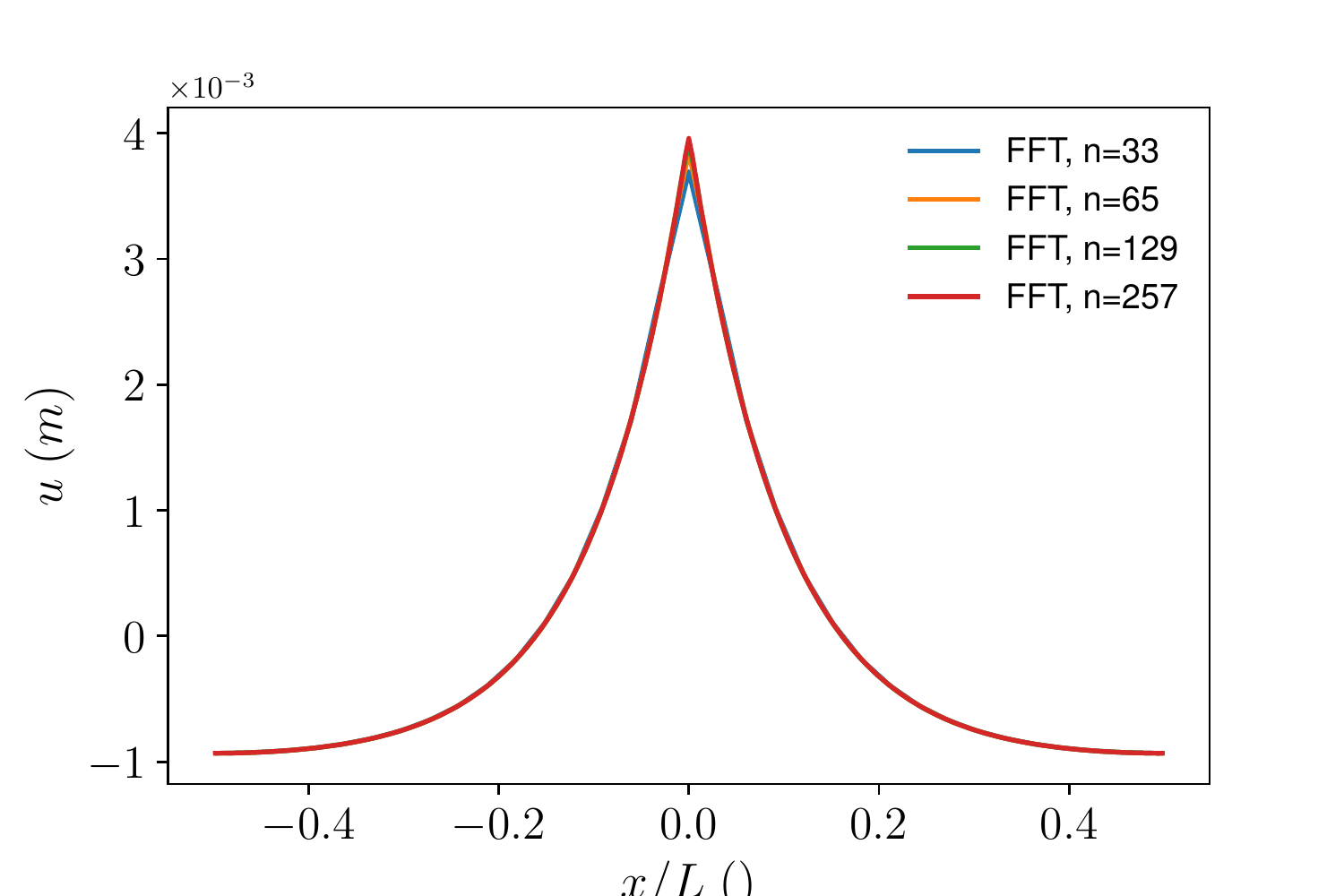}
    \caption{1D Green operator, numerical solution with N=33 and analytical solution}
    \label{fig:green3D}
\end{figure}

The linear operator of Eq. \eqref{eq:ND_consv_lm_Newmark_ord} for a homogeneous medium corresponds, as in the 1D case, to a Helmholtz equation with a negative identity term. To simplify the analysis, an idealized homogeneous medium with stiffness equals to $\mathbb{C}=E\mathbb{I}^s$ is considered. In this case, the problem becomes scalar and corresponds to
\begin{equation}
\nabla^2 u - K^2 u = 0
\label{eq:helmholtz3d}
\end{equation}
with $K=\left(\frac{\rho}{\beta \Delta t^2 E }\right)^{1/2}$. The Green' s function of this problem is
\begin{equation}
g(\mathbf{x},\mathbf{x}') = g(\mathbf{x}-\mathbf{x}')=g(r)=-\left( \frac{K}{8\pi^3r}\right)^{1/2} \text{K}_{1/2}(Kr) \label{eq:green1D}
\end{equation}
with $\text{K}(\cdot)$ the modified Bessel function of the second kind. Contrary to the 1D case, this function is singular at $\mathbf{x}=\mathbf{x}'$, and therefore the displacement caused by a unit force on its application point diverges. This means that $g^{-1}(\mathbf{x},\mathbf{x})$ is not well defined, contrary to the 1D case. As a consequence, the FFT resolution of the linear operator $A$ (Eq.\ref{eq:helmholtz3d} ) equaled to a unit force shows a divergent behavior near the application of the force for increasing number of voxels, Fig. \ref{fig:green3D}(left). This singularity of the Green's function in higher dimensions is a common feature with other linear operators as the Laplace equation.

Nevertheless, this singularity is only relevant for point forces and the 2D/3D algorithm proposed is well posed, because for higher dimensions force densities are applied on a  manifold $\Gamma$ with dimension greater than 0, i.e. a curve in 2D or a surface in 3D. In these cases, the displacement field resolved using FFT on $\Gamma$ result of applying a force density on the same manifold does not diverge with the discretization. To show this non-singular response, the equation \ref{eq:helmholtz3d} for a force concentrated on a plane is solved using FFT and different grid sizes. In Fig. \ref{fig:green3D}(right) the solution $u(\mathbf{x}$ in a line perpendicular to the plane $\Gamma$ is represented as function of the distance to the plane. It can be observed that the displacement is not singular in the loading plane, and the solutions converge when refining the grid. This behavior allows to invert the relation between the applied force density field on $\Gamma$ with the displacement on the manifold, represented with a deconvolution in $\Gamma$ in the continuum case (Eq. \ref{eq:deconvolution}) and with a matrix inversion on the discrete version (Eq. \ref{eq:L-determination_ND})

\section{Analytical solution of 1D elastic wave propagation}\label{anex:1Dsol}
In 1D, the solution to the wave equation Eq. (\ref{eq-wave_eq}) can be written in a general form according to Eq. (\ref{eq-soL_wave_eq}) (d'Alembert solution) \cite{Meyers},  

\begin{equation}
\label{eq-wave_eq}
\dfrac{\partial^2 u}{\partial t^2} = c_0^2 \dfrac{\partial^2 u}{\partial x^2}
\end{equation}

\begin{equation}
\label{eq-soL_wave_eq}
u(x,t) = F (x - c_0 t) + G (x + c_0 t)
\end{equation}

\noindent where $F$ and $G$ are two functions that describe the shape of the pulses propagating, along the $x$ axis, in the positive and negative directions at a velocity $c_0 = \sqrt{E/\rho}$.

In the layered case, when the propagating wave encounters a medium with different impedance $Z = \rho c_{0}$, it is reflected and refracted (transmitted) at the boundary. In this situation, the propagation of the refracted and reflected waves is also governed by Eq. (\ref{eq-soL_wave_eq}) but with different amplitude of the signals (F \& G). Considering the equilibrium and continuity conditions at the boundary Eq. \eqref{eq:wave_boundary}, the amplitude of the transmitted $A_T$ and reflected $A_R$ waves can be calculated according to Eq. \eqref{eq:wave_transmission_reflection}

\begin{equation}
\label{eq:wave_boundary}
\sigma_{I} = \sigma_{R} + \sigma_{T} \ \ ; \ \ \dot{u}_{I} = \dot{u}_{R} + \dot{u}_{T}
\end{equation}

\begin{equation}
\label{eq:wave_transmission_reflection}
A_T = \dfrac{2}{1+(Z_2/Z_1)} A_I \ \ ; \ \ A_R = \dfrac{1-(Z_2/Z_1)}{1+(Z_2/Z_1)} A_I
\end{equation}

\noindent with $A_I$ being the amplitude of the incident wave in the first material \cite{Meyers}

\bibliographystyle{unsrt}

\end{document}